\documentclass[11pt,a4paper]{amsart}
\usepackage[margin=2.5cm,marginpar=2cm]{geometry}
\usepackage{mathtools}
\usepackage{graphicx}
\usepackage[hidelinks]{hyperref}
\usepackage[nameinlink]{cleveref}

\usepackage{multirow}
\usepackage[normalem]{ulem}

\usepackage[abbrev]{amsrefs}
\usepackage{amsthm}
\usepackage{graphicx}
\usepackage{braket}
\usepackage{enumitem}
\usepackage{rotating}

\usepackage[dvipdfx,cmyk,table]{xcolor}
\usepackage{tikz}
\usetikzlibrary{positioning}
\usetikzlibrary{matrix,arrows,decorations.pathreplacing,decorations.pathmorphing,cd,calc}
\usepackage{comment}

\usepackage{amssymb}
\usepackage{amsmath}
\usepackage{mathrsfs}

\usepackage{bbm}

\usepackage[abbrev]{amsrefs}

\definecolor{refkey}{gray}{.3}
\definecolor{labelkey}{gray}{.3}

\newcommand{\Res}{\mathrm{Res}}

\def\bfone{\mathbf{1}}
\def\bbone{\mathbbm{1}}

\newcommand{\lb}{\left(} \newcommand{\rb}{\right)}

\def\Stab{{\rm Stab}}

\def\id{{\rm id}}

\def\sgn{{\rm sgn}}

\makeatletter
\def\rinn#1#2{\langle 
\def\ta{#1}\def\tb{#2}
\ifx\ta\@empty{\;} \else {\ta}\fi ,
\ifx\tb\@empty{\;} \else {\tb}\fi
\rangle} 

\def\inn#1#2{\left\langle 
\def\ta{#1}\def\tb{#2}
\ifx\ta\@empty{\;} \else {\ta}\fi ,
\ifx\tb\@empty{\;} \else {\tb}\fi
\right\rangle}

\makeatother

 

\def\innv#1#2{\inn{#1}{#2}_{V}}
\def\innvl#1#2{\inn{#1}{#2}_{V_l}}

\def\innvp#1#2{\inn{#1}{#2}_{V'}}

\def\Rep{{\rm Rep}}

\def\det{{\rm det}}

\usepackage{xparse}
\def\usecsname#1{\csname #1\endcsname}
\def\useLetter#1{#1}
\def\usedbletter#1{#1#1}

\ExplSyntaxOn 

\def\mydefbar#1#2#3{\expandafter\def\csname bar#3{#1}\endcsname{\bar{#2{#1}}}}
\def\mydefhat#1#2#3{\expandafter\def\csname hat#3{#1}\endcsname{\hat{#2{#1}}}}
\def\mydefwh#1#2#3{\expandafter\def\csname wh#3{#1}\endcsname{\widehat{#2{#1}}}}
\def\mydeft#1#2#3{\expandafter\def\csname t#3{#1}\endcsname{\tilde{#2{#1}}}}
\def\mydefu#1#2#3{\expandafter\def\csname u#3{#1}\endcsname{\underline{#2{#1}}}}
\def\mydefr#1#2#3{\expandafter\def\csname r#3{#1}\endcsname{\mathrm{#2{#1}}}}
\def\mydefb#1#2#3{\expandafter\def\csname b#3{#1}\endcsname{\mathbb{#2{#1}}}}
\def\mydefwt#1#2#3{\expandafter\def\csname wt#3{#1}\endcsname{\widetilde{#2{#1}}}}
\def\mydeff#1#2#3{\expandafter\def\csname f#3{#1}\endcsname{\mathfrak{#2{#1}}}}
\def\mydefbf#1#2#3{\expandafter\def\csname bf#3{#1}\endcsname{\mathbf{#2{#1}}}}
\def\mydefc#1#2#3{\expandafter\def\csname c#3{#1}\endcsname{\mathcal{#2{#1}}}}
\def\mydefsf#1#2#3{\expandafter\def\csname sf#3{#1}\endcsname{\mathsf{#2{#1}}}}
\def\mydefs#1#2#3{\expandafter\def\csname s#3{#1}\endcsname{\mathscr{#2{#1}}}}
\def\mydefcks#1#2#3{\expandafter\def\csname cks#3{#1}\endcsname{{\check{
\csname s#2{#1}\endcsname}}}}
\def\mydefckc#1#2#3{\expandafter\def\csname ckc#3{#1}\endcsname{{\check{
\csname c#2{#1}\endcsname}}}}
\def\mydefck#1#2#3{\expandafter\def\csname ck#3{#1}\endcsname{{\check{#2{#1}}}}}

\NewDocumentCommand{\doGreek}{m}
{
\clist_map_inline:nn {alpha,beta,gamma,Gamma,delta,Delta,epsilon,varepsilon,zeta,eta,theta,vartheta,Theta,iota,kappa,lambda,Lambda,mu,nu,xi,Xi,pi,Pi,rho,sigma,varsigma,Sigma,tau,upsilon,Upsilon,phi,varphi,Phi,chi,psi,Psi,omega,Omega,tG} {#1{##1}{\usecsname}{\useLetter}} 
}

\NewDocumentCommand{\doAtZ}{m}
{
\clist_map_inline:nn {A,B,C,D,E,F,G,H,I,J,K,L,M,N,O,P,Q,R,S,T,U,V,W,X,Y,Z} {#1{##1}{\useLetter}{\useLetter}} 
}
\NewDocumentCommand{\doatz}{m}
{
\clist_map_inline:nn {a,b,c,d,e,f,g,h,i,j,k,l,m,n,o,p,q,r,s,t,u,v,w,x,y,z} {#1{##1}{\useLetter}{\usedbletter}} 
}

\NewDocumentCommand{\doallAtZ}{}
{
\clist_map_inline:nn {mydefsf,mydeft,mydefu,mydefwh,mydefhat,mydefr,mydefwt,mydeff,mydefb,mydefbf,mydefc,mydefs,mydefck,mydefcks,mydefckc,mydefbar} {\doAtZ{\csname ##1\endcsname}}
}
\NewDocumentCommand{\doallatz}{}
{
\clist_map_inline:nn {mydefsf,mydeft,mydefu,mydefwh,mydefhat,mydefr,mydefwt,mydeff,mydefb,mydefbf,mydefc,mydefs,mydefck,mydefbar} {\doatz{\csname ##1\endcsname}}
}

\NewDocumentCommand{\doallGreek}{}
{
\clist_map_inline:nn {mydefck,mydefwt,mydeft,mydefwh,mydefbar,mydefu} {\doGreek{\csname ##1\endcsname}}
}

\NewDocumentCommand{\doGroups}{m}
{
\clist_map_inline:nn {GL,Sp,rO,rU,fgl,fsp,foo,fuu,fkk,fuu,ufkk,uK} {#1{##1}{\usecsname}{\useLetter}} 
}
\NewDocumentCommand{\doallGroups}{}
{
\clist_map_inline:nn {mydeft,mydefu,mydefwh,mydefhat,mydefwt,mydefck,mydefbar} {\doGroups{\csname ##1\endcsname}}
}

\doallAtZ
\doallatz
\doallGreek
\doallGroups
\ExplSyntaxOff

\def\GL{\mathrm{GL}}

\def\Isom{{\rm Isom}}

\def\Im{{\rm Im\,}}

\def\Ind{{\rm Ind}}

\DeclareMathOperator{\Ad}{Ad}

\DeclareMathOperator{\Hom}{Hom}
\DeclareMathOperator{\End}{End}
\DeclareMathOperator{\Irr}{Irr}
\DeclareMathOperator{\Herm}{Herm}

\DeclareMathOperator{\Norm}{Norm}
\DeclareMathOperator{\Gal}{Gal}

\long\def\qcl#1{{{\color{green}QCL: #1}}}
\long\def\zjl#1{{{\color{blue}ZJL: #1}}}

\long\def\delete#1{}

\newcommand{\trivial}[2][]{\ignorespaces\if\relax\detokenize{#1}\relax{\color{red} \vspace{0em} {[} #2 {]}}\else\ifx#1h\ifcsname showtrivial\endcsname{\color{orange}\vspace{0em}{[}#2{]}}\fi\else{\red Wrong argument!}\fi\fi\ignorespaces}

\def\Hom{{\mathrm{Hom}}}

\def\Stab{{\mathrm{Stab}}}

\def\cW{{\mathcal{W}}}

\def\half{{\frac{1}{2}}}

\NewDocumentCommand\cent{m o}{ 
\IfValueTF{#2}{\mathrm{Z}_{#1}({#2})
}{\mathrm{Z}({#1})}
}

\DeclareFontFamily{U} {MnSymbolC}{}
\DeclareFontShape{U}{MnSymbolC}{m}{n}{
<-6> MnSymbolC5
<6-7> MnSymbolC6
<7-8> MnSymbolC7
<8-9> MnSymbolC8
<9-10> MnSymbolC9
<10-12> MnSymbolC10
<12-> MnSymbolC12}{}
\DeclareFontShape{U}{MnSymbolC}{b}{n}{
<-6> MnSymbolC-Bold5
<6-7> MnSymbolC-Bold6
<7-8> MnSymbolC-Bold7
<8-9> MnSymbolC-Bold8
<9-10> MnSymbolC-Bold9
<10-12> MnSymbolC-Bold10
<12-> MnSymbolC-Bold12}{}

\DeclareFontFamily{U} {MnSymbolD}{}
\DeclareFontShape{U}{MnSymbolD}{m}{n}{
<-6> MnSymbolD5
<6-7> MnSymbolD6
<7-8> MnSymbolD7
<8-9> MnSymbolD8
<9-10> MnSymbolD9
<10-12> MnSymbolD10
<12-> MnSymbolD12}{}
\DeclareFontShape{U}{MnSymbolD}{b}{n}{
<-6> MnSymbolD-Bold5
<6-7> MnSymbolD-Bold6
<7-8> MnSymbolD-Bold7
<8-9> MnSymbolD-Bold8
<9-10> MnSymbolD-Bold9
<10-12> MnSymbolD-Bold10
<12-> MnSymbolD-Bold12}{}

\DeclareSymbolFont{MnSyC} {U} {MnSymbolC}{m}{n}
\DeclareSymbolFont{MnSyD} {U} {MnSymbolD}{m}{n}

\DeclareMathSymbol{\medstar}{\mathord}{MnSyC}{130}
\DeclareMathSymbol{\boxslash}{\mathord}{MnSyC}{114}
\DeclareMathSymbol{\boxbackslash}{\mathord}{MnSyC}{115}
\DeclareMathSymbol{\smblksquare}{\mathord}{MnSyC}{105}
\DeclareMathSymbol{\nequiv}{\mathord}{MnSyD}{121}

\DeclareMathSymbol{\otimes}{\mathrel}{MnSyC}{97}
\DeclareMathSymbol{\boxtimes}{\mathrel}{MnSyC}{117}

\def\triangle{{\mathrm{\Delta}}}

\def\dd{{{d}}}

\def\mydefdb#1#2#3{\expandafter\def\csname db#3{#1}\endcsname{\breve{#2{#1}}}}
\doGreek{\mydefdb}
\doAtZ{\mydefdb}
\doatz{\mydefdb}

\def\TTidx#1#2{\,{}^{#1}\hspace{-0.1em}#2} 

\def\any{\smblksquare}
\def\mydefTT#1#2#3{
\expandafter\def\csname ii#3{#1}\endcsname{\TTidx{i}{#2{#1}}}
\expandafter\def\csname jj#3{#1}\endcsname{\TTidx{j}{#2{#1}}}
\expandafter\def\csname zz#3{#1}\endcsname{\TTidx{0}{#2{#1}}}
\expandafter\def\csname ll#3{#1}\endcsname{\TTidx{l}{#2{#1}}}
\expandafter\def\csname aa#3{#1}\endcsname{\TTidx{a}{#2{#1}}}
\expandafter\def\csname bb#3{#1}\endcsname{\TTidx{b}{#2{#1}}}
\expandafter\def\csname oo#3{#1}\endcsname{\TTidx{1}{#2{#1}}}
\expandafter\def\csname ss#3{#1}\endcsname{\TTidx{\boxslash}{#2{#1}}}
\expandafter\def\csname dg#3{#1}\endcsname{\TTidx{\boxbackslash}{#2{#1}}}
\expandafter\def\csname any#3{#1}\endcsname{\TTidx{\any}{#2{#1}}}
}
\def\usecsname#1{\csname #1\endcsname}
\def\useLetter#1{#1}
\def\usedbletter#1{#1#1}

\doGreek{\mydefTT}
\doAtZ{\mydefTT}
\doatz{\mydefTT}

\mydefTT{cTp}{\usecsname}{\useLetter}
\mydefTT{sB}{\usecsname}{\useLetter}
\mydefTT{sL}{\usecsname}{\useLetter}
\mydefTT{tD}{\usecsname}{\useLetter}
\mydefTT{tSigma}{\usecsname}{\useLetter}
\mydefTT{Sigmap}{\usecsname}{\useLetter}
\mydefTT{ckGamma}{\usecsname}{\useLetter}
\mydefTT{Gammap}{\usecsname}{\useLetter}
\mydefTT{OmegaK}{\usecsname}{\useLetter}
\mydefTT{etanD}{\usecsname}{\useLetter}
\mydefTT{teta}{\usecsname}{\useLetter}
\mydefTT{dbkappa}{\usecsname}{\useLetter}
\mydefTT{bomega}{\usecsname}{\useLetter}
\mydefTT{biota}{\usecsname}{\useLetter}
\mydefTT{nrho}{\usecsname}{\useLetter}
\mydefTT{nnrho}{\usecsname}{\useLetter}
\mydefTT{dbrho}{\usecsname}{\useLetter}
\mydefTT{bfbb}{\usecsname}{\useLetter}
\mydefTT{sL}{\usecsname}{\useLetter}
\mydefTT{sLp}{\usecsname}{\useLetter}
\mydefTT{lD}{\usecsname}{\useLetter}
\mydefTT{nD}{\usecsname}{\useLetter}
\mydefTT{tnD}{\usecsname}{\useLetter}
\mydefTT{Vp}{\usecsname}{\useLetter}
\mydefTT{bfW}{\usecsname}{\useLetter}
\mydefTT{dbK}{\usecsname}{\useLetter}
\mydefTT{dbeta}{\usecsname}{\useLetter}
\mydefTT{dbG}{\usecsname}{\useLetter}
\mydefTT{dbJ}{\usecsname}{\useLetter}
\mydefTT{SK}{\usecsname}{\useLetter}
\mydefTT{ckG}{\usecsname}{\useLetter}
\mydefTT{End}{\usecsname}{\useLetter}
\mydefTT{fgg}{\usecsname}{\useLetter}
\mydefTT{fggp}{\usecsname}{\useLetter}
\mydefTT{sH}{\usecsname}{\useLetter}

\def\II{\mathcal{I}}

\def\cf{{\it cf.\ }}

\newcommand{\BTB}[2][]{\if\relax\detokenize{#1}\relax
\cB(#2)
\else 
\cB(#2,#1)
\fi
}
\newcommand{\rBTB}[2][]{\if\relax\detokenize{#1}\relax
\cB_{\rm red}(#2)
\else 
\cB_{\rm red}(#2,#1)
\fi
}

\newcommand{\remove}[1]{\relax}

\newcommand{\tr}{\mathrm{tr}}
\newcommand{\diag}{\mathrm{diag}}
\newcommand{\red}{\color{red}}

\crefformat{equation}{(#2#1#3)}

\crefformat{enumi}{(#2#1#3)}
\crefformat{enumi}{(#2#1#3)}
\Crefformat{enumi}{part (#2#1#3)}
\crefmultiformat{enumi}{(#2#1#3)}{ and~(#2#1#3)}{, (#2#1#3)}{ and~(#2#1#3)}
\Crefmultiformat{enumi}{part~(#2#1#3)}{ and~(#2#1#3)}{, (#2#1#3)}{ and~(#2#1#3)}
\crefrangeformat{enumi}{(#3#1#4) to~(#5#2#6)}
\Crefrangeformat{enumi}{part~(#3#1#4) to~(#5#2#6)}
\crefformat{rmk}{remark~#2#1#3}

\crefformat{thmM}{main theorem}
\Crefformat{thmM}{Main Theorem}
\newtheorem{thm}{Theorem}[section]
\newtheorem{lemma}[thm]{Lemma}

\newtheorem{remark}[thm]{Remark}

\newtheorem*{lemma*}{Lemma}
\newtheorem{prop}[thm]{Proposition}
\newtheorem{cor}[thm]{Corollary}

\newtheorem*{claim*}{Claim}

\theoremstyle{definition}

\newtheorem*{remark*}{Remark}

\crefformat{del}{{\mbox{\red Deleted reference!}}}
\Crefformat{del}{{\mbox{\red Deleted reference!}}}

\newlist{enumR}{enumerate}{1} 
\setlist[enumR]{wide,label=\arabic*.}

\newlist{enumC}{enumerate}{3} 
\setlist[enumC]{label=(\alph*)}
\newlist{enumP}{enumerate}{3} 
\setlist[enumP]{label=(\roman*)}
\setlist[enumP,1]{label=(\roman*)}
\setlist[enumP,2]{label=(\alph*)}
\setlist[enumP,3]{label=(\arabic*)}
\newlist{enumPF}{enumerate}{3}
\setlist[enumPF]{label=(\roman*),wide}
\setlist[enumPF,1]{label=(\roman*),wide}
\setlist[enumPF,2]{label=(\alph*)}
\setlist[enumPF,3]{label=\arabic*).}
\newlist{enumS}{enumerate}{3} 
\setlist[enumS]{label=\roman*)}
\setlist[enumS,1]{label=\roman*)}
\setlist[enumS,2]{label=\alph*)}
\setlist[enumS,3]{label=\arabic*.}
\newlist{enumIL}{enumerate*}{1} 
\setlist*[enumIL]{label=\roman*)}
\newlist{enumI}{itemize}{1} 
\setlist*[enumI]{label=\textbullet}

\newlist{enumST}{enumerate}{2} 
\setlist*[enumST,1]{wide,font=\bfseries,label=Step~\arabic*,leftmargin=2em}
\setlist*[enumST,2]{wide,font=\bfseries, label*=.\arabic*,leftmargin=2em}

\def\AA{\sfA}

\def\boxtimes{\otimes}

\author{Jia-jun Ma}
\address{
School of Mathematical Sciences, Xiamen University,
Xiamen China
}

\email{hoxide@xmu.edu.cn}

\author{Congling Qiu}
\address{
Department of Mathematics, Yale University, New Haven, USA 
}
\email{congling.qiu@yale.com}

\author{Jialiang Zou}
\address{
Department of Mathematical Sciences, National University of Singapore, Singapore 
}
\email{zoujialiang@u.nus.edu}


\numberwithin{equation}{section}


\title[]{generic Hecke algebra and theta correspondence over finite fields}

\providecommand{\Span}{\mathrm{Span}}

\def\barcAkk{\overline{\cA}_{k}}

\def\dww{\dot{w}}

\def\SSS{\subsection{}}

\def\tcV{\widetilde{\cV}}

\def\Akk{A_{k}}
\def\cAkk{\cA_{k}}
\def\Ak'{A'_{k-1}}
\def\cAk'{\cA'_{k-1}}
\def\Vpdag{V'^{\dag}}

\def\cVpdV{\cV'^{\diamond}_V}
\def\wtVpdag{\wtV'^{\dag}}
\def\Vpan{V'_{\mathrm an}}
    
\def\I{\triangle_l}

\newcommand{\com}{\,}

\makeatletter
\newsavebox{\@brx}
\newcommand{\llangle}[1][]{\savebox{\@brx}{\(\m@th{#1\langle}\)}%
  \mathopen{\copy\@brx\kern-0.5\wd\@brx\usebox{\@brx}}}
\newcommand{\rrangle}[1][]{\savebox{\@brx}{\(\m@th{#1\rangle}\)}%
  \mathclose{\copy\@brx\kern-0.5\wd\@brx\usebox{\@brx}}}

\def\iinn#1#2{\llangle 
\def\ta{#1}\def\tb{#2}
\ifx\ta\@empty{\;} \else {\ta}\fi ,
\ifx\tb\@empty{\;} \else {\tb}\fi
\rrangle} 

\makeatother

\begin{document}
\begin{abstract}
We study the Hecke algebra modules arising from theta correspondence between certain Harish-Chandra series for type I dual pairs over finite fields. For the product of the pair of Hecke algebras under consideration, we show that there is a  generic Hecke algebra module whose specializations at prime powers give the Hecke algebra modules and whose specialization at $1$ can be explicitly described. As an application, we prove the conservation relation on  the first occurrence indices for all irreducible representations. As another application, we generalize the results of Aubert-Michel-Rouquier and Pan on theta correspondence between the Harish-Chandra series.


\end{abstract}

\def\sfPp{\sfP_{V'^+}}
\def\xiVpp{\xi_{V'^+}}
\def\sfbb{\bW}
\def\disc{\mathrm{disc}}

\def\lVp{l'}
\def\lWp{l}
\def\ZV{Z'}
\def\ZW{Z}
\def\PV{P'}
\def\PW{P}
\def\LV{L'}
\def\LW{L}
\def\UV{U'}
\def\UW{U}

\def\cksigmaP{\cksigma_{P}}
\def\cksigmaPp{\cksigma'_{P'}}
\def\sigmaP{\sigma_{l}}
\def\sigmaL{\sigma_{L}}
\def\tsigmaL{\widetilde{\sigma}_{L}}
\def\sigmaPp{\sigma'_{l'}}

\def\tsgn{\widetilde{\sgn}}
\def\sigmaPW{\sigma_{\PW}}
\def\sigmaPV{\sigma'_{\PV}}
\def\eV{\varepsilon_V}

\def\Mp{\mathrm{Mp}}

\def\deltaan{\delta}

\def\Weil#1{{{\omega}_{\ensuremath{#1}}}}

\maketitle
\setcounter{tocdepth}{1}
\tableofcontents 

\section{Introduction}\label{introduction}
\SSS
Let $F$ be a finite field of characteristic $p$ and let $q$ be the cardinality of $F$. 
Let $\tau$ be an element in $\Gal(F/\bF_p)$ such that $\tau^2$ is the identity and $F_0$ be the fixed field of  $\tau$.
We fix a non-trivial additive character $\psi\colon F_0\rightarrow \bC^\times$ of $F_0$ throughout the paper. 
When $p$ is odd, let $ \xi$ denote the unique quadratic character of $F^\times$, i.e.,  
\begin{equation}\label{xi}
   \xi\colon F^\times \rightarrow \set{\pm 1},\quad a\mapsto a^{\half(q-1)}. 
\end{equation}

\SSS\label{dual pair cases}
For  $\epsilon\in \set{\pm 1}$, an $\epsilon$-Hermitian space $V$ (over $F$)
is an $F$-vector space equipped with a non-degenerate form $\innv{}{}:V\times V\rightarrow F$ such that 
\[
\innv{av_1+v_2}{v_3} = a \innv{v_1}{v_3} + \epsilon\innv{v_3}{v_2}^\tau, 
\quad \forall a\in F \text{ and } v_1,v_2, v_3\in V. 
\]
Let $n=\dim V$ and  let
\[
\rU(V):=\set{g\in \End_F(V)|\innv{g v_1}{g v_2}=\innv{v_1}{v_2},\ \forall v_1,v_2\in V}
\]
be the isometry group of $\inn{}{}_V$. If $F=F_0$ and $\epsilon=1$, we take a basis $\{e_1,\cdots, e_n\}$ of $V$ and define the discriminant of $V$ by 
\begin{equation}\label{discriminant}
    \disc(V):=(-1)^{\frac{n(n-1)}{2}}\det\left((\langle e_i, e_j \rangle_V)_{i,j}\right)\in F^{\times} /F^{\times 2}.
\end{equation}
Note that $\disc(V)$ is independent of the choice of the basis $\{e_1,\cdots, e_n\}$ (cf. \cite[Chapter 2, \S 2.1]{SW}).

\SSS\label{different cases}
  A dual pair $(V,V')$ consists of an $\epsilon$-Hermitian space $V$
and an $\epsilon'$-Hermitian space  $V'$ 
where $\epsilon, \epsilon'\in \set{\pm 1}$ such that $\epsilon\epsilon' = -1$. Given such a dual pair $(V,V')$, $V\otimes_F V'$ is naturally a $(-1)$-Hermitian space with
\[
\inn{v_1\otimes v'_1}{v_2\otimes v'_2}_{V\otimes_F V'} := \innv{v_1}{v_2}\innvp{v'_1}{v'_2},
\]
and there is a natural map 
\begin{equation}\label{eq:embedding}
\rU(V)\times \rU(V')\longrightarrow \rU(V\otimes_F V'). 
\end{equation}
\trivial[h]{
We view $V$ as a right $F$-vector space by $v\, a := \tau(a) v$ for all $a\in F$ and $v\in V$. 
Then $V\otimes_F V'$ is a $F$-vector space with scalar multiplication 
$a (v_1\otimes v'_1) := v_1\otimes (a v'_1)$ and $(-1)$-Hermitian form 
\[
\inn{v_1\otimes v'_1}{v_2\otimes v'_2}_{V\otimes_F V'} := \innv{v_1}{v_2}^\tau\innvp{v'_1}{v'_2},
\]
for all $a\in F$, $v_1,v_2\in V$ and $v'_1,v'_2\in V'$. 
}
According to different choices of $F$ and $\epsilon$, we have the following five cases: 
  \begin{enumerate}[label=(\Alph*),wide=0pt]
\item[($A$)] $F\neq F_0$. In this case, $(\rU(V), \rU(V'))$ is a unitary dual pair; 
\item[($B$)] $F=F_0, \epsilon=1$ and $\dim V$ is odd. In this case, $(\rU(V), \rU(V'))$ is an odd orthogonal-symplectic dual pair;
\item[($C$)] $F=F_0, \epsilon=-1$ and $\dim V'$ is even. In this case, $(\rU(V), \rU(V'))$ is a symplectic-even orthogonal dual pair; 
\item[($\widetilde{C}$)] $F=F_0, \epsilon=-1$ and $\dim V'$ is odd.
In this case, $(\rU(V), \rU(V'))$ is a symplectic-odd orthogonal dual pair. We view $\rU(V)$ as a 
``metaplectic group'' which is a notion borrowed from the local theta correspondence;
\item[($D$)] $F=F_0, \epsilon=1$ and $\dim V$ is even. In this case, $(\rU(V), \rU(V'))$ is an even orthogonal-symplectic dual pair. 
\end{enumerate}
In this paper, we assume that $p\neq 2$ in Case $B,C,\widetilde{C}$ and $D$, where $p$ is the characteristic of $F$.  

\delete{
Let $\xi_{V,V'}$ be a character of $\rU(V)\times \rU(V')$ defined by 
\[
\xi_{V,V'}\coloneqq
\begin{cases}
\bbone\boxtimes \bbone \quad &\text{in Case $A$}\\
\xi_V^{\frac{\dim V'}{2}}\boxtimes \mathbbm 1\quad & \text{in Case $B$ and $D$}\\
\mathbbm 1\boxtimes \xi_{V'}^{\frac{\dim V}{2}}\quad & \text{in Case $C$ and $\widetilde{C}$},
\end{cases}
\]
where $\xi_V$ is the quadartic character of $\rU(V)$ defined by 
\[
\xi_V(g)\coloneqq \det_V(g)^{\half (q-1)}
\]
and likewise for $\xi_{V'}$. 
}

\subsection{}
Let $\omega_{V\otimes_F V'}$ be the Weil representation of $\rU(V\otimes_F V')$ attached to the additive character $\psi$ constructed in G\'erardin~\cite[Theorem 3.3]{Ger} for Case $(A)$ and in G\'erardin~\cite[Theorem 2.4]{Ger} for other cases ($p\neq 2$ here).
The (modified) Weil representation $\omega_{V,V'}$ of $\rU(V)\times \rU(V')$ is defined by 
\begin{equation}\label{eq:Modified Weil representation}
\begin{split}
 &    \omega_{V,V'}\coloneqq   
   \begin{cases}
   \omega_{V\otimes_F V'}
   & \text{ in case $(A)$,}\\
   \left((\xi\circ \det_{\rU(V)})^{\half\dim_F V'}\boxtimes \bbone_{ \rU(V')} \right)\otimes
   \omega_{V\otimes_F V'}
   & \text{ in case $(B)$ $(D)$,}\\
   \left(\bbone _{\rU(V)}\boxtimes (\xi\circ \det_{ \rU(V')})^{\half\dim_F V} \right)\otimes
   \omega_{V\otimes_F V'}
   & \text{ in case $(C)$ $(\wtC)$,}\\
   \end{cases}
   \end{split}
\end{equation}
where $\omega_{V\otimes_F V'}$ is viewed as a $\rU(V)\times \rU(V')$-representation
via \eqref{eq:embedding},
 and $\det_{\rU(V)}$ (resp. $\bbone _{\rU(V)}$)  denotes the determinant (resp. trivial) character of  $\rU(V)$.

\SSS
Let $\Rep(G)$ denote the category of representations of $G$ over $\bC$ 
and $\Irr(G)$   the set of irreducible representations of $G$. 
For each dual pair  $(V,V')$,  we define a functor 
$\Theta_{V,V'}$ from $\Rep(\rU(V))$ to $\Rep(\rU(V'))$ by  
\[
\Theta_{V,V'}(\pi) \coloneqq \Hom_{\rU(V)}(\pi, \omega_{V,V'}). 
\]
When $\pi \in \Irr(\rU(V))$, the $\pi$-isotypic component of $\omega_{V,V'}$ is canonically isomorphic to  \[
\pi\boxtimes \Theta_{V,V'}(\pi).
\]
Therefore, our definition agrees with the conventional definition of the big theta lift map.

\SSS \label{definition of Witt tower}
The hyperbolic $\epsilon'$-Hermitian plane $\rH_{\epsilon'}$ is the unique (up to isomorphism) $2$-dimensional $\epsilon'$-Hermitian space 
containing a non-zero isotropic vector. 
 A set  $\cV'$ of  $\epsilon'$-Hermitian spaces is called 
a Witt tower if  there is an anisotropic 
$\epsilon'$-Hermitian space
$V'_{\rm an}$  and   $ \cV' $  consists of all $\epsilon'$-Hermitian spaces $\Vpdag$ such that
\[
\Vpdag \cong \Vpan\oplus \underbrace{\rH_{\epsilon'}\oplus\cdots \oplus\rH_{\epsilon'}}_{r\text{-terms}}.
\]
The integer $r$ above is called the split rank of $\Vpdag$ and the parity of $\dim \Vpan$ is called the parity of $\cV'$. If $F=F_0$ and $\epsilon'=1$, we have 
\[
\disc(\Vpdag)=\disc(\Vpan) \quad \text{for all } \Vpdag\in \cV'.
\]
Therefore, we  define $\disc(\cV'):=\disc(V'_{\rm an})$ and call it the discriminant of $\cV'$ in this case.

For a dual pair $(V,V')$, let $\cV'$ be the Witt tower containing $V'$ and 
we define the generalized Witt tower 
$\cV'_{V}$ to   be the collection 
\[
\cV'_{V}=\Set{\left(\Vpdag,\Weil{V,\Vpdag}\right)| \Vpdag\in \cV'},
\]
where $\Weil{V,\Vpdag}$ is the modified Weil representation defined in \Cref{eq:Modified Weil representation} (\cf \cite{SZ}*{\S~3}).

We define the companion generalized Witt tower $\tcV'_V$  as follows: 
\begin{itemize}
    \item In Case $(A)$, there exists a unique Witt tower $\tcV'$ such that parity of $\tcV'$ and $\cV'$ are different.  
    We define 
    \[
    \tcV'_V \coloneqq \Set{\left(\wtVpdag,\Weil{V,\wtVpdag}\right)| \wtVpdag \in \tcV'}.
    \]
    \item In Case $(C)$ and $(\wtC)$, there exists a unique Witt tower $\tcV'$ such
    that $\tcV'$ and $\cV'$ have the same parity and different discriminants. 
    We define
    \[
    \tcV'_V \coloneqq \Set{\left(\wtVpdag,\Weil{V,\wtVpdag} \right)
    | \wtVpdag \in \tcV'}.
    \]
 \item In Case $(B)$ and $(D)$, there is only one Witt tower of symplectic spaces.
Let $\tcV' \coloneqq \cV'$ and 
 \[
    \tcV'_V \coloneqq \Set{\left(\wtVpdag,\omega_{V,\wtVpdag}\otimes 
    (\det_{\rU(V)}\boxtimes \bbone_{\rU(\wtVpdag)})\right)|\wtVpdag\in \tcV'}.
 \]
\end{itemize}


\SSS
For a generalized Witt tower $\cVpdV\in  \set{\cV'_V,\tcV'_V}$ of $\epsilon'$-Hermitian spaces, 
we define a function $n_{\cVpdV}$ from $\Irr(\rU(V))$ to the set $\bN$ of non-negative integers by sending $\pi \in \Irr(\rU(V))$ to 
\begin{equation}\label{First occurence}
  n_{\cVpdV}(\pi)\coloneqq \min\Set{\dim \Vpdag|\left(\Vpdag,\omega \right)\in \cVpdV
  \text{ and } \Hom_{\rU(V)}(\pi, \omega)\neq 0 }.
\end{equation}
The number $n_{\cVpdV}(\pi)$ 
is called the first occurrence index of $\pi$ with respect to $\cVpdV$, which is  
well-defined by the non-vanishing of theta lift when the split rank of $\Vpdag\in \cVpdV$ is greater than $\dim V$, see \cite{MVW}*{3.IV}.

In the rest   of the introduction, we will retain the notations $V$, $\cV'_V$ etc. defined above.  
All definitions and notations extend naturally to $V'$ by switching the roles of $V$ and $V'$. 

\SSS\label{theta cuspidal}
Suppose that $V$ has split rank $r$. 
We fix a chain of $\epsilon$-Hermitian spaces $\Set{V_i| -r\leq i \in \bZ}$ 
in the Witt tower 
containing $V$ such that  $V_{i}\subset V_{i+1}$ and 
$\dim V_{i} = \dim V + 2i$ for each integer $i\geq -r$. 
Note that $V_0 =V$. For each integer $0 \leq k \leq r$, let
$Q_k=M_kN_k$ be a parabolic subgroup of $\rU(V)$ with Levi component  
\[
M_k\cong \GL_k(F) \times \rU(V_{-k})
\]  
 unipotent radical $N_k$.   
Let $J_{\rU(V)}^{M_k}$ be the Jacquet functor from $\rU(V)$ to $M_k$ 
by taking $N_k$-coinvariants. 
Set 
\begin{equation}\label{chi}
\chi_{V,V'}\coloneqq
    \begin{cases}
    \xi & \text{in Case $(\wtC)$},\\ 
    \bbone & \text{otherwise}. 
    \end{cases}
\end{equation}
We abbreviate $\chi_{V,V'}$ to be $\chi$ if there is no confusion caused.  
For each $\pi\in \Irr(\rU(V))$, define 
\[
c(\pi) := \max\Set{k|\Hom_{\GL_k}(J_{\rU(V)}^{M_k}(\pi), \chi\circ \det))\neq 0 }.
\]
\trivial[h]{
$\sigma'' \leftrightarrow \alpha''\times \beta''$
$\Ind_{Q_r}^G \chi\circ {\det}\otimes \sigma'' \leftrightarrow \Ind_{S_r\times W_{l-r}}^{W_l} 1\otimes (\alpha\times \beta) = 
$
}
  Clearly if $\sigma$ is a cuspidal representation of $\rU(V)$, then  $c(\sigma) = 0$.
In general, we call a representation $\sigma\in \Irr(\rU(V))$ \emph{theta cuspidal} (with respect to the cases listed in \Cref{dual pair cases}) when $c(\sigma)=0$. 

The behavior of theta cuspidal representations are similar to cuspidal representations under theta correspondence. 
Indeed,
by the argument in \cite{Ku} (see also \cite{MVW}*{3.IV}) adapted to our situation, one can see that, when $\sigma$ is theta cuspidal,   $\Theta_{V,V'}(\sigma)$ is irreducible or zero for every $V'\in \cV'_V$. 
Moreover, $\Theta_{V,V'}(\sigma)$ is theta cuspidal if and only if $\dim V'=n_{\cV'_V}(\sigma)$. 


\SSS \label{Coservation section}
We prove the following conservation relation on the first occurrence indices.
Let 
\begin{equation}\label{deltadelta'}
\deltaan=\begin{cases}
1-\epsilon \quad &\mbox{if $F=F_0$},\\
1  \quad &\mbox{if $F\neq F_0$},
\end{cases} \quad\mbox{and}  \quad
\delta'=\begin{cases}
1+\epsilon \quad &\mbox{if $F=F_0$},\\
1 \quad &\mbox{if $F\neq F_0$}.
\end{cases}
\end{equation}

\begin{thm}\label{thm:conservation}
For $\pi\in \Irr(\rU(V))$,  we have 
\begin{equation}\label{eq:conservation relation}
n_{\cV'_V}(\pi)+n_{\tcV'_V}(\pi) + c(\pi) = 
2\dim V + \deltaan .
\end{equation}
\end{thm}

We will first prove \Cref{thm:conservation} for theta-cuspidal representations (\Cref{Conservation}).
The general case will be  proved in  \Cref{thm:conservationproof}, using \Cref{thm:main} below (whose proof relies on the  theta-cuspidal case of \Cref{thm:conservation}).  

For cuspidal representations, Pan obtained \Cref{eq:conservation relation} 
by a reduction to the $p$-adic case \cite[Theorem 12.3]{Pan02J}.
Recently Pan also obtained the general case by a reduction to the unipotent case \cite{Pan22}. 

\SSS\label{genericHeckealgebra} 
Now we  introduce another important player, the generic Hecke algebra.
Fix an integer $l\geq 0$.
Let $\sfW_l$ be  the Coxeter group of type  $BC_l$ with simple reflections
\[
\I:=\set{\sfss_1, \sfss_2, \cdots, \sfss_{l-1}, \sftt_l},
\]
where   the subgroup  $\sfS_l$  of $\sfW_l$ generated by $\{\sfss_1,\cdots, \sfss_{l-1}\}$   is isomorphic to the symmetric group on $l$ letters. 
Let $l\colon \sfW_l\rightarrow \bN$ be the length function on $\sfW_l$. 

Let $\nu$ be an indeterminate and $R:= \bZ[\nu^{\half},\nu^{-\half}]$
  the ring of polynomials with coefficients in $\bZ$ generated by $\nu^{\pm\half}$. For  $\mu\in  \half \bZ$  and a  non-negative integer  $l$,
the generic Hecke alegebra $\sfH_{l,\mu}$ of $\sfW_l$ 
is the unique associative algebra over $R$ with $R$-basis $\set{\sfT_{\sfww,\mu}  | \sfww\in \sfW_l}$
(we abbreviate $\sfT_{\sfww}:=\sfT_{\sfww,\mu} $),
  such that
\begin{enumerate}[label=(\alph*),wide=0pt]
\item $(\sfT_{\sfss_i}+1)(\sfT_{\sfss_i}-\nu)=0, i=1,\cdots, l-1$,
\item $(\sfT_{\sftt_l}+1)(\sfT_{\sftt_l}-\nu^\mu)=0$,
\item $\sfT_{\sfww_1} \sfT_{\sfww_2} = \sfT_{\sfww_1\sfww_2}$ if 
$l(\sfww_1\sfww_2) = l(\sfww_1)+l(\sfww_2)$. 
\end{enumerate}

 
Let $\bC_q$ be the $R$-module whose underlying $\bZ$-module is $\bC$  
and $\nu^{\half}$ acts by $\sqrt{q}$. 
Tensoring of an $R$-algebra or $R$-module with $\bC_q$ is called the \emph{specialization} at $\nu=q$. 
We write $\sfH_{l,\mu,\nu=q}:= \sfH_{l,\mu}\otimes_R \bC_q$ for the specialization of an $R$-algebra $\sfH_{l,\mu}$
at $\nu=q$ and $\sfM_{\nu=q}:=\sfM\otimes_R \bC_q$ for the $\sfH_{l,\mu,\nu=q}$-module obtained by specialization of an $\sfH_{l,\mu}$-module $\sfM$ at $\nu=q$.  

\SSS
\label{chain} Let
$P_l=L_lU_l$  be a  parabolic subgroup of $\rU(V_l)$ with unipotent radical $\rU_l$ and Levi component  
\[
L_l\cong \underbrace{\GL_1(F)\times\cdots \times \GL_1(F)}_{l\text{-terms}}
\times \rU(V_0). 
\] 
Let $\sigma$ be an irreducible theta cuspidal representation of $\rU(V_0)$ and let 
\begin{equation}\label{sigmaldef}
\sigmaP:=\chi^{\boxtimes l}\boxtimes \sigma
\end{equation}
be an irreducible representation of $L_l$. By abuse of notation, we also use $\sigmaP$ to denote its inflation to $P_l$. 
Consider the intertwining algebra 
\begin{equation}\label{eq:cHl}
\cH_l:= \End_{\rU(V_l)}\left(\Ind_{P_l}^{\rU(V_l)}(\sigmaP^\vee)\right),
\end{equation}
where $\sigmaP^\vee$ denotes the contragredient of $\sigmaP$. 

\begin{thm}\label{thm:norm}
There is an  isomorphism  $\cH_l\cong \sfH_{l,\mu(\sigma),\nu=q}$ with 
\begin{equation}\label{eq:mu}
\mu(\sigma) := \half (n_{\cV'_V}(\sigma) - n_{\tcV'_V}(\sigma)).
\end{equation}
\end{thm}
The isomorphism in \Cref
{thm:norm} will be explicit constructed in \Cref{structure of cHl}. 
Note that we  also have $\cH_l\cong \sfH_{l,-\mu(\sigma),\nu=q}$ by switching the role of $\cV'_V$ and $\tcV'_V$ in \Cref{eq:mu}.  These two isomorphism  
are related by a natural isomorphism 
$
\kappa: \sfH_{l,-\mu(\sigma)}\cong \sfH_{l,\mu(\sigma)}
$, see \Cref{Different normalization}.


 In the case that $\sigma $ is cuspidal, by  Howlett-Lehrer \cite[Theorem 3.14]{HL}, Lusztig \cite[8.6]{L2} and Geck \cite[Corollary 2]{Geck}, it was already known  that $\cH_l= \sfH_{l,\mu,\nu=q}$ for some $\mu \in \half\bZ$.
 Thus the content of \Cref{thm:norm} is to relate  the parameter of the generic Hecke algebra with the first occurrence indexes of $\sigma$.  


\subsection{}
Let $(V,V', \sigma,\sigma')$ be a quadruple such that   
  $(V,V')$ is a dual pair and 
 $(\sigma,\sigma')\in \Irr(\rU(V))\times \Irr(\rU(V'))$ is a pair of 
 theta cuspidal representations that satisfy
 $\Theta_{V,V'}(\sigma) = \sigma'$. 
 By \Cref{thm:conservation} in the  theta-cuspidal case (\Cref {Conservation}) and \Cref{thm:norm}, we have 
\begin{equation}\label{eq:mu1}
\mu(\sigma) = \dim V' - \dim V - \half \delta
\end{equation}
and 
\[
\quad\mu(\sigma') = \dim V - \dim V' - \half \delta' = -1 -\mu(\sigma).
\]
%

 In view of \eqref{deltadelta'} and \eqref{eq:mu1}, let  \begin{equation}\label{eq:mu.range}
 \mu\in 
\begin{cases}
 \bZ+\half \quad &\mbox{in Case $(A)$,}\\
2 \bZ +1\quad &\mbox{in Case $(B)$ and $(C)$,}\\
 2\bZ \quad &\mbox{in Case $(\wtC)$ and $(D)$}.
\end{cases}
\end{equation}
We say that  $\mu $  and 
  a   quadruple $(V,V', \sigma,\sigma')$ as above are   \emph{relevant} 
if \begin{equation}\label{eq:mu11}  
\mu=\dim V' - \dim V - \half \delta.
\end{equation} 
We prove that
every $\mu$ in \eqref{eq:mu.range} is relevant  to some quadruple as above (\Cref{abundance}).
  
\SSS\label{mainsec}
Below in this introduction, we use Hecke algebras to analyze theta lifts. 
Fix $\mu$ in \eqref{eq:mu.range} and integers $l,l'\geq 0$.
Let $(V,V', \sigma,\sigma')$ be a quadruple relevant to $\mu$.
By \cite[\S 3, IV.4]{MVW}, we know that the irreducible components of the theta lifts of elements in the Harish-Chandra series
\[
\cE(\rU(V_l),\sigma):=\set{\pi \in \Irr(\rU(V_l)) |\pi \text{ occurs in } \Ind_{P_l}^{\rU(V_l)}(\sigmaP)}
\]
to $\rU(V'_{l'})$ lie in the Harish-Chandra series $\cE(\rU(V'_{l'}),\sigma')$.
From the definition \eqref{eq:cHl} and \Cref{thm:norm}, 
$\cE(\rU(V_l),\sigma)$ is  
parameterized by simple modules of  $ \cH_l$. Likewise $\cE(\rU(V'_{l'}),\sigma')$ is parameterized by simple modules of $ \cH'_{l'}$ which is defined as in \eqref{eq:cHl}. Thus the decomposition of the   multiplicity space
 \begin{equation}\label{eq:multsp}
    \cM:= \Hom_{\rU(V_l)\times \rU(V'_{l'})}\left(
    \Ind_{P_l}^{\rU(V_l)}(\sigmaP)\boxtimes \Ind_{P_l}^{\rU(V'_{l'})}(\sigmaPp),
    \omega_{V_l,V'_{l'}}\right)
\end{equation}
 as a  \emph{left} $  \cH_l\otimes \cH'_{l'}$-module will give rise to the theta correspondence between $\cE(\rU(V_l),\sigma)$ and $\cE(\rU(V'_{l'}),\sigma')$.
 We show that $\cM$ lifts to a generic Hecke algebra module.
 
 \trivial[h]{
 \begin{equation}\label{nu to 1}
 \begin{split}
  \sfH_{\nu=1}&\cong \bC[\sfW_l]\\
 \sfT_{\sfww}\otimes 1&\mapsto \sfww    \quad \sfww\in \sfW_l,
 \end{split}
 \end{equation} 
 } 
 
  Let $\sfH=\sfH_{l,\mu}$ and $\sfH'=\sfH_{l',-1-\mu}$.   Fix the isomorphisms $\cH_{l}\cong \sfH_{\nu=q}$ and $\cH'_{l'}\cong \sfH'_{\nu=q}$ from \Cref{thm:norm}.  Fix the natural isomorphism 
  $\bC[\sfW_l]\cong \sfH_{\nu=1}$ by sending $\sfww\in \sfW_l$ to $\sfT_{\sfww}\otimes_R 1$,
and likewise
$\bC[\sfW_{l'}]\cong \sfH'_{\nu=1}$.
\trivial[h]{$\sfH'_{\nu=1}\cong \bC[\sfW_{l'}]$}
Define the ``signature'' representation of $\sfH$ by 
 \begin{equation}\label{varepsilon of cHl}
   \varepsilon_l(\sfT_{\sfss_i})=\nu\mbox{ for }  i=1,\cdots, l-1,\quad 
   \mbox{and}\quad \varepsilon_l(\sfT_{\sftt_l})=-1.  
 \end{equation}
 To ease the notation, 
 we still denote the specialization of   $\varepsilon_l$  at $\nu=q$ or $\nu=1$ by the same symbol. 
 The same definition  works for $\sfH'$ and its specializations as well. 
 
  For each non-negative integer $k$, 
we have natural subgroups 
$\sfW_{l-k}\times \sfW_k \subset \sfW_{l}$ and  $\sfW_{k}\times\sfW_{l'-k}\subset \sfW_{l'}.$
Let $\triangle \sfW_{k}\subset (\sfW_{l-k}\times \sfW_k)\times(\sfW_{k}\times \sfW_{l'-k}) \subset \sfW_{l}\times  \sfW_{l'}$ be the diagonal embedding.
%
 
\begin{thm}\label{thm:main}
There is an   $\sfH\otimes \sfH'$-module $\sfM$, free over $R$  such that
\begin{enumerate}[label=(\alph*),wide=0pt]
\item as $\sfH_{\nu=q}\otimes\sfH'_{\nu=q}\cong \cH_l\otimes \cH'_{l'}$-modules, 
$\sfM_{\nu=q} \cong \cM $; 
\item
 as $\sfH_{\nu=1}\otimes \sfH'_{\nu=1} \cong\bC[\sfW_l]\otimes\bC[\sfW_{l'}]$-modules,
\[
   \sfM_{\nu=1}\cong 
\sum_{k=0}^{\min\set{l,l'}} \Ind_{\sfW_{l-k}\times \triangle \sfW_{k}\times  \sfW_{l'-k}}^{\sfW_{l}\times \sfW_{l'}}\left(
\varepsilon_{l-k}\otimes\varepsilon_k \otimes \varepsilon_{l'-k}\right).
\]
\end{enumerate}
\end{thm}
We will explicitly construct the module $\sfM$ in \Cref{thm:mainproof}. 
\begin{remark*}
    The generic Hecke algebras $\sfH$, $\sfH'$, and the module $\sfM$ only depend on the triple $(\mu, l, l')$.
\end{remark*}

\SSS \label{result interplation}
By Tits' deformation \cite[10.11]{C2}, we have a natural bijection between simple modules of $\cH_l$ and $\Irr(\sfW_l)$ under our fixed isomorphisms  $\cH_l\cong \sfH_{\nu=q}$  
and $\bC[\sfW_l]\cong \sfH_{\nu=1}$
. Thus we can identify  $\cE(\rU(V_l),\sigma)$ with $\Irr(\sfW_l)$
. Likewise, we identify  $\cE(\rU(V'_{l'}),\sigma')$ with $\Irr(\sfW_{l'})$.
We now recall some basic facts in the representation theory
 of $\sfS_l$ and $\sfW_l$.
Let $\cS_a$ be the Grothendieck group of finite dimensional representations of $\sfS_a$.
Then $\cS:=\bigoplus_{a\in \bN} \cS_a$ is a graded commutative unital ring under the multiplication

\begin{equation}\label{Symm}
\alpha\cdot \beta :=  \Ind_{\sfS_a\times \sfS_b}^{\sfS_{a+b}} (\alpha \boxtimes \beta), 
\quad \forall \alpha\in \cS_a, \beta\in \cS_b. 
\end{equation} 
Let $\bfone_l$ be the trivial representation of $\sfS_l$. Define $
 \sfX_l \colon \cS \rightarrow \cS$ by $\alpha \mapsto  \bfone_l\cdot \alpha$
and let $\sfX_l^*\colon  \cS\rightarrow \cS$ be the adjoint of $\sfX_l$, under  the  non-degenerate pairing  that assigns $\inn{\alpha}{\beta} = \delta_{\alpha,\beta}$  (Kronecker delta)
if both $\alpha$ and $\beta$ are irreducible representations.  The operators $\sfX_l$ and $\sfX^*_l$ can be computed by Pieri's rule (cf. \cite[4.44]{Fu}). 
Similarly, we define $\cW:= \bigoplus_{a\in \bN} \cW_a$ with $\cW_a$ being
the Grothendieck group of finite dimensional representations 
of $\sfW_a$. Then $\cW$ is also a commutative unital ring under 
the obvious analog of \eqref{Symm}. For $\alpha\in \Irr(\sfS_a)$, let $\wtalpha$ denote the inflation of $\alpha$ to $\sfW_a$ via the natural quotient map $\sfW_l\rightarrow \sfS_l$. 
It is a classical result that the map 
\begin{equation*}\label{Hyperocta}
\Irr(\sfS_a)\times \Irr(\sfS_b) \ni (\alpha,\beta) \mapsto 
\alpha\times \beta\coloneqq 
\Ind_{\sfW_a\times \sfW_b}^{\sfW_{a+b}} (\wtalpha \boxtimes (\wtbeta \otimes \varepsilon_b))\in \cW_{a+b}.
\end{equation*}
extends to a graded ring isomorphism 
$\cS\otimes \cS \rightarrow \cW$, see \cite{Z}*{Section~7.3}. 


\begin{thm}\label{cor:AMR} 
Suppose $\pi\in \cE(\rU(V_l),\sigma)$ and corresponds to $\alpha\times\beta\in \Irr(\sfW_l)$. Then $\Theta(\pi)$ is a multiplicity-free combination of
irreducible representations in $\cE(\rU(V'_{l'}), \sigma')$ corresponding to 
\begin{equation}\label{Thm 1.4}
   \sum_{k=0} ^{\min\set{l,l'}}  \sfX_{l-k}^* (\beta) \times \sfX_{l'-k} (\alpha)\in \cW_{l'}. 
\end{equation}

\end{thm}
  \Cref{cor:AMR} will be proved in  \Cref{proof of theorem 2} using \Cref{thm:main}. This theorem covers  the results of Aubert-Michel-Rouquier \cite[Theorem 3.10, Conjecture 3.11]{AMR} and 
Pan \cite[Theorem 3.4]{PanUni} for theta lifts 
between unipotent representations and quadratic-unipotent representations (called $\theta$-representations in \cite{LW}*{Theorem~3.3}), 
see \Cref{comparison normalization}.

\subsection*{Acknowledgment}
The first author thanks Anne-Marie Aubert, Tomasz Przebinda, and Binyong Sun for helpful discussions. 
He is supported in part by the National Natural Science Foundation of China (Grant No. 11701364 and Grant No. 11971305) and  Xiamen University Malaysia Research Fund (Grant No. XMUMRF/2022-C9/IMAT/0019). The second author is partially supported by the NSF grant DMS-2000533. The third author thanks Yifeng Liu and Zhiwei Yun for encouraging discussions. He is partially supported by a Singapore government MOE Tier 1 grant R-146-000-320-114.  

\section{Preliminary: Parabolic subgroup, Weyl group, and Weil representation} 
We explicate some parabolic subgroups of $\rU(V_l)$ and the corresponding relative Weyl groups. We also recall the mixed model of the Weil representation and   Kudla's filtration.

\subsection{Parabolic subgroup and Weyl group}\label{Parabolic subgroup}
Recall that for $l\geq 0$,  $V_l=V\oplus\rH_{\epsilon}^l$. 
Fix an isotropic basis $v_i,v_{- i}$ of the $i$-th $\rH_{\epsilon}$ such that $\innvl{v_i}{v_{-i}}=1$.  For each $0\leq k\leq l$, we define 
\[
 V_k^+=\Span\set{ v_{1},\cdots, v_k} ,\quad  V_k^-=\Span\set{  v_{-1},\cdots, v_{-k}},
\] 
\begin{equation*}\label{Vl-k+} 
  \widehat{V}_{l-k}^+ = \Span\Set{v_{k+1}, \cdots, v_{l}},\  
  \widehat{V}_{l-k}^-= \Span\set{v_{-(k+1)}, \cdots, v_{-l}}.
\end{equation*}
Note that $V_l^+=V_k^+\oplus \widehat{V}_{l-k}^+$. We also define 
\begin{equation}\label{Vl-k}
  \widehat{V}_{l-k}= \widehat{V}_{l-k}^+\oplus V \oplus \widehat{V}_{l-k}^-.   
\end{equation}
Then 
\begin{equation}\label{V_l}
V_l=V^+_k \oplus \widehat{V}_{l-k} \oplus V^-_k.    
\end{equation}
The  discussion and notations here and below apply to $V $ and $V'$ in obviously parallel ways.  

Using these  subspaces, we  describe the parabolic subgroup of $\rU(V_l)$ that appeared in \ref{chain}.
For $v\in V_l$, let $\langle v \rangle$ be its span.
We identify  
$\GL_1(\langle v_i\rangle)\subset\GL(V_l^+)\subset\rU(V_l^+\oplus V_l^-)$ and 
$\rU(V_0)$ as  subgroups of
$ \rU(V_l)$ via the above choices of bases and decompositions.   
Specify the parabolic subgroup  $P_l=L_lU_l$  of $\rU(V_l)$ in \ref{chain} to be the one stabilizing the   flag 
$
0\subset V_1^+\subset \cdots \subset V_l^+,
$ and specify its Levi subgroup $L_l$ to be the one stabilizing the   flag 
$
0\subset V_1^-\subset \cdots \subset V_l^-.
$ Let 
$$T_l= \GL_1(\langle v_1\rangle )\times...\times  \GL_1(\langle v_l\rangle )\subset\GL(V_l^+) .$$
 It is direct to check the following 
  equality
$$\Norm_{\rU(V_l)}(L_l) = \Norm_{\rU(V_l^+\oplus V_l^-)}(T_l)  \times \rU(V_0)$$
 between the normalizers.
 In particular, the natural embedding yields the following isomorphism 
\[
\Norm_{\rU(V_l^+\oplus V_l^-)}(T_l)/T_l\xrightarrow{\ \ \ \cong \ \ \ } \Norm_{\rU(V_l)}(L_l)/L_l. 
\]
\trivial[h]{
Note that $V_l^+\oplus V_l^-$ is the sum of vectors in $V_l$ for which the action of $L_l$ factor through $\GL_1(F)$. Therefore, if $g\in \Norm_{\rU(V_l)}(L_l)$, then $g$ stabilize $V_l^+\oplus V_l^-$ and $(V_l^+\oplus V_l^-)^\perp =V_0$, which implies $g\in U(V_l^+\oplus V_l^-)\times \rU(V_0)$.  \zjl{Here we use $\GL_1\not \cong \rU(V_0)$. What if $\dim V_0=1$, $F=2$ and $\rU(V_0)=\mathrm O(V_0)$.} \qcl{Maybe we can just say this is obvious}
}

Now let us consider the Weyl group and mark the simple reflections of $\sfW_l$
on the following Dynkin diagram of type $BC_l$: 
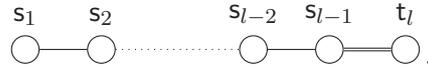
\[
\tikzstyle{vertex}=[circle, draw, minimum size=0pt]
\newcommand{\vertex}{\node[vertex]}
\begin{tikzpicture}
\vertex[label={$\sfss_1$}] (s1) at (-2,0) {};
\vertex[label={$\sfss_2$}] (s2) at (-1,0) {};
\vertex[label={$\sfss_{l-2}$}]  (s3) at (1,0) {};
\vertex[label={$\sfss_{l-1}$}]  (sl) at (2,0) {};
\vertex[label={$\sftt_{l}$}]  (tl) at (3,0) {};
\draw[] (s1) -- (s2);
\draw[dotted] (s2) -- (s3);
\draw[] (s3) -- (sl);
\draw[double] (sl) -- (tl);
\end{tikzpicture}.
\]
For $k=1,\dots,l-1$, let
\begin{equation}\label{sfttk}
\sftt_k := \sfss_k \cdots \sfss_{l-1}\, \sftt_l \, \sfss_{l-1} \cdots \sfss_k. 
\end{equation} 
Recall that $\sfS_l$ is the subgroup of $\sfW_l$ generated by $\{\sfss_1,\cdots, \sfss_{l-1}\}$ which is  isomorphic to the symmetric group on $l$ letters.  
Let $\sfD_l$ be the subgroup of $\sfW_l$ generated by $\{\sftt_1,\cdots, \sftt_l\}$ which is isomorphic to $(\bZ/2\bZ)^l$. 
We have 
\begin{equation}\label{sfW semi}
 \sfW_l=\sfD_l\rtimes \sfS_l .   
\end{equation}
We make the following choices of lifts of elements in $\sfW_l$ to $\Norm_{\rU(V_l^+\oplus V_l^-)}(T_l)$ according to Tits \cite{Ti}: 
\begin{itemize}
\item For $i=1,\cdots, l-1$, the lift of $\sfss_i$ is
\begin{equation}\label{Tits lifts si}
 s_i: \begin{cases}
v_i &\mapsto -\epsilon v_{i+1},\\
v_{i+1}&\mapsto \epsilon v_{i},\\
v_{-(i+1)}&\mapsto \epsilon v_{-i},\\
v_{-i}&\mapsto -\epsilon v_{-(i+1)},\\
v_j &\mapsto v_j\quad \mbox{if $-l\leq j\leq l$ and $j\neq\{ i,i+1,-i,-(i+1)\}$}.
\end{cases}   
\end{equation}
\item The lift of $\sftt_l$ is 
\begin{equation}\label{Tits lifts tl}
    t_l: \begin{cases}
v_l &\mapsto -v_{-l},\\
v_{-l}&\mapsto -\epsilon v_{l},\\
v_j &\mapsto v_j\quad \mbox{if $-l<j<l$}.
\end{cases}
\end{equation}

\item If $\sfww=\sfww_{1}\cdots \sfww_{j}$ is an reduced expression of $\sfww\in \sfW_l$ with $\sfww_i\in \I$, take the lift of $\sfww$ to be $w=w_{1}\cdots w_{j}$, where $w_i, 1\leq i\leq j$ is the lift of $\sfww_i$  in \Cref{Tits lifts si} and \Cref{Tits lifts tl}. 
\end{itemize}
The definition $\sfww \mapsto w$ above does not depend on the choice of reduced expression of $\sfww$  since  the lifts in \Cref{Tits lifts si} and \Cref{Tits lifts tl} satisfy the following braid relation  
\begin{equation}\label{braid equation}
\begin{split}
    s_is_j&=s_js_i,\quad \mbox{for all $j\neq i\pm 1$ and $1\leq i,j\leq l-1$},\\
    s_i s_{i+1}s_i&=s_{i+1}s_i s_{i+1}\quad \mbox{for all $1\leq i\leq l-1$},\\
    s_it_l&=t_l s_i,\quad \mbox{for $1\leq i\leq l-2$},\\
    s_{l-1}t_ls_{l-1}t_l&=t_ls_{l-1}t_ls_{l-1}.
\end{split}
\end{equation}
The lifts above induce the natural identification 
 \[\sfW_l = \Norm_{\rU(V_l^+\oplus V_l^-)}(T_l)/T_l=\Norm_{\rU(V_l)}(L_l)/L_l. \]

By \eqref{sfttk}, for $k=1,\dots,l-1$, the lift of $\sftt_k$ equals  
\[
t_k: \begin{cases}
v_k &\mapsto (-1)^{l-k+1} v_{-k},\\
v_{-k}&\mapsto \epsilon (-1)^{l-k+1} v_{k},\\
v_j &\mapsto -v_j\quad \mbox{if $j>k$ or $j<-k$},\\
v_j &\mapsto v_j\quad \mbox{if $-k<j<k$}.
\end{cases}
\]

We end this subsection by introducing some maximal parabolic subgroups. For   $1\leq k\leq l$, 
by \Cref{V_l}, we identify $\GL(V_k^+)$ and 
$\rU(\widehat{V}_{l-k})$ as subgroups of
$ \rU(V_l)$ naturally.
Let 
$P(V_k^+,V_l)=L(V_k^+,V_l)N(V_k^+,V_l)$ be the parabolic subgroup of $\rU(V_l)$ stabilizing $V^+_k$, where $L(V_k^+,V_l)=\GL(V^+_k)\times \rU(\widehat{V}_{l-k})$ is the Levi component of $P(V_k^+,V_l)$ stabilizing $V^-_k$.

 Now we specify the unipotent $N(V_k^+,V_l)$.
 For $T\in \Hom(V_l,V_l)$, we define $T^*\in \Hom(V_l,V_l)$ such that 
 \begin{equation}\label{adjoint of T}
  \innvl{Tv}{v'}=\innvl{v}{T^*v'},\quad \mbox{for every }v,v'\in V_l.  
 \end{equation}
 Consider $\Hom(\widehat{V}_{l-k},V^+_k)\subset \Hom(V_l,V_l)$ and $\Hom\left(V^-_k, V^+_k\right)\subset \Hom(V_l,V_l)$ via the decomposition \eqref {V_l}. For each $b\in \Hom(\widehat{V}_{l-k},V^+_k)$ and $c\in \Hom\left(V^-_k, V^+_k\right)$, we define
 \begin{equation}\label{Nkl}
  u(b, c)=1+b-b^*+c .  
\end{equation}
Then 
\[
N(V_k^+,V_l)=\{u(b,c)|b\in \Hom(\widehat{V}_{l-k},V^+_k),c\in \Hom\left(V^-_k, V^+_k\right),c+c^*+b b^*=0 \}.
\]
In particular, the center of $N(V_k^+,V_l)$ is 
\[
Z(V_k^+,V_l)\coloneqq \{u(0,c)|\Herm(V^-_k,V^+_k)\},
\]
where 
$$\Herm(V^-_k, V^+_k):=\left\{c \in \Hom\left(V^-_k, V^+_k\right) | c^{*}+c=0\right\}.$$
\trivial[h]{
We have an exact sequence 
\[
\begin{tikzcd}[row sep=0em,ampersand replacement=\&]
0 \ar[r] \& Z(V_k^+,V_l) \ar[r] \& N(V_k^+,V_l) \ar[r,"{u(b,c)\mapsto b}"] \& \Hom(\widehat{V}_{l-k},V^+_k) \ar[r] \& 0.
\end{tikzcd}
\]
}


\subsection{The Mix Model for Weil representation} \label{Mix model section}
We recall some explicit formulas for the Weil representation by G\'erardin~\cite[\S 2, \S3]{Ger}. 

Let $\mathbb V$ be a $(-1)$-Hermitian space over $F$ with Hermitian form $\inn{}{}_\mathbb V$.
For all $v_1,v_2\in \mathbb V$, we define
\begin{equation*}\label{diamond inner}
\iinn{v_1}{v_2}_\mathbb V:= \begin{cases}
\half \inn{v_1}{v_2}_\mathbb V,  & \text{if $F=F_0$},\\
\inn{v_1}{v_2}_\mathbb V,  & \text{otherwise}, 
\end{cases}
\quad \text{for all $v_1,v_2\in \mathbb V$}.
\end{equation*}
Note that we assumed $p\neq 2$ if $F=F_0$. The Heisenberg group attached to $\mathbb V$ is given by  
\[
H(\mathbb V):=\Set{(v,t)\in \mathbb V\oplus F|t-t^\tau = \iinn{v}{v}_\mathbb V}
\]
with the group law
\[
(v_1,t_1)\cdot (v_2, t_2)=\left(v_1+v_2, t_1+t_2+\iinn{v_1}{v_2}_\mathbb V \right). 
\]
The additive group $F_0$ is identified with the center of $H(\mathbb V)$ under the following exact sequence:
\[
\begin{tikzcd}[row sep=0em,ampersand replacement=\&]
0 \ar[r] \& F_0 \ar[r, "{t\mapsto (0,t)}"] \& H(\mathbb V)  \ar[r,"{(v,t)\mapsto v}"] \& \mathbb V \ar[r] \& 0.
\end{tikzcd}
\]
%
We denote by $\rho_{\mathbb V}$ the unique irreducible representation of $H(\mathbb V)$ with  central character $\psi$, which is called the Heisenberg representation. In particular, if we have a polarization $\mathbb V=\mathbb V^+\oplus \mathbb V_0 \oplus  \mathbb V^-$, then 
$\rho_{\mathbb V}$ can be realized on the space
$
\bC[\mathbb V^-]\otimes \rho_{\mathbb V_0}
$
of $ \rho_{\mathbb V_0}$-valued functions on $\mathbb V^-$. For $f\in \bC[\mathbb V^-]\otimes \rho_{\mathbb V_0}$ and $x\in \mathbb V^-$, the action is given as follows:
\begin{equation}\label{Heisenberg representation}
\begin{aligned}
\left(\rho_{\mathbb V}(y^+,0)f\right)(x) & =
\psi(\tr_{F/F_0}(\inn{x}{y^+}_\mathbb V) ) f(x)\quad & y^+&\in \mathbb V^+ \\
\left(\rho_{\mathbb V}(y^-,0)f\right)(x) & = f(x+y^-)& y^- &\in \mathbb V^-\\
\left(\rho_{\mathbb V}(v_0,t)\right)f(x)& =\rho_{\mathbb V_0}(v_0,t) f(x) & (v_0,t) &\in H(\mathbb V_0). 
\end{aligned}
\end{equation}
Here $\tr_{F/F_0}\colon F\rightarrow F_0$ is the trace map. 

\trivial[h]{
A computation for the mixed model. Take a partial polarization $\mathbb V=\mathbb V^+\oplus \mathbb V_0\oplus \mathbb V^-$. Let $P(\mathbb V^+,\mathbb V)=M(\mathbb V^+,\mathbb V)N(\mathbb V^+,\mathbb V)$ be the parabolic subgroup of $\rU(\mathbb V)$ stabilizing $\mathbb V^+$. Let $(\mathbb V^+)^\perp=\mathbb V^+\oplus \mathbb V_0$ and $H((\mathbb V^+)^\perp)$ be the pull back of $(\mathbb V^+)^\perp$ along the surjective map 
\[
\begin{tikzcd}[row sep=0em,ampersand replacement=\&]
H(\mathbb V)  \ar[r,"{(v,t)\mapsto v}"] \& \mathbb V 
\end{tikzcd}
\]
Then we have an exact sequence 
\begin{equation}\label{Heisenbgerg mix}
 0\longrightarrow \mathbb V^+ \longrightarrow H((\mathbb V^+)^\perp) \longrightarrow H(\mathbb V_0)\longrightarrow 0.    
\end{equation}
The action of $P(\mathbb V^+,\mathbb V)$ on $(\mathbb V^+)^\perp/\mathbb V^+\cong \mathbb V_0 $ gives a surjective map 
\begin{equation}\label{projection to Levi}
 P(\mathbb V^+,\mathbb V)\longrightarrow \rU(\mathbb V_0).
\end{equation}
Let $\omega_0$ be the Heisenberg-Weil representation of $\rU(\mathbb V_0)\times H(\mathbb V_0)$. Pull back $\rho_0$ to $P(\mathbb V^+,\mathbb V)\rtimes H((\mathbb V^+)^\perp)$ via \Cref{Heisenbgerg mix} and \Cref{projection to Levi}, which we still denote by $\omega_0$. (Note that the map 
\[
 P(\mathbb V^+,\mathbb V)\rtimes H((\mathbb V^+)^\perp) \longrightarrow \rU(\mathbb V_0)\rtimes H(\mathbb V_0)
\]
given by projection maps \Cref{Heisenbgerg mix} and \Cref{projection to Levi} is not a group homomorphism. One need to be careful here. ) 
Consider the induced representation 
\[
\Ind_{P(\mathbb V^+,\mathbb V)\rtimes H((\mathbb V^+)^\perp)}^{P(\mathbb V^+,\mathbb V)\rtimes H(\mathbb V)}(\omega_0)
\]
of $P(\mathbb V^+,\mathbb V)\rtimes H(\mathbb V)$, which we denote by $\omega'$ (bad notation ohh). We have the following vector space isomorphism 
\[
\Ind_{P(\mathbb V^+,\mathbb V)\rtimes H((\mathbb V^+)^\perp)}^{P(\mathbb V^+,\mathbb V)\rtimes H(\mathbb V)}(\omega_0)
\cong \Ind_{H((\mathbb V^+)^\perp)}^{H(\mathbb V)}(\omega_0)\cong \bC[\mathbb V^-]\otimes  \omega_0
\]
given by restrictions. We denote the map from the first to the third by $T$. Now we compute the action (denote by $\omega$, bad notation ohh) of $P(\mathbb V^+,\mathbb V)$ on $\bC[\mathbb V^-]\otimes \omega_0$ via the isomorphism $T$. We focus on the action of $N(\mathbb V^+,\mathbb V)$. For $u(b,c)\in N(\mathbb V^+,\mathbb V)$ and $f\in \Hom(\bC[\mathbb V^-], \omega_0)$ and $x\in \mathbb V^-$, we have 
\[
\begin{split}
  \omega(u(b,c)) f(x)= &(T \omega'(u(b,c))T^{-1}f)(x)= (\omega'(u(b,c))T^{-1}f)(x)\\
  =&(T^{-1}f)(x u(b,c))= (T^{-1}f)(u(b,c)u(b,c)^{-1}x u(b,c))\\
  =& (T^{-1}f)(u(b,c)^{-1}\cdot x ) \quad \mbox{(Here $\cdot$ means action)}\\
  =& (T^{-1}f)(u(-b,c^*)\cdot x )\\
  =&  (T^{-1}f)(((c^*x, b^*x, x),0))  \quad \mbox{Here we write element as $H(\mathbb V)=\mathbb V\oplus F=(\mathbb V^+\oplus \mathbb V_0\oplus \mathbb V^-)\oplus F$ }\\
  =& (T^{-1}f)\left( ((c^*x, b^*x,0 ),\iinn{c^*x}{-x}_{V})+((0,0,x),0)\right) \\
  =&  \omega_0\left( ((c^*x, b^*x,0 ),\iinn{c^*x}{-x}_{\mathbb V}) \right)(T^{-1}f)\left(((0,0,x),0)\right)\\
  =& \omega_0\left( ((0, b^*x,0 ),\iinn{c^*x}{-x}_{\mathbb V}) \right)(T^{-1}f)\left(((0,0,x),0)\right)\\
  =&\omega_0\left( ((0, b^*x,0 ),\iinn{c^*x}{-x}_{\mathbb V}) \right)f(x).
\end{split}
\]
\begin{itemize}
    \item In the symplectic group case, we may further simplify it as 
\[
\begin{split}
 \omega(u(b,c)) f(x)=& \psi(\iinn{c^*x}{-x}_{\mathbb V})\omega_0\left( ((0, b^*x,0 ),0) \right)f(x)\\
 =& \psi(\half \inn{c^*x}{-x}_{\mathbb V})\omega_0\left(b^*x) \right)f(x)
 \end{split}
\]
Choose $c=c^*=-\half b b^*$, we have 
\[
\psi(\half \inn{c^*x}{-x}_{\mathbb V})= \psi(\frac{1}{4}  \inn{b^*x}{b^*x}_{\mathbb V})=1.
\]
We have 
\[
\omega(u(b,c)) f(x)=\omega_0\left(b^*x) \right)f(x).
\]
\item In the unitary group case, if $p\neq 2$, we can use the isomorphism
\[
\begin{split}
    H(\mathbb V) &\longrightarrow H(\Res_{F/F_0}(\mathbb V))\\
     (v,t) &\mapsto (v,t- \half \langle v, v\rangle_\mathbb V ).
\end{split}
\]
to transfer the action to $H(\Res_{F/F_0}(\mathbb V))$. It corresponds to
\[
\omega_0^{'}\left( ((0, b^*x,0 ),\iinn{c^*x}{-x}_{V}-\half \iinn{ b^*x} {b^*x }) \right).
\]
Here $\omega_0^{'}$ denotes the Heisenberg representation of $H(\Res_{F/F_0}(\mathbb V)$. If we choose $c=c^*=-\half b b^*$, then 
\[
\iinn{c^*x}{-x}_{\mathbb V}-\half \iinn{b^*x}{b^*x} =0. 
\]
We have 
\[
\omega(u(b,c)) f(x)=\omega'_0\left(b^*x) \right)f(x).
\]
\end{itemize}
In both cases, the action of the center $Z(\mathbb V^+,\mathbb V)$ is given by 
\[
 \omega(u(0,c)) f(x)=  \omega(u(0,c)) f(x)=\psi( \iinn{c^*x}{-x}_{\mathbb V})f(x)= \psi( \iinn{cx}{x}_{\mathbb V})f(x).
\]
Note that $c+c^*=0$ in this case.
}

\trivial[h]{Some computations on the isomorphism $V\otimes V'\cong \Hom(V,V')$. If we define $V\otimes V'$ to be the vector space over $F$ with multiplication $a(v\otimes v')=av \otimes v'=v\otimes av'$ and $(-1)$-Hermition form 
\[
\inn{v_1\otimes v'_1}{v_2\otimes v'_2}_{V\otimes V'}=\inn{v_1}{v_2}_{V}\inn{v'_1}{v'_2}_{V'}.
\]

Following Geradin, we define $\Hom(V,V')$ be the set of semi-linear map from $V$ to $V'$, i.e.,
\[
\Hom(V,V')=\{f(av)=a^\tau f(v) |v\in V \}.
\]
We can make $\Hom(V,V')$ an $F$-vector space by defining $(a\cdot f)(v)=af(v)$. Then the map 
\[
\begin{split}
  \mathcal L:  V\otimes V'&\longrightarrow \Hom(V,V')\\
    v\otimes v'& \mapsto \left( \cdot \mapsto \innv{v}{\cdot}v' \right)
\end{split}
\]
is a linear isomorphism. We denote the image by $T_{v\otimes v'}$. For each $T\in \Hom(V,V')$, we define $T^\dag\in \Hom(V',V)$ (understand as semi-linear map from $V'$ to $V$) by 
\[
\inn{T v}{v'}_{V'}=\epsilon \inn{T^\dag v'}{v}_{V}=\inn{v}{T^\dag v'}^\tau_{V} \quad \mbox{for all $v\in V, v'\in V'$.}
\]
Then $T^\dag_{v\otimes v'}$ is the semi-linear map $(\cdot \mapsto \epsilon \inn{v'}{\cdot}_{V'} v)$
One can check that 
\[
\inn{T_{v\otimes v'}\widetilde v }{\widetilde v'}_{V'}=\inn{v}{\widetilde v}_V \inn{v'}{\widetilde v'}_{V'}=\epsilon \inn{T^\dag_{v\otimes v'} \widetilde v'}{\widetilde v}_{V}\quad \mbox{for all $\widetilde v\in V$ and $\widetilde v'\in V'$}.
\]
We want to have 
\[
\inn{v_1\otimes v_1'}{v_2\otimes v_2'}_{V\otimes V'}=\tr(T^\dag_{v_1\otimes v_1'}T_{v_2\otimes v_2'})
\]
The left hand side is 
\[
\inn{v_1}{v_2}_{V}\inn{v'_1}{v'_2}_{V'}.
\]
The right hand side $T^\dag_{v_1\otimes v_1'}T_{v_2\otimes v_2'}$ is 
\[
\cdot \mapsto \inn{v_2}{\cdot}_{V}^\tau v_2' \mapsto \epsilon \inn{v_2}{\cdot}^\tau \inn{v_1'}{v_2'} v_1= \inn{\cdot }{v_2} \inn{v_1'}{v_2'} v_1 .
\]
Therefore, the trace is 
\[
\inn{v_1}{v_2}\inn{v_1'}{v_2'}. 
\]
Under the isomorphism $\mathcal L$, we translate the action of $\rU(V)\times \rU(V')$ on $V\otimes V'$ to $\Hom(V,V')$, it is given by 
\[
(g,g')f(v)= g'f(g^{-1}v) \quad g\in \rU(V), g'\in \rU(V'), v\in V, f\in \Hom(V,V'). 
\]
}

\trivial[h]{Translate the mixed model via $\mathcal L$: 
Compute the action of $N(V^+,V)$ on $\bC[\Hom(V^+, V')]\otimes \omega_0$. For $x\in \Hom(V^+,V')$ and $u(b,c)\in N(V^+,V)$, we have 
\[
\begin{split}
  \omega(u(b,c)) f(x)= &(T \omega'(u(b,c))T^{-1}f)(x)= (\omega'(u(b,c))T^{-1}f)(x)\\
  =&(T^{-1}f)(x u(b,c))= (T^{-1}f)(u(b,c)u(b,c)^{-1}x u(b,c))\\
  =& (T^{-1}f)(u(b,c)^{-1}\cdot x ) \quad \mbox{(Here $\cdot$ means action)}\\
  =& (T^{-1}f)(x \cdot  u(b,c) )  \quad \mbox{(Here $\cdot$ means compose)}\\
  =&  (T^{-1}f)((x, x\cdot b, x\circ c),0))  \quad \mbox{$\Hom(V,V')=\Hom(V^+,V')\oplus \Hom(V_0,V')\oplus \Hom(V^-,V')$ }\\
  =& (T^{-1}f)\left( ((0, x\cdot b,x\cdot c ),-\tr((x\cdot c)^\dag x))+((x,0,0),0)\right) \\
  =&  \omega_0\left( ((0, x\cdot b,0 ),,-\tr((x\cdot c)^\dag x)) \right)(T^{-1}f)\left(((x,0,0),0)\right)\\
  =&\omega_0\left( ((0, x\cdot b,0 ),,-\tr((x\cdot c)^\dag x)) \right)f(x)\\
  =& \omega_0\left( (x\cdot b, -\tr(c^* x^\dag x)) \right)f(x).
\end{split}
\]
The last step we use $(x\cdot c)^\dag= c^* \cdot x^\dag$. If $b=0$, then 
\[
\omega(u(0,c)) f(x)=\omega(u(0,c)) f(x)= \psi\left( -\tr(c^* x^\dag x)) \right)f(x)= \psi\left( \tr(c x^\dag x)) \right)f(x). 
\]
}

Recall that $V_l\otimes_F  V'_{l'}$ is a $-1$-Hermitian vector space. We define $\Hom(V_l,V'_{l'})$ to be the set of all $F$-conjugate linear maps from $V_l$ to $V'_{l'}$, i.e., 
\[
\Hom(V_l,V'_{l'})=\{T: V_l\rightarrow V'_{l'}|T(av)=a^\tau T(v),\mbox{for all $a\in F, v\in V_l$} \}.
\]
We can make $\Hom(V_l,V'_{l'})$ an $F$-vector space by defining $(a\cdot T)(v)=aT(v)$. Then the map 
\begin{equation}\label{identify hom tensor}
\begin{split}
   V_l\otimes V'_{l'}&\longrightarrow \Hom(V_l,V'_{l'})\\
    v\otimes v'& \mapsto \left( \cdot \mapsto \innv{v}{\cdot}v' \right)
\end{split}
\end{equation}
is a linear isomorphism. 
Transporting via this isomorphism, the $(-1)$-Hermitian form on $\Hom(V_l,V'_{l'})$ is given by 
\[
\inn{T_1}{T_2}_{\Hom(V_l,V'_{l'})} =  \tr(T_1^\dag\, T_2) \quad \text{for all $T_1,T_2\in \Hom(V_l,V'_{l'})$}.
\]
Here for $T\in \Hom(V_l,V'_{l'})$, we define $T^\dag\in \Hom(V'_{l'},V_l)$ by requiring 
\[
 \rinn{T v}{v'}_{V'_{l'}} = \inn{v}{T^\dag v'}^\tau_{V_{l}}\quad \text{for all $v\in V_l$ and  $v'\in V'_{l'}$}.
\]

\trivial[h]{
I take another approach. 
We identify $V_l\otimes_F V'_l$ with the left $F$-vector space \zjl{Why we need to emphasis left here. } $\Hom(V_l,V'_l)$  such that $v\otimes v'\mapsto (a\mapsto \innv{a}{v}v')$ for all $v\in V_l$ and  $v'\in V'_{l'}$. 
By transporting the structures, the $(-1)$-Hermitian form on $\Hom(V_l,V'_{l'})$ is given by 
\[
\inn{x_1}{x_2}_{V_l\otimes V'_{l'}} =  \tr(x_1^\dag\, x_2) \quad \text{for all $x_1,x_2\in \Hom(V_l,V'_{l'})$}.
\]
Here $x_1^\dag\in \Hom(V'_{l'},V_l)$ is determined by
\[
\inn{T v}{v'}_{V'}=\epsilon \inn{T^\dag v'}{v}_{V}=\inn{v}{T^\dag v'}^\tau_{V} \quad \mbox{for all $v\in V, v'\in V'$.}
\]
}

Taking the partial polarization 
\[
\Hom(V_l,V'_{l'})=\Hom(V_l^+,V'_{l'})\oplus \Hom(V_0,V'_{l'})\oplus \Hom(V_l^-,V'_{l'}), 
\]
the mixed model of $\Weil{V_l,V_l'}$ can be realized on the space 
\begin{equation}\label{Mix model}
\Weil{V_l,V_l'}=\mathbb C[\Hom(V_l^+,V'_{l'})]\otimes \Weil{V_0,V_{l'}'}
\end{equation}
of $\Weil{V_0,V_{l'}'}$-valued functions on $\Hom(V_l^+,V'_{l'})$. For $f\in \Weil{V_l',V_l}$ and $x\in \Hom(V_l^{+},V'_{l'})$, the action of 
$P(V_l^+,V_l)\times \rU(V'_{l'})$ is given as follows:
\begin{equation}\label{mix model equation}
\begin{aligned}
\big(\Weil{V_l,V_{l'}'}(g') f\big)(x) = & \Weil{V_0,V_{l'}'}(g')  f((g')^{-1}\com  x) & &  g' \in \rU(V'_{l'}),\\
\big(\Weil{V_l,V_{l'}'}(a) f\big)(x) = & \chi(\det_{V_l^+} (a)) f( x\com a) & & a \in \GL(V_l^+),\\
\big(\Weil{V_l,V_{l'}'}(g_0)f\big)(x) = & \Weil{V_0,V_{l'}'}(g_0)  f(x) & & g_0\in \rU(V_0),\\
\big(\Weil{V_l,V_{l'}'}(u(b,c))f\big)(x) = & \rho_{V_0,V_{l'}'}\left(x\com b, -\iinn{x\com c}{x}_{\Hom(V_l,V'_{l'})}\right) f(x) & &    \begin{array}{l}b  \in   \Hom(V_0,V_l^+), \\
 c \in \Hom (V_l^-,V_l^+).\\\end{array}\\
\end{aligned}
\end{equation}


\begin{itemize}
\item  $\chi$ is a quadratic character of $F^\times $ defined in \Cref{chi};
\item $u(b,c)$ is defined in \Cref{Nkl};
\item $\rho_{V_0,V_{l'}'}$ is the Heisenberg representation of the group
$H(\Hom(V_0, V'_{l'}))$ with central character $\psi$, see \Cref{Heisenberg representation};
\end{itemize}
Next, we describe the action of $t_k$ for $1\leq k\leq l$. Write $x=(x_1,\cdots,x_l)$ for $x\in \Hom(V_l^+,V'_{l'})$ and $x_i\in  \Hom(\langle v_i\rangle, V'_{l'})$. Then the action of $t_k$ is given by the partial Fourier transform on $\Hom(\langle v_k\rangle, V'_{l'})$: 
\begin{equation}\label{Partial Fourier}
\begin{split}
& \left(\Weil{V_l,V'_{l'}}(t_k)f\right)(x)  \\
&  :=  \gamma_{V'_0} \displaystyle\int_{\Hom(\langle v_k\rangle, V'_{l'})} \hspace{-4em} f(x_1,\cdots, x_{k-1},y,x_{k+1},\cdots, x_l) \psi\left(\tr_{F/F_0} \left(\inn{y}{ x_k\com t_k}_{\Hom(V_l, V'_{l'})}\right)\right)dy.
\end{split}
\end{equation}
\begin{itemize}
    \item Here and below, we choose the Haar measure on any vector space over $F$ so that the volume of the vector space is the square root of its cardinality. 
\item $\gamma_{V'_0}$ is defined by 
\begin{equation}\label{Weil index}
\gamma_{V'_0}\coloneqq \begin{cases}
(-1)^{\dim V'_0} \quad &\mbox{in Case $A$},\\
1 \quad &\mbox{in Case $B$ and $D$}, \\
\xi\left(\disc( V'_0)\right) \quad &\mbox{in Case $C$}, \\
\gamma_{\psi}(1) \xi\left(\disc( V'_0)\right) \quad &\mbox{in Case $\widetilde{C}$}, \\
\end{cases}
\end{equation}
where $\xi$ is  the unique quadratic character of $F^\times$ as in \Cref{xi}, and 
in case $\widetilde{C}$
\[
\gamma_\psi(1)=\int_{F} \psi\left(\half x^2 \right) dx 
\]
is the usual Weil index. 
\end{itemize}


Finally we describe the action of $w$ for $\sfww\in \sfW_l$, where $w$ is the fixed lift of $\sfww$ in \Cref{Parabolic subgroup}. Write $\sfww=\sfdd \sfss$ for $\sfdd\in \sfD_l, \sfss\in \sfS_l$ by \Cref{sfW semi}. Since $\sfD_l$ is the set of distinguished representatives of the right coset $\sfW_l/ \sfS_l$, we have $w=d s$, where $d$ and $s$ are lifts of $\sfdd$ and $\sfss$. We may further write $\sfdd=\sftt_{i_1}\cdots \sftt_{i_k}$ uniquely for $1\leq i_1<\cdots < i_k\leq l$. It is routine to check that 
\[
l(\sfdd)=l(\sftt_{i_1})+\cdots +l(\sftt_{i_k}).
\]
Therefore, we have $s\in \GL(V_l^+)$, $d= t_{i_1}\cdots t_{i_k}$, and the action of $w$ can be deduced from \Cref{mix model equation} and \Cref{Partial Fourier}. We describe it explicitly below. Let
\trivial[h]{Here we explain that if $\sfdd=\sftt_{i_1}\cdots \sftt_{i_k}$ for $i_1<\cdots < i_k$, then $l(\sfdd)=\sum_{j=1}^k l(\sftt_{i_j})$. We do it by direct computation, i.e., count how many positive is sent into negative roots. First we have 
\[
l(\sftt_{i})= 2i+1.
\]
The positive roots which send to negative by $\sftt_i$ are $\{\epsilon_{i}\pm \epsilon_{k}, 2\epsilon_i| k=i+1,\cdots, l\}$, which we denote by $R_{\sftt_i}$. Then it is easy to check that  the positive roots which send to negative by $\sfdd$ is $\bigsqcup_{j=1}^k R_{\sftt_{i_j}}$. Therefore, we have 
$l(\sfdd)=\sum_{j=1}^k l(\sftt_{i_j})$. 
} 
\begin{equation*}\label{Vw+}
 V_{\sfww}^+= w^{-1} V_l^+\cap V_l^+, \quad V_{\sfww}^-= w^{-1} V_l^-\cap V_l^+   
\end{equation*}
and 
\begin{equation}\label{bxw}
 \bX_{\sfww}^+ : =\Hom(V_{\sfww}^+, V'_{l'}) ,\quad \bX_{\sfww}^- :=\Hom(V_{\sfww}^-,V'_{l'}).   
\end{equation}
Write an element in $\Hom(V^+_l,V'_{l'})$ as  
$(x^+,x^-)$ with $x^+\in \bX_{\sfww}^+$ and $x^-\in \bX_{\sfww}^-$. Then the action of $w$ is given by:
\begin{equation}\label{Partial Fourier d}
\begin{split}
& \left(\Weil{V_l,V'_{l'}}(w)f\right)(x^+,x^-)\\
& = (\gamma_{V_0'})^{\iota(\sfww)}  \int_{\bX_{\sfww}^-}f(x^+\com s,y) \psi\left(\tr_{F/F_0}\left(\inn{y}{x^-\com  w}_{\Hom (V_l, V'_{l'})}\right)\right)dy,
\end{split}
\end{equation} 
where  
\begin{equation}\label{iotaw}
    \iota(\sfww):=\dim V_{\sfww}^-.
\end{equation}

\subsection{Kudla's filtration}  \label{Kudla's filtration}
Using the mixed model \Cref{Mix model} of $\Weil{V_l,V'_{l'}}$ , we compute the Jacquet module of $\left(\Weil{V_l,V'_{l'}}\right)$ with respect to the parabolic subgroup $P(V_l^+,V_l)=L(V_l^+,V_l)N(V_l^+,V_l)$. Let 
\[
\cZ = \set{x\in \Hom(V_l^{+}, V'_{l'})|\text{$\Im(x)$ is an isotropic subspace of $V'_{l'}$}}.
\]
For $0\leq k\leq l$, let
 \begin{equation}\label{Zk}
\cZ_k := \Set{x\in \cZ|\text{$x$ has rank $k$ }.}
\end{equation}
It is non-zero only if $k\leq  \text{split-rank of } V_{l'}'$.
 
\begin{prop}[Kudla's filtration]\label{Kudla's filtrationlem}
\begin{enumerate}[label=(\arabic*),wide=0em]
\item We have the following decomposition as $(P(V_l^+,V_l)/Z(V_l^+,V_l))\times \rU(V'_{l'})$-representations: 
\begin{equation*}\label{Kudla1}
\Weil{V_l,V'_{l'}}^{Z(V_l^+,V_l)}= \bigoplus_{k}  \bC[\cZ_k]\otimes \Weil{V_0,V'_{l'}},
\end{equation*}
where the sum is over 
$0\leq k\leq \min\set{l,\text{split-rank of
} V'_{l'} }.$ 
\item For $0\leq k\leq \min\set{l,\text{split-rank of
} V'_{l'} }$, we have 
\[
\begin{split}
& \lb\bC[\cZ_k]\otimes \Weil{V,V'_{l'}}\rb^{N(V_l^+,V_l)/Z(V_l^+,V_l)} 
\\
& \cong \Ind_{P(V_{l-k}^+,V_l^+) \times   \rU(V_0)\times P(V'^+_{k},V'_{l'})}^{\GL(V_l^{+})\times \rU(V_0)\times \rU(V'_{l'})} \left( (\chi\circ \det_{ V_{l-k}^+})\otimes \bC[\Isom(\widehat{V}^+_k,V'^+_k)]\otimes
\omega_{V_0, V'_{l'-k}}\right)
\end{split}
\]
as a $ \GL(V_l^{+})\times \rU(V_0)\times \rU(V'_{l'})$-representation, where 
\begin{itemize}
\item  $ \widehat{V}^+_k,V'^+_k$ are defined  as in \Cref
{Parabolic subgroup} if $k\leq l'$, and if $k>l'$,    $V'^+_k$ is the direct sum of $V'^+_{l'}$ and a maximal totally isotropic subspace of the orthogonal complement of $V'_{l'-k}$ in $V'_{0}$;
\item  $P(V_{l-k}^+,V_l^+)$ is the maximal parabolic subgroup of $\GL(V_l^{+})$ stabilizing $V^+_{l-k}$;
\item $P(V'^+_{k},V'_{l'})$ is the parabolic subgroup of $\rU(V'_{l'})$ defined in the same way as in  the end of \Cref{Parabolic subgroup};
\item  $\Isom(\widehat{V}_k^+,V_k'^+)$ is the set of invertible $F$-conjugate linear maps from $\widehat{V}_k^+$ to $V'^+_k$;
\item  $\GL(\widehat{V}_{k}^{+})$ and   $\GL(V'^+_{k} )$ acts on $\bC[\Isom(\widehat{V}^+_k,V'^+_k)]$ as
\[
\left((g,h)\cdot f\right)(x)= \chi(\det_{\widehat{V}_{k}^{+}}(g))f(h^{-1}\com x\com g)
\]
for $(g,h)\in \GL(\widehat{V}_{k}^{+})\times \GL(V'^+_{k} )$, $f\in \bC[\Isom(\widehat{V}^+_k,V'^+_k)]$ and $x\in \Isom(\widehat{V}^+_k,V'^+_k)$;
\end{itemize}
\end{enumerate}
\end{prop}
\begin{proof} The proof goes in the same way as \cite{Ku}*{Theorem~2.8}, see also \cite{MVW}*{3.IV}.  
\end{proof}
Let $\sigma$ and $\sigma'$ be irreducible theta cuspidal representations of $\rU(V_0)$ and $\rU(V_0')$ such that $\Theta_{V_0,V_0'}(\sigma)=\sigma'$. Then $(V_0',\sigma')$ is the first occurrence of $\sigma$ (see \ref{theta cuspidal}).
Let $P'_{l'}$ be the parabolic subgroup of $\rU(V'_{l'})$ defined in the same way as  $P_l$ in \Cref{Parabolic subgroup}.
For $0\leq k\leq  \min\set{l,\text{split-rank of
} V'_{l'} }$,  let
\begin{equation}\label{definition of cFk}
 \cF_k\coloneqq \Hom_{P_l\times P'_{l'}}\left(
\sigmaP\otimes \sigmaPp, \bC[\cZ_k]\otimes \Weil{V_0,V'_{l'}}\right).   
\end{equation} 

By  Frobenius reciprocity, $\cM$ defined in \eqref{eq:multsp} is naturally identified with  
\begin{equation}\label{eq:intersp}
\left(\Weil{V_l,V_{l'}'}\right)_{\sigmaP\boxtimes \sigmaPp}:=\Hom_{P_l\times P'_{l'}}(\sigmaP \boxtimes \sigmaPp,  \Weil{V_l,V'_{l'}})
\end{equation}
\begin{cor}\label{cor:Kudla}
We have 
\begin{equation*}
 \left(\Weil{V_l,V_{l'}'}\right)_{\sigmaP\boxtimes \sigmaPp}  = \bigoplus_{k=0}^{\min\set{l,l'}}\cF_k.
\end{equation*}
\end{cor}
\begin{proof}
Under the assumption of the first occurrence of $\sigma$, we have $ \cF_k=0$ for $k>l'$ and the equation follows.  
\end{proof} 

\section{Hecke algebra and first occurrence index} \label{Normalization of Hecke operator}

\trivial[h]{	\qcl{need revision. SHorter} \zjl{the following sentence need modified}
Let $\sigma$ be an irreducible theta cuspidal representation of $\rU(V_0)$. In \Cref{Harish Chandra series and Hecke algebra}, we introduce the unnormalized Hecke operator $T_{\sftt_l}\in \cH_l$. In this section, we will construct two eigenspaces of $T_{\sftt_l}$ and compute their eigenvalues by using theta lifts. As a consequence, we obtain a normalization of $T_{\sftt_l}$ and relate the first occurrence indexes $n_{\cV'_V}(\sigma)$ and $n_{\tcV'_V}(\sigma)$ with the parameter for $T_{\sftt_l}$. }
 
We related the structure  of  the Hecke algebra for a   theta  cuspidal representation to its first occurrence index.
This allows us to prove the conservation relation for theta cuspidal representations. 

\subsection{Hecke algebra and induction} 
  Let $G$ be a reductive group over $F_0$.
Let  $P$ be a parabolic subgroup of $G$ and $\sigma_P$ be an irreducible representation of  $P$. Following \cite{BK}, we define the Hecke algebra 
\[
\begin{split}
& \cH(G,P,\sigma_P)\\
=&\End_{G}\left(  \Ind_P^{G}(\sigma_P^{\vee})\right) \\   
=&\{F: G\rightarrow \End_{\mathbb C}(\sigma_P^\vee)|F(p_{1}gp_{2})  = \sigma_P^{\vee}(p_{1}) F(g)\sigma_P^{\vee}(p_{2}),\ \ \forall p_1,p_2\in P\} 
\end{split}
\]
with the product given by the convolution 
\[
F_1*F_2(g)=\int_G F_1(h)F_2(h^{-1}g)dh
\]
for $F_1,F_2\in \cH(G,P,\sigma_P)$. 
Here we normalize the Haar measure on $G$ so that $P$ has volume $1$. We define the Hecke algebra $\cH(G,P,\sigma_P^{\vee})$ in the same way by identifying $(\sigma_P^{\vee})^{\vee}$ with $\sigma_P$. For $T\in \End_{\mathbb C}(\sigma_P^\vee)$, we denote by $T^\vee \in \End_{\mathbb C}(\sigma_P)$ the adjoint of $T$. Then we have an anti-isomorphism   
\begin{equation}\label{anti-isomorphism}
\begin{split}
  \cH(G,P,\sigma_P)&\cong \cH(G,P,\sigma_P^{\vee})\\
  F &\mapsto F^{\vee}(g)= F(g^{-1})^{\vee}.
\end{split}
\end{equation}
  For a representation  $\pi$ of $G$, let $\pi_{\sigma_P}=\Hom_P(\sigma_P, \pi)$. Following \cite[p 591]{BKt}, we know that $\pi_{\sigma_P}$ is a left $\cH(G,P,\sigma_P)$-module via the action 
\begin{equation}\label{Hecke-action}
F\cdot \phi: v\mapsto \int_{G} \pi(g)\phi(F^{\vee}(g^{-1})v) dg, \quad \mbox{for $F\in \cH(G,P,\sigma_P), \phi\in \pi_{\sigma_P}, v\in \sigma_P $}.
\end{equation}

\trivial[h]{
Recall that the Harish Chandra series attached to $(P,\sigma_P)$ is the set of all irreducible summands of $\Ind_P^{G}(\sigma_P)$. The following theorem can be found in \cite{BKt} and \cite{BK}.

\begin{thm}
\begin{enumerate}
    \item The functor $\mathfrak M_{\sigma}:\pi \mapsto \pi_{\sigma_P}$ from the category of $G$-modules to the 
category of $\cH(G,P,\sigma_P)$-modules
induces a bijection from the Harish-Chandra series associated to $(P,\sigma_P)$ to simple $\cH(G,P,\sigma_P)$-modules. For any irreducible module $\eta$ of $\cH(G,P,\sigma_P)$, we denote by $\sigma_\eta$ the Harish-Chandra series corresponding to $\eta$.
\item The map 
\[
\mathfrak M_{\sigma}: \Hom_{G}(\sigma_\eta, \pi)\longrightarrow \Hom_{\cH(G,P,\sigma_P)}(\eta,\pi_{\sigma_P})
\]
is an isomorphism. Therefore, the dimension of the multiplicity space of $\sigma_\eta$ in $\pi$ equals to the dimension of the multiplicity space of $\eta$ in $\pi_{\sigma_P}$.
\end{enumerate}
\end{thm}
}
The following compatibility result of Hecke algebra modules and parabolic induction is well known, see \cite{BKt}*{(8.4)}. 
\begin{prop}\label{lem:Hind}
Let   $ P  \leq Q  $  be  two parabolic subgroups in $G$ with  corresponding
  Levi subgroups  $L, M$.
Assume that $L\subset M$ and $P_M = P\cap M$  is  a parabolic subgroup of $M$ with Levi component $L$. In particular, we have a natural inclusion of $N_{M}(L)/L$ to $N_{G}(L)/L$. Assume that $\sigma$ is an irreducible representation of $L$. 
Let $\sigma_P$ and $\sigma_{P_M}$ be the inflations of  $\sigma$ to $P$ and $P_M$ respectively.
\begin{enumerate}[label=(\arabic*),wide=0em]
\item 
 The natural map
$$ t:\cH(M,P_M,\sigma_{P_M})\rightarrow \cH(G,P,\sigma_P)$$ 
that sends  $F\in \cH(M,P_M,\sigma_{P_M})$ supported on $P_M s P_M$ for $s\in N_{M}(L)$
to the unique element in $\cH(G,P,\sigma_P)$ supported on 
$P s P$ such that $(t(F))(s) = F(s)$ is an injective ring homomorphism. 
\item  
For a representation  $\pi$ of $M$,  let $\pi_Q$ be the inflation of $\pi$ to $Q$.
Identify $\pi$ with the sections in $\Ind_{Q}^G \pi$ supported on $Q$. 
This identification is $P_M$-equivariant and so induces a natural embedding 
$
\pi_{\sigma_{P_M}} \longrightarrow (\Ind_{Q}^G \pi_Q)_{\sigma_P}. 
$
   Then the naturally induced map 
\[
\cH(G,P,\sigma_P)\otimes_{\cH(M,P_M,\sigma_{P_M})} \pi_{\sigma_{P_M}} 
\longrightarrow (\Ind_{Q}^G \pi_Q)_{\sigma_P} 
\]
is an isomorphism of $\cH(G,P,\sigma_P)$-modules. 
\end{enumerate}
\end{prop}
\trivial[h]{
\begin{proof}

The proof  consists of the analogs of \cite[(6.12),(7.9),(7.12)]{BKt},  and is explained as follows.
We use the notations and definitions in \cite[(6.1)-(6.12)]{BKt}, but  with $P_u=Q,N_u=U'$ and $J=P$, and $\tau=\sigma_P$. Then $J_u=U'$ and $J_l=\{e\}$. 
First, 
according to the definition in  \cite[(6.1),(6.5)]{BKt}, $(J,\sigma)$  is decomposed with respect to  $(M,P_u)$, and 
every element in $M$ is positive. 
Then (1) follows from the same proof as  \cite[(6.12)]{BKt}. 
Second, we do not need the discussion from \cite[(6.13)-(7.7)]{BKt}. The  $t_u$ defined  in the paragraph above \cite[(7.8)]{BKt} is just our $t$ in our situation. Then the  results in \cite[(7.9)]{BKt} hold in our setting (especially without the supposition (b) in \cite[(7.9)]{BKt}). More precisely, \cite[(7.9)(1)]{BKt}   holds by the same computation, since every element   in $M$ is positive. And \cite[(7.9)(1)]{BKt}      follows by definition and that $\tau|_{U}$ is trivial. 
Then   \cite[(7.11), (7.12)]{BKt} hold  in our setting  and we have the commutative diagrams in \cite{BKt}*{(8.4)} by the same proof. The second commutative diagram  gives (2).
\end{proof}
}

\subsection{Hecke algebras of type $BC_l$}\label{Harish Chandra series and Hecke algebra}
Let $\sigma\in \Irr(\rU(V_0))$ be theta cuspidal and $\sigmaP$ be as in \eqref {sigmaldef}. 
Let 
\[
\cH_l:=\cH(\rU(V_l),P_l,\sigmaP).
\]
We study the structure of $\cH_l$. For   $\sfww\in \sfW_l$, let $w$ be the lift of $\sfww$ fixed in \Cref{Parabolic subgroup} and define the unnormalized Hecke operator $T_{\sfww}\in \cH_l$ by 
\begin{equation}\label{unnormalized intertwining operator 2}
T_{\sfww}(g) := \begin{cases}
\sigmaP^{\vee}(p_1)\sigmaP^{\vee}(p_2) \quad &\mbox{if $g=p_1  w p_2$ for some $p_1,p_2\in P_l$,}\\
0 \quad &\mbox{if $g\notin P_l w P_l$}.
\end{cases}   
\end{equation}
Here, since  $\sigmaP^\vee (w p  w^{-1})=\sigmaP^\vee (  p  )$ for
 $ p\in P_l\cap w^{-1} P_l  w 
$,
  it is easy to check that
  $T_{\sfww}(g)$ is well defined and does not depend on  the choice of  $p_1,p_2$ if $g\in P_l w P_l$.
\begin{lemma}\label{basis of Hl}
The set $\set{T_{\sfww}|\sfww\in \sfW_l }$ forms a  basis of $\cH_l$ (as a vector space over $\bC$).
\end{lemma}
\begin{proof}
This follows by the argument \cite[pp. 635-636]{Spr} adapted to our situation, together with the analog of ``basic geometric lemma'' in \cite[\S 3, 1.2]{Be}. 
\trivial[h]{
By Frobenious 
\[
\begin{split}
\cH_l\cong \End_{\rU(V_l)}\left( \Ind_{P_l}^{\rU(V_l)}(\sigmaP^{\vee}),\Ind_{P_l}^{\rU(V_l)}(\sigmaP^{\vee})\right)\cong \End_{P_l}\left( \Ind_{P_l}^{\rU(V_l)}(\sigmaP^{\vee})|_{P_l},\sigmaP^{\vee}\right)
\cong \End_{L_l}\left(J_{\rU(V_l)}^{L_l}\circ \Ind_{P_l}^{\rU(V_l)}(\sigmaP^{\vee}),\sigmaP^{\vee}\right)
\end{split}
\]
It follows by the Geometric lemma (cite Berstein) that 
\[
J_{\rU(V_l)}^{L_l}\circ \Ind_{P_l}^{\rU(V_l)}(\sigmaP^{\vee})\cong \bigoplus_{w\in P_l\backslash \rU(V_l)/P_l} \Ind_{L_l^{w}}^{L_l}\circ \Ad(w)\circ J_{L_l}^{(L_{l})_w}(\sigmaP^{\vee})
\]
(or introduce the Wely group and use $\sfww\in \sfW_l\backslash \sfW/\sfW_l$)
where $L_l^w= L_l\cap w^{-1}L_lw$ and $(L_l)_w=L_l\cap w L_l w^{-1}$ and $\Ad(w): \Rep((L_l)_w)\longrightarrow \Rep(L_l^w)$ is defined by $\Ad(w)\sigma(l)=wlw^{-1}$ for $\sigma\in \Rep((L_l)_w)$ and $l\in (L_l)_w)$. Then we have 
\[
\cH_l\cong \bigoplus_{w\in P_l\backslash \rU(V_l)/P_l}  \End_{L_l}\left(\Ind_{L_l^{w}}^{L_l}\circ \Ad(w)\circ J_{L_l}^{(L_{l})_w}(\sigmaP^{\vee}),\sigmaP^{\vee}\right)
\]
(use second Frobenious ?)
}
\end{proof}
It is  easy to check that 
\begin{equation}\label{unnormalized compose}
  T_{\sfww_1}T_{\sfww_2}=T_{\sfww_1\sfww_2} \quad  \mbox{if  $l(\sfww_1\sfww_2)=l(\sfww_1)+l(\sfww_2)$}.
\end{equation}
 Moreover, for every $\sfww\in \I$, $T_{\sfww}^2$ is a linear combination of $T_{\sfww}$ and $T_{\sfee}$, i.e., we have a quadratic relation
\[
T_{\sfww}^2 =  C_{\sfww} T_{\sfww} + A_{\sfww}\, T_{\sfee}
\]
for some $C_{\sfww},A_{\sfww}\in \mathbb C$. The following lemma gives the quadratic relation for $T_{\sfss_i}$ for  $i=1,...,l-1$. 
\begin{lemma}\label{quadratic relation 1}
For $1\leq i< l-1$, we have 
\[
(T_{\sfss_i}+1)(T_{\sfss_i}-q) = 0.
\]
\end{lemma}
\begin{proof}
This follows by a $\GL_2$-calculation, see \cite{L}*{Proposition~5.8}. 
\end{proof}
We will compute the quadratic relation for $T_{\sftt_l}$ using theta lifts in \Cref{The quadratic relation}.

\subsection{Explicit formula of $T_w$ actions}
\trivial[h]{For $\sfww\in \sfW_l$, let $U^+_{\sfww}=U_l\cap w^{-1}U_lw$ and $ U^-_{\sfww}=U_l\cap w^{-1}U^{op}_lw$, where $U^{op}_l$ is the unipotent subgroup of the opposite parabolic $P_l^{op}$ of $P_l$ (cf. add). Then we have $U_l=U_{\sfww}^+\cdot U_{\sfww}^-$ and the map 
\begin{equation*}\label{Double coset}
    \begin{split}
        U^{-1}_{\sfww}\times P_l&\longrightarrow P_l w P_l \\
        (u,p)& \mapsto u w p
    \end{split}
\end{equation*}
is a bijection.}
\trivial[h]{
For $\sfww\in \sfW_l$, let 
\begin{equation}
U_{\sfww}=U_l\cap w^{-1}U^{op}_lw,    
\end{equation}
where $U^{op}_l$ is the unipotent subgroup of the opposite parabolic $P_l^{op}$ of $P_l$ (cf. add). 
}
Fix a $BN$-pair of $\rU(V_l)$ such that $B\subseteq P_l$. Let  $\Sigma$ be the root system and  $\Sigma^+\subseteq \Sigma$ the set of positive roots. 
For $\sfww\in \sfW_l$, let 
\begin{equation*}\label{definition Uw}
U_{\sfww}:=\prod_{\alpha\in \Sigma^+,\ \sfww^{-1}\alpha<0,\ U_\alpha\in U_l} U_{\alpha}   
\end{equation*}
where $U_{\alpha}$ is the one parameter subgroup of $\rU(V_l)$ corresponds to $\alpha$. Then the map
\begin{equation}\label{Double coset}
    \begin{split}
        U_{\sfww}\times P_l \longrightarrow P_l w P_l, \quad 
        (u,p)  \mapsto u w p
    \end{split}
\end{equation}
is a bijection. For a representation $\pi$ of $\rU(V_l)$ and    $\phi\in \pi_{\sigmaP}=\Hom_{P_l}(\sigmaP,\pi)$, it follows directly by  \Cref{Hecke-action}, \Cref{unnormalized intertwining operator 2} and \Cref{Double coset} that
\trivial[h]{
: \qcl{define $U_w$ here is better}
\begin{equation}
\begin{split}
(T_w\cdot \phi)(v)=& \int_{\rU(V_l)} \pi(g)\phi(T_w^{\vee}(g^{-1})v) dg
= \int_{\rU(V_l)}  \pi(g)\phi((T_w(g))^{\vee} v) dg\\
=&\sum_{u\in U_w} \int_P \pi(u\dot w p)\phi((T_w(u\dot w p))^{\vee} v)   dp=\sum_{u\in U_w} \int_P \pi(u\dot wp) \phi((\sigmaP^{\vee}(p))^{\vee} v)   dp\\
= &\sum_{u\in U_w} \int_P \pi(u\dot wp) \phi(\sigmaP(p^{-1}) v)  dp=\sum_{u\in U_w} \int_P \pi(u\dot wp)\pi(p^{-1}) \phi( v)   dp\\
=&\sum_{u\in U_w} \pi(u\dot w)\phi(v).
\end{split}
\end{equation}
for any $\phi\in \pi_{\sigmaP}$ and $v\in \sigmaP$. Therefore, we have 
}
\begin{equation}\label{actformula}
T_{\sfww}\cdot \phi=\sum_{u\in U_{\sfww}} \pi(u w)\phi.
\end{equation}
\trivial[h]{
For $\sfww \in \sfW_l$, let
$
U_{\sfww} 
$ be a product of  one-parameter unipotent subgroups such that map 
\begin{equation*}\label{Double coset wrong}
    \begin{split}
        U_{\sfww}\times P_l&\longrightarrow P_l w P_l \\
        (u,p)& \mapsto u w p
    \end{split}
\end{equation*}
is a bijection, as in \cite[3.2]{BT1}.
\zjl{Maybe it is better to give a explicit definition of $U_w$? . For $w \in \sfW_l$, let 
\begin{equation}\label{defintion of Uw}
 U_w=\prod_{\alpha >0, w^{-1}\alpha<0} U_\alpha,   
\end{equation}
where $\alpha$ runs over all the positive roots of $\rU(V_l)$ (Need to fix a Borel blabla)and $U_\alpha$ is the  one-parameter unipotent subgroups corresponds to $\alpha$. Then we have 
\begin{equation*}\label{Double coset}
    \begin{split}
        U_{w}\times P_l&\longrightarrow P_l\dww P_l \\
        (u,p)& \mapsto u \dww p
    \end{split}
\end{equation*}
is a bijection (cf. \cite[3.2]{BT1}). 
}}

 \trivial[h]{
In \Cref{Some explicit formulas on the Hecke algebra actions}, we will use \Cref{actformula} to study the structure of $\left(\Weil{V_l.V_l'}\right)_{\sigmaP\otimes \sigmaPp}$ as $(\cH_l\times \cH(\rU(V'_{l'}),P',\sigmaP'))$-module. To simplify the notation, we denote by  $\cH_l=\cH_l$ and $\cH'_{l'}=\cH(\rU(V'_{l'}),P',\sigmaP'))$. 
}

\subsection{First occurrence}\label{Baby case}
Assume that  the first occurrence of $\sigma$ with respect to generalized Witt tower  $\cV'_{V}$ (see \ref{definition of Witt tower} and recall that $V_0=V$) is achieved on $V'_0$. Let  $\sigma'=\Theta_{V_0,V'_0}(\sigma)$. By the argument in \cite{MVW}*{3.IV} adapted to our situation, one can see that $\sigma'$ is irreducible and theta cuspidal. In particular, the multiplicity space 
\begin{equation}\label{iota0}
   \left(\Weil{V_0,V'_0}\right)_{\sigma\boxtimes \sigma'}=
\Hom_{\rU(V_0)\times \rU(V'_0)}(\sigma\boxtimes \sigma', \Weil{V_0,V'_0})
\end{equation}
is 1-dimensional and we fix a non-zero element $\iota_0$ inside it.
  We use the mixed model \Cref{Mix model} of $\Weil{V_l,V'_0} $ and view 
the multiplicity space
\[
\big(\Weil{V_l,V'_0}\big)_{\sigmaP \boxtimes \sigma'} = \Hom_{P_l\times \rU(V'_0)}(\sigmaP \boxtimes \sigma',\bC[\Hom(V_l^+, V'_0)]\otimes \Weil{V_0,V'_0})
\] 
as a subspace of $\bC[\Hom(V_l^+, V'_0)]\otimes \left(\Weil{V_0,V'_0}\right)_{\sigma\boxtimes \sigma'}$.
For a subset $\cS \subset \Hom(V_l^+, V'_0)$, let $I_{\cS}$ be the characteristic function on $\cS$.
If $\cS=\{v\}$ is a singleton, we abbreviate $I_\cS$ to $I_v$.

\begin{lemma}\label{lem:baby}
The multiplicity space $\big(\Weil{V_l,V'_0}\big)_{\sigmaP \boxtimes \sigma'}$ is   spanned by $I_{0}\otimes \iota_0$. 
Moreover,  
\begin{equation}
\begin{aligned}
T_{\sfss_i}\cdot (I_{0}\otimes \iota_0 )&=q (I_{0}\otimes \iota_0) \quad i=1,\cdots, l-1,\\
T_{\sftt_l}\cdot (I_{0}\otimes \iota_0) &= \gamma_{V'_0} q^{\dim V_0-\half \dim V'_0+\half \delta} (I_{0}\otimes \iota_0),
\end{aligned}
\end{equation}
 where $\gamma_{V'_0}$ is defined in \Cref{Weil index}.
\end{lemma}
\begin{proof}
By \Cref{cor:Kudla}, $\big(\Weil{V_l,V'_0}\big)_{\sigmaP \boxtimes \sigma'} = \cF_0$.
Then the first part of the lemma is a consequence of the fact that $\iota_0$ spans $\big(\Weil{V_0,V'_0}\big)_{\sigma \boxtimes \sigma'}$. 

Next, we compute the action of $T_{\sfss_i}$ on $I_{0}\otimes \iota_0$. 
\[
\begin{split}
  T_{\sfss_i}\cdot (I_{0}\otimes \iota_0) =&\sum_{u\in U_{\sfss_i}} \omega_{V_l,V'_0}(u s_i) (I_{0}\otimes \iota_0)\quad \mbox{by \Cref{actformula}} \\
  =&\sum_{u\in U_{\sfss_i}}\chi(\det_{V_l^+}(us_i))(I_{  0\com us_i}\otimes \iota_0) \quad \mbox{by  \Cref{Mix model} and $U_{\sfss_i}\subseteq \GL(V_l^+)$}\\
  =&q (I_{0}\otimes \iota_0)\quad \mbox{by $|U_{\sfss_i}|=q$ and $\det_{V^+_l}(u\, s_i)=1$}.
\end{split}
\]
Let $\lambda_{\sftt_l}\in \bC$ be the scalar such that  
$T_{\sftt_l}\cdot (I_{0}\otimes \iota_0)=\lambda_{\sftt_l} (I_{0}\otimes \iota_0)$. 
To compute $\lambda_{\sftt_l}$, it is sufficient to compute $\left(T_{\sftt_l}\cdot(I_{1}\otimes \iota_0)\right)(0)$. 
Note that $U_{t_l}\subseteq N(V_l^+,V_l)$ and we have 
\begin{equation}\label{shape of U_t}
U_{t_l}= \{u(b,c)|b\in \Hom(V_0,\langle v_l \rangle ), c\in \Hom(\langle v_{-l} \rangle ,\langle v_l \rangle ),c+c^*+bb^*=0 \},
\end{equation} 
and $ |U_{t_l}|=q^{\dim V_0+\half\deltaan}$ by \Cref{Nkl}.
Then it follows from  \Cref{mix model equation}, \Cref{Partial Fourier}, \Cref{actformula} and \Cref{shape of U_t} that 
\[
\begin{split}
& \left(T_{\sftt_l}\cdot(I_{0}\otimes \iota_0)\right)(0)\\
& = \sum_{u(b,c)\in U_{t_l}}\left(\omega_{V_l,V_0'}(u(b, c)\, t_l)(I_{0}\otimes \iota_0)\right) (0)\\
& = \gamma_{V'_0} q^{-\half \dim V'_0}\sum_{u(b,c)\in U_{t_l}}\left(\omega_{V_l,V_0'}(u(b,c))(I_{\Hom(\langle v_l\rangle, V'_0)}\otimes \iota_0)\right) (0)\\
& =\gamma_{V'_0} q^{-\half \dim V'_0} \sum_{u(b,c)\in U_{t_l}}\rho_{V_0,V'_0}((0\com b, \iinn{0\com c}{c}_{\Hom(V_l,V'_0)}  )  \iota_0 \\
& = \gamma_{V'_0} q^{\dim V_0-\half \dim V'_0+\half \delta}\iota_0
\end{split}
\]
Then the equation for the action of $T_{\sftt_l}$ follows.
\trivial[h]{
Therefore, we deduce that 
$
T_{\sftt_l}\cdot (I_{0}\otimes \iota_0)
= \gamma_{V'_0}q^{\dim V_0+\half(\deltaan-\dim V'_0)} (I_{0}\otimes \iota_0).
$}
\end{proof}

\subsection{The quadratic relation of $T_{\sftt_l}$}\label{The quadratic relation}
 Suppose that $(\widetilde{V}_0',\wtomega)$  is the first occurrence of $\sigma$ in the companion generalized Witt tower $\wtV'_{V}$ (see \Cref{definition of Witt tower}). 
 Then $\wtsigma':=\Hom_{\rU(V)}(\sigma,\wtomega)$ is theta cuspidal. 
In case $B$ and $D$, let $\gamma_{\wtV'_0}= - \gamma_{V'_0}$. 
Otherwise,  let $\gamma_{\wtV'_0}$ be defined  as in \Cref{Weil index}.  
By the definition of companion generalized Witt tower, 
a straightforward computation shows that 
\begin{equation}\label{gamma difference}
    \gamma_{V'_0}/ \gamma_{\widetilde{V}'_0}= -1.
\end{equation} 
Set 
\begin{equation}\label{eq:lambda}
\lambda_{\sftt_l} :=\gamma_{V'_0} q^{\dim V_0-\half \dim V'_0+\half \delta} 
\quad \mbox{and}\quad 
\wtlambda_{\sftt_l} := \gamma_{\widetilde{V}'_0} q^{\dim V_0-\half \dim \widetilde{V}'_0+\half \delta}. 
\end{equation}
The same proof of \Cref{lem:baby} gives the following result.
\begin{lemma}\label{lem:tbaby}
The multiplicity space $\wtomega_{\sigmaP \boxtimes \wtsigma'}$ is one-dimensional
and $T_{\sftt_l}$ acts on it by the scalar $\wtlambda_{\sftt_l}$.
\end{lemma}

\trivial[h]{
\begin{equation}\label{gamma sum}
  \gamma_{V'_0} \gamma_{\widetilde{V}'_0} =-\chi(-1). 
\end{equation}
\begin{proof}
This follows from a direct calculation using the definition of companion generalized Witt tower (cf. \Cref{definition of Witt tower}), \Cref{Weil index},  \Cref{determinant of tk} (in Case $B$ and $D$) and \Cref{Weil index square} (in Case $\widetilde{C}$). We omit the details here. 

We give the proof according to the types of $\rU(V)$. 
\begin{itemize}
    \item Type A: By \Cref{Weil index}, we have 
    \[
    \gamma_{V'_0}=(-1)^{\dim V'_0}\quad \mbox{and}\quad \gamma_{\widetilde{V}'_0}=(-1)^{\dim \widetilde{V}'_0}
    \]
    Since $\dim V'_0\not\equiv \dim \widetilde{V}'_0\bmod 2$, we deduce \Cref{gamma difference} and \Cref{gamma sum}.
    \item Type B and D: By \Cref{Weil index}, we have $\gamma_{V'_0}=1$. 
    On the other hand, it follows from \Cref{definition of Witt tower}, \Cref{determinant of tk} and \Cref{Partial Fourier} that 
 $\gamma_{\widetilde{V}'_0}=-1$. The rest is clear. 
\item Type C: By \Cref{Weil index}, we have  
\[
\gamma_{V'_0}=\xi\left(\disc( V'_0)\right)\quad  \mbox{and}\quad \gamma_{\widetilde{V}'_0}=\xi\left(\disc( \widetilde{V}'_0)\right).
\]
Since $\disc V'_0\not\equiv \disc \widetilde{V}'_0\bmod F^{\times 2}$, we deduce  \Cref{gamma difference} and \Cref{gamma sum}.
\item Type $\widetilde{C}$: By \Cref{Weil index}, we have  
\[
\gamma_{V'_0}=\gamma_{\psi}(1)\cdot \xi\left(\disc( V'_0)\right)\quad  \mbox{and}\quad \gamma_{\widetilde{V}'_0}=\gamma_{\psi}(1)\cdot\xi\left(\disc( \widetilde{V}'_0)\right).
\]
Since $\disc V'_0\not\equiv \disc \widetilde{V}'_0\bmod F^{\times 2}$, we deduce  \Cref{gamma difference} and \Cref{gamma sum}.
\end{itemize}
\end{proof}
}

 By \Cref{gamma difference}, we have
\begin{equation}\label{C_1/C_2 1}
\widetilde{\lambda}_{\sftt_l} /\lambda_{\sftt_l}= -q^{\half (\dim V'_{0}-\dim \widetilde{V}'_{0})}.
\end{equation}
In particular, $\lambda_{\sftt_l}\neq \widetilde{\lambda}_{\sftt_l}$. 
Combining \Cref{lem:baby} and \Cref{lem:tbaby}, we have the following proposition. 
\begin{prop}\label{quadratic relation2}
The quadratic relation of $T_{\sftt_l}$  is 
\begin{equation}\label{quadratic Ttl}
  (T_{\sftt_l}-\lambda_{\sftt_l})(T_{\sftt_l}-\widetilde{\lambda}_{\sftt_l})=0.  
\end{equation} 
\end{prop}

  \subsection{Conservation relation for theta cuspidal representations}
 By \Cref{eq:lambda}, we have 
\begin{equation}\label{C_1C_2 1}
\lambda_{\sftt_l}\widetilde{\lambda}_{\sftt_l}=  
 \gamma_{V'_0} \gamma_{\wtV'_0}q^{2\dim V_0-\half (\dim V'_{0}+\dim \widetilde{V}'_{0})+\delta}.
\end{equation}
On the other hand, by evaluating \Cref{quadratic Ttl}
at the identity $e\in \rU(V_l)$ and using \Cref{unnormalized intertwining operator 2}, we deduce  
\[
\begin{split}
-\lambda_{\sftt_l}\widetilde{\lambda}_{\sftt_l} \id_{\sigma_l^\vee} = T_{\sftt_l}^2(e)  & = \int_{\rU(V_l)} T_{\sftt_l}(h) T_{\sftt_l}(h^{-1}) dh \\  
& =\int_{U_{\sftt_l}} \int_{P_l} T_{\sftt_l}(ut_lp) T_{\sftt_l}(p^{-1}t_l t_{l}^{-2}u^{-1}) dudp \\
& = |U_{\sftt_l}| \sigma_l^{\vee}(t_l^{-2}).
\end{split}
\]
Note that
\[
t_l^2=\begin{cases}
e \quad & \mbox{if $\epsilon=1$},\\
\diag(1,\cdots, 1, -1, -1, 1, \cdots, 1)\quad &\mbox{if $\epsilon=-1$}. 
\end{cases}
\]
Therefore, 
$
\sigma_l^{\vee}(t_l^{-2})= \chi(-1) \id_{\sigma_l^{\vee}}$.
By the above calculations and \eqref{shape of U_t}, we deduce 
\begin{equation}\label{C_1C_2 2}
 \lambda_{\sftt_l}\widetilde{\lambda}_{\sftt_l} = 
-\chi(-1) |U_{t_l}| = -\chi(-1) q^{\dim V_0+\half \delta} .   
\end{equation}
Note that $\gamma_{V'_0}$ and $\gamma_{\wtV'_0}$ are roots of unity. 
Combining \Cref{C_1C_2 1} and \Cref{C_1C_2 2}, we get 
\begin{equation}\label{gamma sum}
  \gamma_{V'_0} \gamma_{\widetilde{V}'_0} =-\chi(-1). 
\end{equation}
and the conservation for theta cuspidal representations: 
\begin{prop}\label{Conservation} 
If $\sigma$ is theta cuspidal, then \Cref{thm:conservation} holds,  i.e., 
 \[
n_{\cV'_V}(\sigma) + n_{\tcV'_V}(\sigma) =  2\dim V + \deltaan. 
 \]
 \trivial[h]{
Let $\sigma$ be a theta cuspidal representation of $\rU(V_0)$. We have 
\begin{equation}\label{conservation relation 2}
\dim V'_0+\dim \widetilde{V}'_0=2
\dim V_0+\deltaan=
\begin{cases}
2\dim V_0+1\quad & \mbox{if $G=\rU(V_0)$ is of type $A$};\\
2\dim V_0 \quad &\mbox{if $G=\rU(V_0)$ is of type $B$ or $D$};\\
2\dim V_0+2 \quad &\mbox{if $G=\rU(V_0)$ is of type $C$ or $\widetilde{C}$}.\\
\end{cases}
\end{equation}
\qcl{ the whole purpose of setting up these notations is to avoid listing all the cases. }
}
\end{prop}

\subsection{Normalization of $\cH_l$}\label{structure of cHl}
Define the normalized Hecke operators
\begin{equation}\label{nomralized operator}
\begin{split}
   \cT_{\sfss_i} &:= T_{\sfss_i}\quad \text{for } i=1,\cdots, l-1,\quad \mbox{and}\\ 
   \cT_{\sftt_l} &:= -\lambda_{\sftt_l}^{-1}T_{\sftt_l}=-\gamma_{V'_0}^{-1} q^{-\dim V_0+\half \dim V'_0-\half \delta} T_{\sftt_l}.
\end{split}
\end{equation}
For a general $\sfww\in \sfW_l$, let $\sfww=\sfww_{1}\cdots \sfww_{r}$ be a reduced expression of $\sfww\in \sfW_l$ ($\sfww_{1},\cdots, \sfww_{r}\in \I$) and we define the normalized Hecke operator 
\begin{equation}\label{Normalized operator general}
\cT_{\sfww} := \cT_{\sfww_{1}}\cdots \cT_{\sfww_{i}}.  
\end{equation}
By \Cref{unnormalized compose}, we know that the definition does not depend on the choice of the reduced expression of $\sfww$ and 
\begin{equation}\label{Normalized compose}
\cT_{\sfww_1\sfww_2}=\cT_{\sfww_1} \cT_{\sfww_2}\quad \mbox{if $
l(\sfww_1\sfww_2)=l(\sfww_1)+l(\sfww_2)$}.
\end{equation}
\begin{proof}[Proof of \Cref{thm:norm}] 
By  \Cref{quadratic relation 1}, \Cref{quadratic relation2}, \Cref{C_1/C_2 1} and \Cref{nomralized operator},  we have 
$
(\cT_{\sfss_i}+1)(\cT_{\sfss_i}-q) = 0$ for $i=1,\cdots, l-1 , $ and $
   (\cT_{\sftt_l}+1)(\cT_{\sftt_l}-q^{\mu(\sigma)})=0.$
  Combining with \Cref{Normalized compose},
  the isomorphism 
  \[
  \sI \colon  \cH_l\xrightarrow{\ \ \cong \ \ } \sfH_{l,\mu(\sigma), \nu=q}. 
  \]
  is given by $\cT_\sfww\mapsto \sfT_{\sfww}\otimes_R 1$ for $\sfww\in \sfW_l$. 
\end{proof}

\def\wtsI{\widetilde{\sI}}
\begin{remark}\label{Different normalization}
\begin{enumerate}
\item There is a natural isomorphism $\kappa$ from $\sfH_{l,-\mu}$ to $\sfH_{l,\mu}$ defined by
\begin{enumerate}[label=(\alph*),wide=0pt]
\item $\kappa( \sfT_{\sfss_i,-\mu}) = \sfT_{\sfss_i,\mu} $ for $1\leq i\leq l-1$, and 
\item $\kappa( \sfT_{\sftt_l,-\mu}) = - \nu^{-\mu} \sfT_{\sftt_l,\mu}$.  
\end{enumerate}
 
    \item Switching the role of $\widetilde{V}_0'$ and $V'_0$, we get another normalization of
$T_{\sftt_l}$ by $\widetilde \cT_{\sftt_l}=-\widetilde{\lambda}_{\sftt_l}^{-1}T_{\sftt_l}$ and the quadratic relation 
\[
(\widetilde \cT_{\sftt_l}+1)(\widetilde \cT_{\sftt_l}-q^{-\mu(\sigma)})=0.
\]
Then we have an isomorphism $\wtsI\colon \cH_l\xrightarrow{\ \cong \ }\sfH_{l,-\mu(\sigma),\nu=q}$ by sending $\cT_{\sfss_{i}}$ to the specialization of $\sfT_{\sfss_i,-\mu}$  for $i=1,\cdots, l-1$ and $\widetilde \cT_{\sftt_{l}}$ to  the specialization of $\sfT_{\sftt_l,-\mu}$ at $\nu=q$.
\item  Clearly, we have the following commutative diagram:
\[
\begin{tikzcd}
 & \cH_l \ar[dl,"{\sI}"']\ar[dr,"{\wtsI}"]& \\
 \sfH_{l,\mu(\sigma),\nu=q}\ar[rr,"\kappa"] & & \sfH_{l,-\mu(\sigma), \nu=q}\\
\end{tikzcd}
\]
\end{enumerate}
\end{remark}


\trivial[h]{
Show the cocycle is trivial: Choose a Borel subgroup $B(V_l)=T(V_l)N(V_l)$ such that $B(V_l)\subseteq P_l$. Let $\pi$ be an irreducible substitution of $J_{L_l}^{T(V_l)}(\sigmaP^\vee)$. Therefore, $\pi$ is cuspidal. By \Cref{lem:Hind}, we have an embedding 
\[
\cH_l \hookrightarrow \cH(\rU(V_l),B(V_l),\pi).
\]
By the works of Howlett-Lehrer \cite[Theorem 4.14]{HL}, Lusztig \cite[8.6]{L2} and Geck \cite[Corollary 2]{Geck}, we know that $\cH(\rU(V_l),B(V_l),\pi)$ is isomorphic to a certain Hecke algebra of $\sfW$. Therefore, we deduce that the two cocycle of $\cH_l$ is also trivial. 
}

The rest of this section consists of several results which will be used later. 
By \Cref{lem:baby} and \Cref{nomralized operator}, we have the following:
 \begin{lemma}\label{Quadratic relation}    
The action of $\cH_l$   on $\big(\Weil{V_l,V'_0}\big)_{\sigmaP \boxtimes \sigma'}$  
 is by  the character $\varepsilon_l$ 
defined in \Cref{varepsilon of cHl}. \qed
\end{lemma}
 
By \Cref{sfttk}, \Cref{nomralized operator} and \Cref{Normalized compose}, we have the normalized Hecke operators
 \begin{equation}\label{Normalized t_k2}
 \cT_{\sftt_k} = -\gamma_{V'_0}^{-1} q^{-\dim V_0+\half \dim V'_0-\half \delta} T_{\sftt_k},\quad k=1,\cdots,l.
\end{equation}
Changing the role of $(V_0,\sigma,l)$ and $(V_0',\sigma',l')$, we have the normalized Hecke operators
 \begin{equation}\label{Normalized t_k23}
 \cT'_{\sftt_k}  = -\gamma_{V_0}^{-1} q^{-\dim V'_0+\half \dim V_0-\half \delta'} T'_{\sftt_k},\quad k=1,\cdots,l,
\end{equation}
where 
 $\delta'=2-\delta$, see   \eqref{deltadelta'}.


\subsection{Abundance of theta cuspidal dual pairs}\label{Abundance of cuspidal dual pairs} 
In this section, we fix a case listed in \Cref{different cases} but allow the spaces $V$ and $V'$ to vary. 
We say that $\mu $  in \eqref{eq:mu.range} is \emph{relevant} (with respect to our fixed case) if it  is relevant to some quadruple $(V,V', \sigma,\sigma')$ (see \Cref{eq:mu11}).

\begin{lemma}\label{conservation relevant}
    If $\mu$ is relevant, $-\mu$ and $-\mu-2$ are also relevant (in the same case). 
\end{lemma}
\begin{proof}
We only give the prove for Case (D) and the proofs of other cases are similar.
Let $(V,V', \sigma,\sigma')$  be a quadruple relevant to $\mu$. 
Suppose the first occurrence index of $\sigma$ with respect to the companion  generalized Witt tower $\tcV'_V$  is achieved on $\widetilde{V}'$. 
Then $\big(V,\wtV',\det\otimes \sigma,\Theta_{V,\wtV'}(\det\otimes\sigma)\big)$ is relevant to $-\mu$ by   \Cref{eq:mu}. In particular, $-\mu$ is relevant.

Let $\tcV$ be the Witt tower of even dimensional quadratic spaces which does not contain $V$. 
Suppose that $\wtV$ is the first occurrence of $\sigma'$ with respect to the generalized Witt tower $\set{(\wtV^\dag,\omega_{V',\wtV^\dag}) | \wtV^\dag\in \tcV}$.  Set $\wtsigma := \Theta_{V',\wtV}(\sigma')$. 
Recall $\delta$ and $\delta'$ defined in \eqref{deltadelta'}. 
According to \Cref{Conservation}, we have
\[
\dim \wtV + \dim V = 2 \dim V' + \delta'.
\]
Combining the above equation with \eqref{eq:mu11}, we have 
\[
\mu(\wtsigma) =  \dim V'- \dim \wtV  - \half \delta 
= -\mu -2.  
\]
Therefore, $(\wtV,V', \wtsigma, \sigma')$ is relevant to $-\mu-2$ which proves the lemma.  
\end{proof}

\begin{prop}\label{abundance}
Every $\mu$ in \eqref{eq:mu.range} is relevant.
\end{prop}

\begin{proof}
We only prove Case $(A)$ and other cases are similar. 
By \Cref{conservation relevant}, it is enough to show that $\half$ is relevant in Case $(A)$. Consider the dual pair $(V,V')$ for $\dim V=0$ and $\dim V'=1$. Then $\Theta_{V,V'}(\mathbbm 1_{\rU(V)})=\mathbbm 1_{\rU(V')}$. In other word, $(V,V', \mathbbm 1_{\rU(V)}, \mathbbm 1_{\rU(V')})$ is relevant to  $\dim V'-\dim V-\half \delta =\half$.
This finishes the proof. 
\end{proof}

\section{Multiplicity space as a Hecke module}\label{Some explicit formulas on the Hecke algebra actions}
 
We have related  Hecke algebras and  theta lifts in the last section.
Retrieve the set-up in \ref{mainsec}.
In the rest of the paper, 
we study the $\cH_l\times \cH_{l'}$-module $\cM =  \left(\Weil{V_l,V_{l'}'}\right)_{\sigmaP\boxtimes \sigmaPp}$ (see \eqref{eq:intersp}). 

\trivial[h]{
Now we consider the general case. Let $\mu$ be as in \Cref{abundance} and $\sigma$ and $\sigma'$ be irreducible theta cuspidal representation of $\rU(V_0)$ and $\rU(V_0')$ such that $\Theta_{V_0,V_0'}(\sigma)=\sigma'$ and 
\[
    \dim V'_0-\dim V_0-\half \delta=\mu.
\]
Fix the integers $l,l'\geq 0$, we want to study the structure of  $\left(\Weil{V_l,V_{l'}'}\right)_{\sigmaP\otimes \sigmaPp}$ as $\cH_l\otimes \cH_{l'}$-module.
}

 In this section, we first show that $\left(\Weil{V_l,V_{l'}'}\right)_{\sigmaP\boxtimes \sigmaPp}$  is an induced module for the product of $\cH'_{l'}$ and the parabolic subalgebra of $\cH_{l}$ generated by $\cT_{\sfss_i}$'s, $ i=1,\cdots, l-1$, using Kudla's filtration in \Cref{Kudla's filtrationlem}. 
 Then we use the induced structure to find a basis $\fG$ of $\left(\Weil{V_l,V_{l'}'}\right)_{\sigmaP\boxtimes \sigmaPp}$ and realize the basis explicitly in the mixed model. 
 Finally, we compute the action of $\cT_{\sftt_l}$ on $\fG$, and thus  we get a complete description of $\left(\Weil{V_l,V_{l'}'}\right)_{\sigmaP\boxtimes \sigmaPp}$ as an $\cH_l\times \cH_{l'}$-module.
\subsection{Action of $\cH_{V_l^+}\otimes \cH'_{l'}$}\label{The explicit Hecke action}
Retrieve the notations in \Cref{Parabolic subgroup}. For each $1\leq k\leq l$, let $B(V_{l-k}^+)$ and $B(\widehat{V}_k^+)$ be the Borel subgroups of $\GL(V_{l-k}^{+})$ and $\GL(\widehat{V}_k^+)$ stabilizing the flag 
\[
V_{1}^+\subseteq \cdots \subseteq V_{l-k}^+\quad \mbox{and }\quad \widehat{V}_1^+\subseteq \cdots \subseteq \widehat{V}_k^+ \quad \text{respectively.}
\]
Let $T_{l-k}=\GL_1(\langle v_1\rangle)\times \cdots \times \GL_1(\langle v_{l-k}\rangle)$ and $\widehat{T}_k=\GL_1(\langle v_{l-k+1}\rangle)\times \cdots \times \GL_1(\langle v_{l}\rangle)$ be the maximal tori inside  $B(V_{l-k}^+)$ and $B(\widehat{V}_k^+)$. We denote by $\chi_{B(V_{l-k}^+)}$ and $\chi_{B(\widehat{V}_k^+)}$ the inflations of $\chi^{\boxtimes {l-k}}$ and $\chi^{\boxtimes k}$ on $T_{l-k}$ and $\widehat{T}_k$. Let 
\[
\cH_{V_{l-k}^+}:=\cH(\GL(V_{l-k}^{+}),B(V_{l-k}^+), \chi_{B(V_{l-k}^+)})
\]
and 
\[
\cH_{\widehat{V}_{k}^+}:=\cH(\GL(\widehat{V}_k^{+}),B(\widehat{V}_k^+), \chi_{B(\widehat{V}_k^+)}).
\]
Then $\cH_{V_{l-k}^+}$ and $\cH_{\widehat{V}_k^+}$ are Hecke algebras of type $A_{l-k}$ and $A_{k}$ respectively with parameter $q$ and we identify them as sub-algebras of $\cH_l$ by 
\Cref{lem:Hind}. We define the relevant Hecke algebra $\cH_{V'^+_k}$ for $\GL(V'^+_k)$ similarly. Using the basis $\{v_{l-k+1},\cdots, v_l\}$ of $\widehat{V}_k^+$ and the basis $\{v_1',\cdots, v_k' \}$ of $V'^+_k$ (see \ref{Parabolic subgroup}), we identify $\GL(\widehat{V}_k^+), \GL(V'^+_k)$ with $\GL_k(F)$. 
It induces unique isomorphisms 
\begin{equation}\label{Hecke algenra ismorphism 3}
    \cH_{\widehat{V}_k^+}
    \cong \cH_{V'^+_k} \cong
    \cH_{\GL_k}\coloneqq 
    \cH(\GL_k(F),B_k(F),\mathbbm 1)
\end{equation}
preserving the support. 
Here $B_k(F)$ is the Borel subgroup of $\GL_k(F)$ and $\mathbbm 1$ is the trivial character of $B_k(F)$. 
Via the isomorphisms \Cref{Hecke algenra ismorphism 3}, we define an action of $\cH_{\widehat{V}_k^+}\times \cH_{V'^+_k}$ on $\cH_{\GL_k}$ by 
\begin{equation}\label{hecke type A}
 (\cT_{\sfww_1}, \cT'_{\sfww_2})\cdot  \cT = \cT'_{\sfww_2}\, \cT\, \cT_{\sfww_1^{-1}},  \end{equation}
 where $ \cT\in \cH_{\GL_k},  \cT_{\sfww_1}\in \cH_{\widehat{V}_k^+}, \cT'_{\sfww_2}\in \cH_{V'^+_k} $. 
For each $0\leq k\leq \min\set{l,l'}$, recall that we have defined $\cF_k$ in \Cref{definition of cFk}. It is easy to check that each $\cF_k$ is stable under the action of $\cH_{V_l^+}\otimes \cH'_{l'}$. 
\begin{prop}\label{Half way}
The $\cH_{V_l^+}\otimes \cH'_{l'}$-module
$\cF_k$ is isomorphic to
\begin{equation}\label{induced}
\left(\cH_{V_l^+}\boxtimes \cH'_{l'}\right)\otimes_{\left(\cH_{V^+_{l-k}}\otimes \cH_{\widehat{V}^+_{k}} \right)  \otimes \left( \cH_{V'^+_{k}}\otimes \cH'_{l'-k} \right)} \left( \mathbbm 1_{V_{l-k}^+} \boxtimes \cH_{\GL_k} \boxtimes \varepsilon_{l'-k} \right),
\end{equation}
where $\mathbbm 1_{V_{l-k}^+}$ is the trivial character of $\cH_{V_{l-k}^+}$ and $\varepsilon_{l'-k}$ is the character of $\cH'_{l'-k}$ defined in \Cref{varepsilon of cHl}. 

\end{prop}
\begin{proof}
It follows from \Cref{Kudla's filtrationlem} that 
\[
\begin{split}
\cF_k \cong  \bigg( \Ind_{P(V_{l-k}^+,V_l^+) \times   \rU(V_0)\times P(V'^+_{k'},V'_{l'})}^{\GL(V_l^{+})\times \rU(V_0)\times \rU(V'_{l'})}  \Big( (\chi\circ \det_{ V_{l-k}^+})& \otimes \bC[\Isom(\widehat{V}^+_k,V'^+_k)] \\ 
& \otimes
\omega_{V_0, V'_{l'-k}}\Big)\bigg)_{\chi_{B(V_l^+)}\boxtimes \sigma\boxtimes \sigmaPp}.
\end{split}
\]
We now analyze the $\left(\cH_{V^+_{l-k}}\boxtimes \cH_{\widehat{V}^+_{k}} \right)  \boxtimes \left( \cH_{V'^+_{k}}\boxtimes \cH'_{l'-k} \right)$-module  
\[
\left(
 (\chi\circ \det_{ V_{l-k}^+})\otimes \bC[\Isom(\widehat{V}^+_k,V'^+_k)]\otimes
\omega_{V_0,V'_{l'-k}}\right)_{\left(\chi_{B(V_{l-k}^+)}\boxtimes \chi_{B(\widehat{V}_k^+)}\boxtimes \sigma \right) \boxtimes \left(\chi_{B(V'^+_k)}\boxtimes \sigma'_{l'-k} \right)}.
\]
Note that 
\begin{itemize}
    \item $\left( \chi\circ \det_{ V_{l-k}^+}\right)_{\chi_{B(V_{l-k}^+)}}$ corresponds to the trivial character $\mathbbm 1_{V_{l-k}^+}$ of $\cH_{V^+_{l-k}}$;
    \item $\left( \bC[\Isom(\widehat{V}^+_k,V'^+_k)]\right)_{\chi_{B(\widehat{V}_k^+)}\boxtimes \chi_{B(V'^+_k)}}$ corresponds to $\cH(\GL_k)$ as a module of $\cH_{\widehat{V}^+_{k}}\otimes \cH_{V'^+_{k}}$ defined in \Cref{hecke type A};
    \item Since $(V',\sigma')$ is the first occurrence of $\sigma$, we know that 
    $\left(\Weil{V_0, V'_{l'-k}}\right)_{\sigma\boxtimes \sigma'_{l'-k}}$ is nonzero if and only if $k\leq l'$. Moreover, for $k\leq \min{\set{l,l'}}$, it follows from \Cref{Quadratic relation} that $\left(\Weil{V_0, V'_{l'-k}}\right)_{\sigma\boxtimes \sigma'_{l'-k}}$ corresponds to the character $\varepsilon_{l'-k}$ of $\cH'_{l'-k}$.
\end{itemize}
The lemma follows by \Cref{lem:Hind}. 
\end{proof}
\def\indbasisgen{\II^k_{1,1,1}}
\def\indbasis{\II^k_{\sfdd_1,\sfdd_2,\sfxx}}
\def\Mss{\left(\Weil{V_l,V'_{l'}}\right)_{\sigmaP\boxtimes \sigmaPp}}

Fix an nonzero vector $\indbasisgen\in \cF_k$ in the following one dimensional subspace of \eqref{induced}
\begin{equation}\label{IAk}
\left(\mathbf 1_{V_l^+}\otimes \mathbf 1_{l'}\right)\otimes
\left( \mathbbm 1_{V_{l-k}^+} \otimes \mathbf 1_{\GL_k} \otimes \varepsilon_{l'-k}, \right), 
\end{equation}
where $\mathbf 1_{V_l^+}, \mathbf 1_{l'}$ and $\mathbf 1_{\GL_k}$ are the vector spaces generated by the identity operators in $\cH_{V_l^+}, \cH'_{l'}$ and $\cH_{\GL_k}$ respectively. \trivial[h]{
Fix a non-zero vectors $\bfone_{V^+_{l-k}}$ and $\mathcal E_{l'-k}$ (\zjl{or  $\mathbf E_{l'-k}$ or  $\mathbf I_{l'-k}$})in 
$\bbone_{V^+_{l-k}}$ and $\varepsilon_{l'-k}$.  
Let 
\begin{equation}\label{IAk}
 \indbasisgen:= 
\left(\bfone_{V_l^+}\otimes \bfone_{l'}\right)\otimes
\left( \bfone_{V_{l-k}^+} \otimes \bfone_{\GL_k} \otimes \bfI_{l'-k} \right) \in \Mss 
\end{equation}
where $1_{V_l^+}$, $1_{l'}$ and $\bfone_{\GL_k}$ are the identity elements in $\cH_{V^+_l}$, 
$\cH'_{l'}$ and $\cH_{\GL_k}$ respectively.
}Let $\sfD_{k}$ and $\sfD'_{k}$ be the set of distinguished representatives of the left cosets 
$\sfS_l/(\sfS_{l-k}\times \sfS_k)$ and $\sfW_{l'}/(\sfS_{k} \times \sfW_{l'-k})$. 
For $(\sfdd_1,\sfdd_2,\sfxx)\in \sfD_k \times \sfD'_k\times \sfS_k$, let $\cT_{\sfdd_1}, \cT'_{\sfdd_2}, \cT_{\sfxx}$ and $\cT'_{\sfxx}$ be the normalized Hecke operators in $\cH_{V^+_{l}} , \cH'_{l'}, \cH_{\widehat{V}^+_k}$ and $\cH_{V'^+_k}$ respectively. Here we identify $\cH_{\widehat{V}^+_k}$ and $\cH_{V'^+_k}$ with $\cH_{\GL_k}$ via \Cref{Hecke algenra ismorphism 3}. For each $(\sfdd_1,\sfdd_2,\sfxx)\in \sfD_k \times \sfD'_k\times \sfS_k$, we define 
\begin{equation}\label{standard basis}
 \indbasis\coloneqq \cT_{\sfdd_1}\cT'_{\sfdd_2}\cT'_{\sfxx}\cdot \indbasisgen=   \cT'_{\sfdd_2}\cT_{\sfdd_1}\cT_{\sfxx^{-1}}\cdot \indbasisgen.  
\end{equation}
\begin{cor}
The following set is a basis of $\cF_k$
\begin{equation*}
  \{\indbasis |(\sfdd_1,\sfdd_2,\sfxx)\in \sfD_k \times \sfD'_k\times \sfS_k \}.
\end{equation*}
\end{cor}
In particular,  we get a basis of $\Mss$: 
\begin{equation}\label{geometric basis}
\fG :=   \bigcup_{k\leq \min\set{l,l'}}\fG_k, \quad \text{where}\quad 
\fG_k:=   \Set{ \indbasis| (\sfdd_1,\sfdd_2,\sfxx)\in \sfD_k \times \sfD'_k \times \sfS_k }.
\end{equation}

\begin{prop}\label{Hecke action half}
The action of $\cH_{V_l^+}\otimes \cH'_{l'}$ on $\fG$ is given as follows: 
\begin{enumerate}
    \item For $i=1,\cdots, l-1$, we have
\[
\cT_{\sfss_i}\cdot \indbasis=\begin{cases}
\II^k_{\sfss_i \sfdd_1,\sfdd_2,\sfxx} \quad & \mbox{if $\sfss_i \sfdd_1\in \sfD_k$ and $l(\sfss_i \sfdd_1)>l(\sfdd_1)$},\\
q\II^k_{\sfss_i \sfdd_1,\sfdd_2,x}+(q-1)\II^k_{\sfdd_1,\sfdd_2,\sfxx}\quad &\mbox{if $\sfss_i \sfdd_1\in \sfD_k$ and $l(\sfss_i \sfdd_1)<l(\sfdd_1)$},\\
q\II^k_{\sfdd_1,\sfdd_2,\sfxx}\quad & \mbox{otherwise.}
\end{cases}.
\]
\item For $i=1,\cdots, l'-1$,  we have 
\[
\cT_{\sfss'_i}\cdot \indbasis=
\begin{cases}
\II^k_{\sfdd_1,\sfss'_i\sfdd_2,\sfxx}\quad & \mbox{if $\sfss_i \sfdd_2\in \sfD'_k$ and $l(\sfss'_i \sfdd_2)>l(\sfdd_2)$},\\
q\II^k_{\sfdd_1,\sfss'_i\sfdd_2,x}+(q-1)\II^k_{\sfdd_1,\sfdd_2,\sfxx} \quad &\mbox{if $\sfss_i \sfdd_2\in \sfD'_k$ and $l(\sfss_i \sfdd_2)<l(\sfdd_2)$},\\
q\II^k_{\sfdd_1,\sfdd_2,\sfxx}\quad & \mbox{otherwise.}
\end{cases}
\]
\item Finally
\begin{equation*}
\cT_{\sftt'_{l'}}\cdot \indbasis=
\begin{cases}
\II^k_{\sfdd_1,\sftt'_l \sfdd_2,\sfxx}\quad & \mbox{if $\sftt'_{l'} \sfdd_2\in \sfD'_k$ and $l(\sftt'_{l'} \sfdd_2)>l(\sfdd_2)$},\\
q^{-1-\mu}\II^k_{\sfdd_1,\sftt'_{l'}\sfdd_2,\sfxx} & \\
\ \ +(q^{-1-\mu}-1)\II^k_{\sfdd_1,\sfdd_2,\sfxx}\quad &\mbox{if $\sftt'_{l'} \sfdd_2\in \sfD'_k$ and $l(\sftt'_{l'} \sfdd_2)<l(\sfdd_2)$},\\
-\II^k_{\sfdd_1,\sfdd_2,\sfxx}\quad & \mbox{otherwise.}
\end{cases}
\end{equation*}
\end{enumerate}
\end{prop}
\begin{proof} 
According to \Cref{Half way}, the proposition follows from the quadratic and braid relations of the Hecke algebras and Deodlhar's lemma stated below. 
 \end{proof}
 
\begin{lemma}[{\cite[Lemma 2.1.2]{GP}}]\label{Deodlhar's lemma}
\begin{enumerate}
    \item Let $\sfdd_1\in \sfD_k$. For $i=1,\cdots, l-1$, either $\sfss_i\sfdd_1\in \sfD_k$ or $\sfss_i \sfdd_1=\sfdd_1\sfss_{j}$ for some $1\leq j\leq l-1, j\neq l-k$. 
    \item Let $\sfdd_2\in \sfD'_k$. 
    \begin{enumerate}
        \item For $i=1,\cdots, l'-1$, either $\sfss'_i\sfdd_2\in \sfD_k$ or $\sfss'_i \sfdd_2=\sfdd_2\sfss'_{j}$ for some $1\leq j\leq l-1, j\neq k$;
        \item  Either $\sftt'_{l'}\sfdd_2\in \sfD_k$ or $\sftt'_{l'} \sfdd_2=\sfdd_2\sftt'_{l'}$. 
    \end{enumerate}
\trivial[h]{
    \begin{enumerate}
        \item If $k<l'$, then for $i=1,\cdots, l'-1$, either $\sfss'_i\sfdd_2\in \sfD_k$ or $\sfss'_i \sfdd_2=\sfdd_2\sfss'_{j}$ for some $1\leq j\leq l-1, j\neq k$, and  either $\sftt'_{l'}\sfdd_2\in \sfD_k$ or $\sftt'_{l'} \sfdd_2=\sfdd_2\sftt'_{l'}$. 
        \item If $k=l'$, then for $i=1,\cdots, l'-1$, either $\sfss'_i\sfdd_2\in \sfD_k$ or $\sfss'_i \sfdd_2=\sfdd_2\sfss'_{j}$ for some $1\leq j\leq l-1$, and $\sftt'_{l'}\sfdd_2\in \sfD_k$. 
    \end{enumerate}
}
\end{enumerate}
\end{lemma}

\subsection{Realization of $\indbasis$ in the mixed model}
 Via the mixed model $ \Weil{V_l,V'_{l'}}= \bC[\Hom(V^+_l,V'_{l'})] \otimes \Weil{V_0,V'_{l'}}$ (see \Cref{Mix model}),  we    view 
$\left(\omega_{V_{l},V'_{l'}}\right)_{\sigmaP\otimes \sigmaP'}$ as a subspace of    
\begin{equation}\label{mixmixmodelspace0}
\begin{split}
\bC[\Hom(V^+_l,V'_{l'})] \otimes   \left(\omega_{V_0,V'_{l'}}\right)_{\sigma\otimes \sigmaP'}.
\end{split}
\end{equation}
By \Cref{mix model equation}, we have the following lemma.
\begin{lemma}\label{support on a single orbit}
Assume that $f\in \Mss$ is supported on a single $ B(V_l^+)\times P'_{l'}$-orbit $\cO$ in $\Hom(V^+_l,V'_{l'})$ via \Cref{mixmixmodelspace0}. Then $f$ is determined by its value on any point $A\in \cO$.  
\end{lemma} 

Applying the mixed model  twice, we identify $\Weil{V_l,V'_{l'}}$ with 
 \begin{equation}\label{mixmixmodel}
\bC[\Hom(V^+_l,V'_{l'})] \otimes \bC[\Hom(V'^+_{l'}, V_0)] \otimes \Weil{V_0,V'_0}.
 \end{equation}
\def\AA#1{I_{\cA_{#1}}}
and view 
$\left(\omega_{V_{l},V'_{l'}}\right)_{\sigmaP\boxtimes \sigmaP'}$ as a subspace of    
\begin{equation}\label{mixmixmodelspace}
\begin{split}
\bC[\Hom(V^+_l,V'_{l'})] \otimes \bC[\Hom(V'^+_{l'}, V_0)] \otimes \left(\omega_{V_0,V'_{0}}\right)_{\sigma\otimes \sigma'}.
\end{split}
\end{equation}

For  $(\sfdd_1,\sfdd_2,\sfxx)\in \sfD_k \times \sfD'_k\times \sfS_k$, we denote by $ d_1, d_2$ and $x$ the lifts  of $\sfdd_1,\sfdd_2$ and $\sfxx$ in $\GL(V^+_l)\subset \rU(V_l)$, $\rU(V'_{l'})$ and $\GL(V'^+_{k})\subset \rU(V'_{l'})$
respectively, as fixed in \Cref{Parabolic subgroup}. From now on, we identify
  \begin{equation}\label{isomorphism V_l+}
 \Hom(V_l^+,V'_{l'})\cong V_l^-\otimes V'_{l'}    
  \end{equation}
 via \Cref{identify hom tensor}. Let
\begin{equation}\label{Akkdef}
\Akk:=\sum_{1\leq i \leq k} v_{-(l-k+i)}\otimes  v'_{i}\in \Hom(V^+_l,V'_{l'}),
\end{equation}
\[
\cO^k_{\sfdd_1,\sfdd_2,\sfxx}:=
P'_{l'}\com (d_2x \com \Akk\com  d_1^{-1})\com B(V_l^+)\subset \Hom(V_l^{+}, V'_{l'}).
\]
For $\sfww\in \sfW_{l'}$, we define 
\begin{equation}\label{bY}
\begin{split}
 V'^-_{\sfww}\coloneqq & w^{-1} V'^-_{l'} \cap V'^+_{l'},\\
 \iota(\sfww)\coloneqq &\dim V'^-_{\sfww}, \quad \text{and} \\
 \bY^-_{\sfdd_2}\coloneqq & \Hom(V'^-_{\sfww} ,V_0)\subseteq \Hom(V'^+_{l'} ,V_0)
 \end{split}
\end{equation}
similarly to \Cref{bxw} and \Cref{iotaw}. 

Let 
\begin{equation}\label{eq:Akk}
    \cAkk:=\cO^k_{1,1,1}= P'_{l'}\com   \Akk\com   B(V_l^+),
\end{equation}
and  let $I_{\cAkk}, I_0$ be  characteristic functions as in \Cref{Baby case}. 
Keeping track of the proof of \Cref{Kudla's filtrationlem} and \Cref{Half way}, 
 it follows that 
the non-zero vector in \eqref{IAk} can be set to  
\begin{equation}\label{realization1}
 \indbasisgen \coloneqq I_{\cAkk}\otimes I_0\otimes \iota_0.
\end{equation}
Here  $\iota_0 $ is a fixed non-zero vector in the one-dimensional space  $\left(\omega_{V_0,V'_{0}}\right)_{\sigma\otimes \sigma'}$, see \Cref{iota0}.
Recall the definition of $\indbasis$ by \cref{standard basis}. 
The rest of this subsection is devoted to proving the following proposition. 
\begin{prop}\label{indbasisrealzation }
 Via \eqref{mixmixmodelspace},
 $\indbasis$ is supported on $\cO_{\sfdd_1,\sfdd_2,\sfxx}$ and determined by 
 \begin{equation}\label{indbasisrealzationeq}
\indbasis ( d_2x\com \Akk
\com d_1^{-1})=\left(-q^{-\dim V_0'-1+\half \deltaan}\right)^{\iota(\sfdd_2)}I_{\bY^-_{\sfdd_2}}\otimes \iota_0,
\end{equation} 
according to  \Cref{support on a single orbit}.
\end{prop}

 \begin{remark} 
 Note that  $\cO^k_{\sfdd_1,\sfdd_2,\sfxx}\subset \cZ_k\subset \Hom(V_l^{+}, V'_{l'})$ 
as predicted by  \Cref{Kudla's filtrationlem}.
\end{remark}
  
In the proof of \Cref{indbasisrealzation }, we need the following description of $\cO_{\sfdd_1,\sfdd_2,\sfxx}$.
 




\def\pcA{\partial\cA}
\def\pcAkk{\partial\cA_{k}}

\trivial[h]{ 
A description of $\sfD_k$ and $\sfD'_k$: 
\begin{itemize}
    \item A description of $\sfD_k$: We identify a permutation  $\sigma
\in \sfS_{l}$ with a sequence $(x_1, x_2,\cdots, x_l)$ of distinct numbers in $\set{1,2,\cdots, l}$ such that $x_i = \sigma(i)$. Then the distinguished coset corresponds to the set 
\[
\set{x_1\cdots x_l| x_1< x_2 \cdots < x_{l-k}, x_{l-k+1}< x_{l-k+2}< \cdots < x_l}. 
\]
\item A description of $\sfD'_k$. We identify $x\in \sfW_{l'}$ with a sequence of elements $x_1\cdots x_{l'}$ in $\set{\pm 1, \cdots, \pm l'}$. Such that $x_i = w(i)$.  Let $J_k = \Sigma - \set{\alpha_k}$ where $\Sigma$ is the set of simple root,
$\alpha_k$ is the $k$-th simple root. 
The distinguished right coset representatives $D'_k$ of
$\sfW_{l'}/\sfS_k\times \sfW_{l'-k}$
can be explicitly described
\[
\begin{split}
   D'_k &= \Set{w\in \sfW_{l'}| w\alpha >0 \text{ for all }  \alpha\in J_k}\\
   & =\Set{x_1\cdots x_r x_{r+1}\cdots x_{k}x_{k+1}\cdots x_{l'}|
   \begin{array}{l}0<x_1 < \cdots < x_r\\ x_{r+1} < \cdots < x_{k}<0\\ 0<x_{k+1} < \cdots < x_{l'}. 
   \end{array}
   }  
\end{split}
\]
\end{itemize}
}

\begin{lemma}\label{geometry action map}
The following map is an injection with image $ \cO^k_{\sfdd_1,\sfdd_2,\sfxx}$:
\[
\begin{split}
 U'_{\sfdd_2\sfxx}\times U_{\sfdd_1}\times \cAkk&\longrightarrow   \Hom(V_l^{+}, V'_{l'}) ,\\
 (u_2, u_1, A)& \mapsto u_2d_2x\com A \com (u_1d_1)^{-1}.
\end{split}
\]
Here $d_1,d_2$ and $x$ are lifts of $\sfdd_1,\sfdd_2$ and $\sfxx$.  
\end{lemma}
\begin{proof} The proof is routine by using Bruhat decomposition and the fact that $\sfD_k$ ad $\sfD'_k$ are the set of distinguished representatives. We omit the details. 
\end{proof}
\trivial[h]{
\zjl{Hide this proof. Use one sentence to illustrate the idea}
We fix some notations that will only be used in this  lemma.
Let  
\[
\alpha_k: P(V_{l-k}^+,V_l^+)\longrightarrow \GL(\widehat{V}^+_k)\cong \GL_k(F),\quad
\beta_k: P(V'^+_{k},V'_{l'})\longrightarrow \GL(V'^+_{k})\cong \GL_k(F)
\]
be the  two natural projections, where $P(V_{l-k}^+,V_l^+)$ and $P(V'^+_{k},V'_{l'})$ are defined in \Cref{Kudla's filtrationlem}, and   the isomorphisms are in \Cref{general linear group isomorphism}.  Let 
\begin{equation*}\label{stablizer Akk}
 R_k:=\Stab_{\GL(V_l^{+})\times \rU(V'_{l'})}{\Akk}=\set{(p,p')\in P(V_{l-k}^+,V_l^+)\times P(V'^+_{k},V'_{l'})| \alpha_k(p)=\beta_k(p')}. 
\end{equation*}
 \begin{equation*}\label{E-k}
E_k:=\Stab_{B(V_l^+)\times P'_{l'}}(\Akk)=R_k\cap (B(V_l^+)\times P'_{l'})=
\set{(p,p')\in B(V_l^+)\times P'_{l'}|\alpha_k(p)=\beta_k(p') }.
\end{equation*}

  Since $\sfdd_1$ is the distinguished representative of $\sfS_l\backslash \sfS_{l-k}\times \sfS_{k}$, $\sfdd_1\sfss_i>0$ for $1\leq i\leq l-k-1$ and $l-k+1\leq i\neq l-1$. Then by the definition \Cref{definition Uw},
$
U_{\sfdd_1^{-1}} $ is contained in the nilpotent of  $P(V_{l-k}^+,V_l^+)$,
which is contained in $ R_k$.  In particular, $
\Akk\cdot U_{\sfdd_1^{-1}}= \Akk.  
$
So
\[
\begin{split}(U'_{\sfdd_2\sfxx} d_2x  )\cdot \cAkk\cdot ( d_1^{-1} U_{\sfdd_1})=&(U'_{\sfdd_2\sfxx} d_2x P'_{l'})\cdot \Akk\cdot (B(V_l^+) d_1^{-1} U_{\sfdd_1})\\
 =& ( P'_{l'}d_2x P'_{l'} )\cdot \Akk\cdot (B(V_l^+) d_1^{-1} B(V_l^+) ) \\
 = & ( P'_{l'}d_2x ) \cdot \Akk\cdot (B(V_l^+) d_1^{-1} B(V_l^+)) \quad  \mbox{(by \Cref{eq:cAk} )}    \\
= & (P'_{l'}d_2x) \cdot \Akk\cdot U_{d_1^{-1}}d_1^{-1} B(V_l^+)=  \cO^k_{\sfdd_1,\sfdd_2,\sfxx}.\\
\end{split}   
\]
Thus  the  map in the lemma has image $ \cO^k_{\sfdd_1,\sfdd_2,\sfxx}$.
Now we prove its injectivity.
Consider
the map 
\begin{equation}\label{fiberek} 
\begin{split}
   (U'_{\sfdd_2\sfxx} d_2x P'_{l'})\times (U_{\sfdd_1}d_1 B(V_l^+))&\longrightarrow \cO^k_{\sfdd_1,\sfdd_2,\sfxx},\\ 
   (u_2 d_2x p', u_1d_1 b )& \mapsto u_2d_2xp'\cdot \Akk \cdot (u_1d_1b)^{-1}.
\end{split}
\end{equation}
which 
is $P'_{l'}\times B(V_l^+)$-equivariant and surjective.
So by a counting argument, we only need to show that the preimage of one  point (and thus every point by the $P'_{l'}\times B(V_l^+)$-equivariance of \eqref{fiberek})
has cardinality   $|E_k|$.
 Consider the point  $d_2x\cdot \Akk \cdot d_1^{-1}$. If $u_2d_2xp'\cdot \Akk \cdot (u_1d_1b)^{-1}=d_2x\cdot \Akk \cdot d_1^{-1}$, then $(d_1^{-1} u_1 d_1 b )\times (\dot (d_2x)^{-1}u_2d_2xp')\in R_k $.  We show that $u_1=u_2=1$ and $(b,p')\in E_k$ in three steps. 

First,  consider the projection of $R_k$ to $P(V^+_{l-k},V^+_{l})$, we have $d_1^{-1} u_1 d_1 b \in P(V^+_{l-k},V^+_{l})$, which implies $d_1^{-1} u_1 d_1\in P(V^+_{l-k},V^+_{l})$. On the other hand,
 It follows from the definition of $U_{\sfdd_1}$ in \Cref{definition Uw} that
\[
U_{\sfdd_1}=\prod_{\alpha\in \Sigma^+_{I}, \sfdd_1^{-1}\alpha<0}U_{\alpha}
\]
and \qcl{I hideded $\Sigma^+_{I}$ }
\[
d_1^{-1}U_{\sfdd_1}d_1=\prod_{\alpha\in \Sigma^+_{I}, \sfdd_1^{-1}\alpha<0}U_{\sfdd_1^{-1}\alpha}.
\]
\zjl{revise here} 
We have $d_1^{-1}U_{\sfdd_1}d_1\in \GL(V_l^+)$. 
Moreover, since $d_1$ is a distinguished representative of We have $d_1^{-1}U_{\sfdd_1}d_1\in \GL(V_l^+)$, we know that $\sfdd_1^{-1}\alpha$ is not the linear span of $-\sfss_i$ for $i=\set{1,\dots, l}\backslash \set{l-k}$. Therefore, we conclude that 
\begin{equation}\label{d1conjugate}
   (d_1^{-1}U_{\sfdd_1} d_1 )\cap P(V^+_{l-k},V^+_{l}) =1.  
\end{equation}

Second, 
consider the projection from $R_k$ to $ P(V'^+_{k},V'_{l'})$, we have $( d_2x)^{-1}u_2d_2xp'\in P(V'^+_{k},V'_{l'})$. Since $x,p'\in P(V'^+_{k},V'_{l'})$, we deduce that $ d_2^{-1}u_2d_2\in P(V'^+_{k},V'_{l'})$. Note that $\sfdd_2$ is a distinguished representative of $\sfW_{l'}/(\sfS_k\times \sfW_{l'-k})$, we have 
$l(\sfdd_2\sfxx)=l(\sfdd_2)+l(\sfxx)$ and  $U'_{\sfdd_2\sfxx}=U'_{\sfdd_2}(U'_{\sfxx})^{d_2}$. 
Write $u_2=u_{2,d_2}d_2u_{2,x}d_2^{-1}$ for $u_{2,d_2}\in U'_{\sfdd_2}$ and $u_{2,x}\in U'_{\sfxx}$. Then 
\[
d_2^{-1} u_2 d_2 =(d_2^{-1}u_{2,d_2} d_2x )u_{2,x}\in P(V'^+_{k},V'_{l}).
\]
Note that $u_{2,x}\in U'_{\sfxx}\subseteq  P(V'^+_{k},V'_{l'})$. Therefore we have $d_2^{-1} u_2 d_2 \in P(V'^+_{k},V'_{l'})$. On the other hand, by a similar argument to \Cref{d1conjugate},  
\trivial[h]{
It follows from the definition of $U'_{\sfdd_2}$ in \Cref{definition Uw} that 
\[
U'_{\sfdd_2}=\prod_{\alpha\in \Sigma'^+\backslash \Sigma'^+_{I}, \sfdd_1^{-1}\alpha<0}U'_{\alpha}.
\]
Then 
\[
x^{-1}d_2^{-1}U_{\sfdd_2}d_2x=\prod_{\alpha\in \Sigma'^+\backslash \Sigma'^+_{I}, \sfdd_2^{-1}\alpha<0}U_{\sfxx^{-1}\sfdd_2^{-1}\alpha}.
\]
Moreover, since $d_2$ is a distinguished representative, we know that 
\[
\sfdd_2^{-1}\alpha \notin \Sigma_J^-.
\]
(Here we use $\sfdd_2\alpha>0$ for $\alpha\in \Sigma_J^+$.
}
We know  that 
\[
d_2^{-1}U'_{\sfdd_2}d_2 \cap P(V'^+_{k},V'_{l'})=1.
\]
Then we conclude that $u_{2,d_2}=1$. 

Finally, by the above two steps, we know that $(d_1^{-1} u_1 d_1 b )\times (\dot (d_2x)^{-1}u_2d_2xp')= b\times (x^{-1}u_{2,x}x p')\in R_k $. We use the condition 
\[
\alpha_k(b)=
\beta_k(x^{-1}u_{2,x}xp')=x^{-1}u_{2,x}x\beta_k(p').
\]
in  \Cref{stablizer Akk}. This implies that $x^{-1}u_{2,x}x\in B_k$. On the other hand, by \Cref{definition Uw}, we know that 
\[
x^{-1}U'_{\sfxx}x\cap B_k= 1.
\]
Then we can conclude that $u_{2,x}=1$ and $(b,p')\in E_k$. 
 }
 \begin{proof}[Proof of  \Cref{indbasisrealzation }]
The key is to compute $T_{\sfdd_1}T'_{\sfdd_2}T'_{\sfxx} \cdot \indbasisgen$.  Since $\sfdd_2\in \sfD'_k$,    $l(\sfdd_2\sfxx)=l(\sfdd_2)+l(\sfxx)$ and $d_2 x$ is the lift of $\sfdd_2\sfxx$, it follows from the mixed model formula  \Cref{mix model equation} and  \Cref{actformula} that  
\begin{equation}\label{action equation}
\begin{split}
 T_{\sfdd_1}T'_{\sfdd_2}T'_{\sfxx}  \cdot \indbasisgen &= T_{\sfdd_1} T'_{\sfdd_2 \sfxx} \cdot (I_{\cAkk}\otimes I_0\otimes \iota_0)  \\&= \sum_{u_2\in U'_{\sfdd_2\sfxx}}\sum_{u_1\in U_{\sfdd_1}}\sum_{A\in \cAkk}\Weil{V_l,V'_{l'}}(u_1d_1u_2d_2x)\left(I_{A}\otimes I_0\otimes \iota_0 \right)\\ 
 & =\sum_{u_2\in U'_{\sfdd_2\sfxx}}\sum_{u_1\in U_{\sfdd_1}}\sum_{A\in \cAkk} I_{u_2d_2x\cdot
  A\cdot (u_1d_1)^{-1}}\otimes \Weil{V_0,V'_{l'}}(u_2d_2x)\left(I_0\otimes \iota_0 \right).
  \end{split}
\end{equation}
 Then by \Cref{geometry action map} and \Cref{Partial Fourier d},
 $T_{\sfdd_1}T'_{\sfdd_2}T'_{\sfxx} \cdot \indbasisgen $ is supported on $\cO^k_{\sfdd_1,\sfdd_2,\sfxx}$ and \begin{equation*}
 \begin{split}
 \left(T_{\sfdd_1}T'_{\sfdd_2}T'_{\sfxx}\cdot \indbasisgen\right) ( d_2x\com \Akk
\com d_1^{-1}) & = \Weil{V_0,V'_{l'}}(d_2x)\left(I_0\otimes \iota_0 \right) \\
& =\left(\gamma_{V_0}q^{-\half \dim V_0}\right)^{\iota(\sfdd_2)}I_{\bY^-_{\sfdd_2}}\otimes \iota_0.
\end{split}
\end{equation*} 
Here we also use the fact that 
$\iota(\sfdd_2\sfxx)=\iota(\sfdd_2)$ and $\bY^{-}_{\sfdd_2\sfxx}=\bY^-_{\sfdd_2}$.
Taking care of the normalizations 
\begin{equation*}
 \cT_{\sfdd_1}=T_{\sfdd_1},\quad \cT'_{\sfdd_2}= \left( -\gamma_{V_0}^{-1}q^{-\dim V_0'+\half\dim V_0-1+\half \deltaan} \right)^{\iota(\sfdd_2)}T'_{\sfdd_2}, \quad   \cT'_{\sfxx}=T'_{\sfxx}
\end{equation*}
given by   \Cref{nomralized operator} and \Cref{Normalized t_k23}, we deduce the  proposition from \eqref{standard basis}. 
\end{proof}
\trivial[h]{
For each $(\sfdd_1,\sfdd_2,\sfxx)\in \sfD_k \times \sfD'_k \times \sfS_k$, we define $\II_{\cO^k_{\sfdd_1,\sfdd_2,\sfxx}}$ to be the unique element in  $\Mss$ such that 
\begin{itemize}
\item $\II_{\cO^k_{\sfdd_1,\sfdd_2,\sfxx}}$ as a function on $\Hom(V^+_{l},V'_{l'})$ is supported on
the $\cO^k_{\sfdd_1,\sfdd_2,\sfxx}$;
\item $\II_{\cO^k_{\sfdd_1,\sfdd_2,\sfxx}}(d_2\cdot \Akk\cdot (\dot d_1 x)^{-1}) = I_{\bX^{-}_{d_2}}\otimes \iota_0$, where 
$\bX^-_{d_2}$ is  defined in \Cref{bxw}.
\end{itemize}
The following proposition follows from \Cref{equation no name}.
\begin{prop}\label{transfer of basis}
We have 
\[
\indbasis  =\left(-q^{-\dim V_0'-1+\half \deltaan}\right)^{\iota(d_2)} \II_{\cO^k_{\sfdd_1,\sfdd_2,\sfxx}}.
\]
\end{prop}
}


\subsection{The action of $\cT_{\sftt_{l}}$}
For $(\sfdd_1,\sfdd_2,\sfxx)\in \sfD_k\times \sfD'_k \times \sfS_k$, 
by \Cref{standard basis} and \Cref{realization1}, we have
\begin{equation}\label{Hecke action l}
\cT_{\sftt_l}  \cdot \indbasis 
 =  \cT_{\sftt_l}  \left( \cT_{\sfdd_1}\cT'_{\sfdd_2} \cT'_{\sfxx}\cdot \indbasisgen \right)
 = \cT'_{\sfdd_2}\cT'_{\sfxx} (\cT_{\sftt_l} \cT_{\sfdd_1}\cdot \indbasisgen).
\end{equation}
We compute $\cT_{\sftt_l} \cT_{\sfdd_1}$ as follows. Note that we have the following double coset decomposition 
\[
 \sfS_{l}= \coprod_{i=1,2} (\sfS_{l-1}\times \sfS_1 )\sfww_i (\sfS_{l-k}\times \sfS_{k}), 
\]
where $\sfww_1:=1$ and $\sfww_2 :=\sfss_{l-1}\cdots \sfss_{l-k}$ are two distinguished double coset representatives. For each $i$, let $\sfD^i_{k}$ be the set of distinguished left coset representatives of $\sfww_i (\sfS_{l-k}\times \sfS_{k})\sfww_i^{-1} \cap (\sfS_{l-1}\times \sfS_1)$ in $\sfS_{l-1}\times \sfS_1$.
\begin{lemma}\label{distigushe double}
We have $
\sfD_{k}=\coprod\limits_{i=1,2}\sfD^i_{k}\sfww_i.
$
In particular, each $\sfdd_1\in \sfD_{k}$ has a unique decomposition $\sfdd_1=\sfyy_i \sfww_i$, where $i\in \set{1,2}, \sfyy_i\in \sfD^i_{k}$ and $l(\sfdd_1)=l(\sfyy_i)+l(\sfww_i)$.
\end{lemma}
\begin{proof}
See \cite[Lemma 2.1.9]{GP}.
\end{proof}
\begin{enumerate}
    \item 
Suppose $\sfdd_1\in \sfD^1_{k}\sfww_1\subseteq \sfS_{l-1}\times \sfS_1$. 
Then 
$\sfdd_1 \sftt_l = \sftt_l \sfdd_1$ and $l(\sftt_l \sfdd_1 ) = l(\sfdd_1 \sftt_l) = l(\sfdd_1) + 1$.
Therefore, 
\begin{equation}\label{d1 in sfS_{l-1}}
\cT_{\sftt_l} \cT_{\sfdd_1}= \cT_{\sfdd_1} \cT_{\sftt_l}. 
\end{equation}
\item Suppose $\sfdd_1\in \sfD^2_{k}\sfww_2$. Write $\sfdd_1 = \sfyy \sfww_2$ with $\sfyy\in \sfD^2_{k}\subseteq \sfS_{l-1}\times \sfS_1$. Then $l(\sfdd_1)=l(\sfyy)+l(\sfww_2)$ and $\cT_{\sfdd_1}=\cT_{\sfyy}\cT_{\sfww_2}$. We have 
\begin{equation*}
\begin{split}
\cT_{\sftt_l}\cT_{\sfdd_1}= \cT_{\sftt_l} \cT_{\sfyy}\cT_{\sfww_2}=  \cT_{\sfyy}\cT_{\sftt_l}\cT_{\sfww_2} = \cT_{\sfyy} \cT_{\sfww_2^{-1}}^{-1} (\cT_{\sfww_2^{-1}} \cT_{\sftt_l}\cT_{\sfww_2}).
\end{split}
\end{equation*}
 Note that $\sfww_2^{-1} \sftt_l \sfww_2 = \sftt_{l-k},$ and $l(t_{l-k}) = 2l(\sfww_2)+1$. Hence $\cT_{\sfww_2^{-1}}\cT_{\sftt_l}\cT_{\sfww_2}=\cT_{\sfww_2^{-1} \sftt_l \sfww_2}=\cT_{\sftt_{l-k}}$ and we deduce
\begin{equation}\label{d1 not in sfS_{l-1}}
\begin{split}
\cT_{\sftt_l}\cT_{\sfdd_1}= \cT_{\sfyy} \cT_{\sfww_2^{-1}}^{-1} \cT_{\sftt_{l-k}}.
\end{split}
\end{equation}
\end{enumerate}

 By \Cref{Hecke action l}, \Cref{d1 in sfS_{l-1}}, \Cref{d1 not in sfS_{l-1}} and  \Cref{Hecke action half}, the computation of $\cT_{\sftt_l}  \cdot \indbasis$ can be reduced to the computations of $\cT_{\sftt_l}\cdot \indbasisgen$ and $\cT_{\sftt_{l-k}}\cdot \indbasisgen$. 
We present the formulas, whose 
proofs will be given in \Cref{appendix A} and \Cref{appendix B}. 
For each $1\leq i\leq j\leq l$, we define $\sfss_{i,j}\in \sfW_l$ by 
\begin{equation}\label{sij}
 \sfss_{i,j}: =\sfss_{j-1}\cdots \sfss_{i}.  
\end{equation}
Here $\sfss_{i,i}=1$ by convention. The notation naturally extends to $\sfW_{l'}$. 
\begin{prop}\label{lem:TlAk}
We have 
\begin{equation*}
\begin{split}
& \cT_{\sftt_{l}}\cdot \indbasisgen \\
& = -q^{k-l'+\mu} \II^k_{1,\sftt'_k,1} \\
& \ \ -(1-q)q^{k-l'} \Big(q^{\mu}\sum_{i=k+1}^{l'} \II^k_{1,\sftt'_i\sfss'_{k,i},1}  - q^{-1} \II^{k-1}_{\sfss_{l-k+1,l},1,1}- q^{-1}\sum_{i=k}^{l'} \II^k_{1,\sfss'_{k,i},1 }\Big),\\
&  \cT_{\sftt_{l-k}} \cdot \indbasisgen \\ 
 &= -q^{2k-l'} \indbasisgen \\
 &\ \ - q^{k-l'}\Big(  \sum_{i=k+1}^{l'}\II^{k+1}_{1,\sfss'_{k,i},\sfss'_{1,k}} -q^{\mu+1}\sum_{i=k+1}^{l'}\II^{k+1}_{1,\sftt'_{i}\sfss'_{k,i},\sfss'_{1,k}} 
 + (q-1)\sum_{i=1}^k  q^{k-i}\II^k_{\sfss_{l-k,l-k+i},1,\sfss'_{1,i}} \Big).\\ 
\end{split}
\end{equation*}
\end{prop}

Combining \Cref{d1 in sfS_{l-1}}, \Cref{d1 not in sfS_{l-1}} and \Cref{lem:TlAk}, we deduce the following corollary. 
\begin{cor}\label{action of Tt_l in general}
For each $(\sfdd_1,\sfdd_2,\sfxx)\in \sfD_k \times \sfD'_k\times \sfS_k$,
there are Laurent polynomials $h^{\sfdd'_1,\sfdd'_2,\sfxx'}_{\sfdd_1,\sfdd_2,\sfxx}$ with coefficients in $\bZ$ and in the indeterminate $\nu^{\half}$, i.e., $h^{\sfdd'_1,\sfdd'_2,\sfxx'}_{\sfdd_1,\sfdd_2,\sfxx}(\nu)\in \bZ[\nu^{\half},\nu^{-\half}]$, such that  
\[
\cT_{\sftt_l} \cdot \indbasis  =\sum_{(\sfdd'_1,\sfdd'_2,\sfxx')} h^{\sfdd'_1,\sfdd'_2,\sfxx'}_{\sfdd_1,\sfdd_2,\sfxx}(q)\, \II^i_{\sfdd'_1,\sfdd'_2,\sfxx'},
\]
where $(\sfdd'_1,\sfdd'_2,\sfxx')$ 
runs over $\bigsqcup \limits_{i=0}^ { \min{(l,l')}}\sfD_i \times \sfD'_i\times \sfS_i$.
\end{cor}

\section{The generic Hecke algebra module  and deformation}
In this section, we prove  \Cref{thm:main} by explicitly constructing the  generic Hecke algebra module $\sfM$, using the explicit formulas in the last section.

\subsection{The construction of  $\sfM$}\label{thm:mainproof}

Retrieve the set-up in \ref{mainsec}.
Let $\sfH^\diamond$ be the subalgebra of $\sfH:=\sfH_{l,\mu}$ generated by $\set{\sfT_{\sfss_i}| i\in {1,\cdots, l-1}}$ so that $\sfH$ is generated by $\sfH^\diamond$ and $\sfT_{\sftt_l}$.
Then $\sfH^\diamond$ is isomorphic to a generic Hecke algebra of type $A_l$ over $R$ with parameter $\nu$,  and the specialization of $\sfH^\diamond$ at $\nu=q$ is isomorphic to $\cH_{V_l^+}$.
For $1\leq k\leq\min\set{l,l'}$, Let $\sfM_k$ be generic version of the induced module $\cF_k$ in  \Cref {induced} with $\cH'_{l'}$ replaced by $\sfH'$ and so on.
 Explicitly, $\sfM_k$ is a free $R$-module with basis \begin{equation}\label{genekbasis}
\{\mbox{symbol } \sfI^k_{\sfdd_1,\sfdd_2,\sfxx}  |  (\sfdd_1,\sfdd_2,\sfxx)\in \sfD_k\times \sfD'_k\times \sfS_k\}, 
\end{equation}
where the action of $\sfH^\diamond\times \sfH'$  is given by the formulas in \Cref{Hecke action half} with ``$\cT$'', ``$\II$'' and ``$q$'' replaced by ``$\sfT$'', $\sfI$'' and  ``$\nu$'' respectively. 
By the compatibility of specialization with induction (see \cite[\S 9.1.5]{GP}), we know that the specialization of $\sfM_k$ at $\nu=q$ is identified with $\cF_k$ as an $\cH_{V^+_l} \times \cH'_{l'}$-module.

The module $\sfM_k$ leads to an $R$-algebra homomorphism:
   \[\Psi_k: \sfH^\diamond\times \sfH'\to \End_R (\sfM_k).
   \]
 Let \[
\sfM:=\bigoplus_{k=1}^{\min\set{l,l'}} \sfM_k, \quad \text{and}\quad \Psi:=\bigoplus_{k=1}^{\min\set{l,l'}} \Psi_k.
\]
We now extend $\Psi$ to $\sfH\times \sfH'$ by defining the image of the $\sfT_{\sftt_l}$.
Recall there are the Laurent polynomials $h^{\sfdd'_1,\sfdd'_2,\sfxx'}_{\sfdd_1,\sfdd_2,\sfxx}\in R $ determined by \Cref{action of Tt_l in general}  for all  $({\sfdd'_1,\sfdd'_2,\sfxx'),(\sfdd_1,\sfdd_2,\sfxx})\in \coprod\limits_{k=1}^{  \min(l,l')}\sfD_k\times \sfD'_k\times \sfS_k$.   
Let $\sT_{\sftt_l}\in \End_R (\sfM)$ be the matrix  with entries  $h^{\sfdd'_1,\sfdd'_2,\sfxx'}_{\sfdd_1,\sfdd_2,\sfxx}$ under the basis \eqref{genekbasis}.   
\trivial[h]{
\begin{thm}
\begin{enumerate}
    \item $\sfM$ is a representation \qcl{module} of the generic Hecke algebra $\sfH\times \sfH'$.
    \item \qcl{remove} The specialization $\sfM_{\nu=q}$ of $\sfM$ at $\nu = q$ is isomorphic to the natural $\cH_l\otimes \cH_{l'}$-module $\Mss$. 
\end{enumerate}
\end{thm}

We first prove (1). A prior, we see that $\sfM$ is a module of the generic Hecke subalgebra $\sfH^\diamond\times \sfH'$ since the specialization is compatible with the parabolic inductions. It remains to check that the action satisfies the quadratic relation 
\begin{equation}\label{equation 1}
   (\cT_{\sftt_l}+1)(\cT_{\sftt_l}-q^{\mu})=0 
\end{equation}
and the braid relation
\begin{equation}\label{equation 2}
  (\cT_{\sftt_l}\cT_{\sfss_{l-1}})^4=1.
\end{equation}
We prove \Cref{equation 1} here and the proof of \Cref{equation 2} is similar. Express the action of $\cT_{\sftt_l}$ as matrix. Then we know that $ (\cT_{\sftt_l}+1)(\cT_{\sftt_l}-q^{\mu})$ is a matrix with each entry a polynomial in $\bZ[\nu^\half, \nu^{-\half}]$. Moreover, we know that the specialization of each entry at $\nu=q$ is $0$. Since we can find infinite many blabla (\zjl{write it properly}) and a non-zero polynomial can only have finite many zeros. We conclude that $ (\cT_{\sftt_l}+1)(\cT_{\sftt_l}-q^{\mu})=0$. (2) is clear by our construction. 
}
We define 
\begin{equation}\label{eq:Psi_Ttl}
\Psi(\sfT_{\sftt_l}):= \sT_{\sftt_l}. 
\end{equation}
\begin{prop} Under the definition \eqref{eq:Psi_Ttl},  
$\Psi$ extends to an $R$-algebra homomorphism:
 \[
 \sfH\times \sfH'\to \End_R (\sfM).
 \]
\end{prop}

\begin{proof}  Since the Hecke algebras can be equivalently defined using the quadratic relations and the braid relations for generators,  it is enough to prove $(\sT_{\sftt_l}+1)(\sT_{\sftt_l}-\nu^{\mu})=0$,
  $ (\sT_{\sftt_l}\Psi (\sfT_{\sfss_{l-1}}))^4=1$
  and 
  $ (\sT_{\sftt_l}\Psi (\sfT_{\sfss_{i}}))^2=1$ ($i\in \set{1,\cdots, l-2}$) in $\End_R (\sfM)$. 
These equations form a system of finite many equations of Laurent polynomials in the indeterminate $\nu^\half$. 
By \Cref{abundance} and \Cref{action of Tt_l in general},  the equations hold at $\nu=q^n$ for every positive integer $n$. 
Therefore, the equations hold.    
 \end{proof}

\subsection{Proof of \Cref {thm:main}} 
 The desired property (a) on the specialization of $\sfM$ at $\nu=q$ follows from the construction. Now we study the specialization of $\sfM$ at $\nu=1$.
\begin{lemma}\label{lemat1} 
For $1\leq k\leq\min\set{l,l'}$, we have
\[
\left(\sfT_{\sftt_{l}}\cdot\sfI^k_{1,1,1} \right)_{\nu=1}=
-\left(\sfT_{\sftt_{k}'}\cdot\sfI^k_{1,1,1}\right)_{\nu=1},
\]
\[
\left(\sfT_{\sftt_{l-k}}\cdot\sfI^k_{1,1,1}\right)_{\nu=1}=
-\left( \sfI^k_{1,1,1}\right)_{\nu=1} \mod \bigoplus_{i>k}\sfM_{i,\nu=1}.
\] 
\end{lemma}
\begin{proof}
These follow from \Cref{Hecke action half} and 
\Cref{lem:TlAk}. 
 \end{proof}
By the construction of $\sfM_k$, as 
an   $\sfH^\diamond_{\nu=1}\otimes \sfH'_{\nu=1} =\bC[\sfS_l]\times\bC[\sfW_{l'}]$-module,
$\sfM_{k,\nu=1}$ is naturally identified with 
 $\Ind_{\sfS_{l-k}\times \triangle \sfS_{k}\times \sfW_{l'-k}}^{\sfS_{l}\times \sfW_{l'}}\left(
 \bbone_{l-k}  \boxtimes \bbone_k \otimes \varepsilon_{l'-k}\right)$, which is
the restriction of the desired module. 
\Cref{lemat1} implies that 
  \[
   \sfM_{\geq k,\nu=1}/\sfM_{\geq {k+1},\nu=1}\cong
   \Ind_{\sfW_{l-k}\times \triangle \sfW_{k}\times \sfW_{l'-k}}^{\sfW_{l}\times \sfW_{l'}}\left(
 \varepsilon_{l-k} \boxtimes \varepsilon_k \boxtimes\varepsilon_{l'-k}\right).\]
Since $\sfM_{\nu=1}$ is a semisimple  $\bC[\sfW_l]\times\bC[\sfW_{l'}]$-module,
property (b) follows.
 \begin{remark}
 We emphasis that $\bigoplus_{i\geq k}\sfM_{i}$ is not stable under the action of $\sfH\times \sfH'$ since 
$\bigoplus_{i\geq k}\sfM_{i,\nu=q}=\bigoplus_{i\geq k}\cF_{i}$ coming from Kudla's filtrations is not stable under the action of $\cH_l\otimes \cH'_{l'}$ (see  \Cref{lem:TlAk}). 
 
 \end{remark}

\subsection{Proof of \Cref{cor:AMR}}\label{proof of theorem 2}
\def\tauL{\tau_L}
\def\tauR{\tau_R}
\def\sigmaL{\sigma_L}
\def\sigmaR{\sigma_R}

Retrieve the notations in \Cref{result interplation}. 
Note that 
\[
\Ind_{\triangle \sfW_{k} }^{\sfW_{k}\times \sfW_{k}}
(\varepsilon_k)=\sum_{\pi\in \Irr(\sfW_k)} \pi\boxtimes (\pi^\vee\otimes \varepsilon_k)
\]
and $\pi\cong \pi^\vee$ for any representation of $\sfW_k$. We have 
\[
\begin{split}
\sfM_{\nu=1}\cong 
&\sum_{k=0}^{\min\set{l,l'}}
 \Ind_{\sfW_{l-k}\times \triangle \sfW_{k}\times \sfW_{l'-k}}^{\sfW_{l}\times \sfW_{l'}}
\left( \varepsilon_{l-k}\boxtimes\varepsilon_k \boxtimes\varepsilon_{l'-k}\right)\\
\cong & \sum_{k=0}^{\min\set{l,l'}}
\sum_{\pi\in \Irr(\sfW_k)}
\Ind_{\sfW_{l-k}\times \sfW_{k}}^{\sfW_{l}}
(\varepsilon_{l-k}\boxtimes\pi ) \boxtimes \Ind_{\sfW_{k}\times  \sfW_{l'-k}}^{\sfW_{l'}}
((\pi\otimes \varepsilon_k) \boxtimes\varepsilon_{l'-k}))\\
\cong& \sum_{k=0}^{\min\set{l,l'}} \sum_{\pi\in \Irr(\sfW_k)}
\Ind_{\sfW_{l-k}\times \sfW_{k}}^{\sfW_{l}}
(\varepsilon_{l-k} \boxtimes \pi) \boxtimes \left(\Ind_{\sfW_{k}\times  \sfW_{l'-k}}^{\sfW_{l'}}
(\pi \boxtimes\bbone_{l'-k}) \otimes \varepsilon_{l'}\right),
\end{split}
\]
 where $\bbone_{l'-k}$ is the trivial character of $\sfW_{l'-k}$.
\begin{lemma}\label{homomorphism}
 Assume that $\gamma\times \eta \in \Irr(\sfW_k)$. Then
 \begin{equation*}\label{Ring morphism}
 \begin{split}
 \Ind_{\sfW_{l-k}\times \sfW_{k}}^{\sfW_{l}}( \varepsilon_{l-k}\boxtimes (\gamma\times \eta))&=\gamma\times \sfX_{l-k}(\eta),\\
  \Ind_{\sfW_k\times \sfW_{l'-k}}^{\sfW_{l'}}((\gamma\times \eta)\boxtimes \bbone_{l'-k})&=\sfX_{l'-k}(\gamma)\times \eta   
  \end{split}
 \end{equation*}
and 
\begin{equation}\label{flip}
 (\gamma\times \eta)\otimes \varepsilon_{k}\cong \eta \times \gamma.
\end{equation}
\end{lemma}
\begin{proof}
See \cite[Proposition 3.5]{AMR}.
\trivial[h]{
Note that $\bbone_b=\bfone_b\times \bfone_0$ and $\varepsilon_b=\bfone_0\times \bfone_b$. By \cite{Z}, the map \Cref{Hyperocta} extends to a natural ring isomorphism $\cS\otimes\cS \rightarrow \cW$. Thus we deduce \Cref{Ring morphism}.  \Cref{flip} follows by a direct computation. 
}
\end{proof}
Now 
\begin{equation*}
\begin{split}
\sfM_{\nu=1}
\cong &\sum_{k=0}^{\min\set{l,l'}} \sum_{\gamma \times \eta \in \Irr(\sfW_k)}
\Ind_{\sfW_{l-k}\times \sfW_{k}}^{\sfW_{l}}
( \varepsilon_{l-k} \otimes (\gamma\times \eta) ) \\
 & \ \hspace{9em} \boxtimes \left(\Ind_{\sfW_{k}\times  \sfW_{l'-k}}^{\sfW_{l'}}
((\gamma\times \eta) \otimes\bbone_{l'-k}) \otimes \varepsilon_{l'}\right)\\
\cong &\sum_{k=0}^{\min\set{l,l'}} \sum_{\gamma \times \eta \in \Irr(\sfW_k)} (\gamma\times \sfX_{l-k}(\eta))
 \boxtimes (\eta \times \sfX_{l'-k}(\gamma)).
\end{split}
\end{equation*}
Note that $\alpha\times \beta$ is an irreducible component of $\gamma\times \sfX_{l-k}(\eta)$ if and only if $\gamma=\alpha$ and 
\[
\inn{\beta}{\sfX_{l-k}(\eta)}=\inn{\sfX_{l-k}^*(\beta)}{\eta}\neq 0,\]
which must equal to $1$ by Pieri's rule (cf. \cite[4.44]{Fu}). Then the theorem follows.
\trivial[h]{
Pieri's rule says that $ \sfX_{l-k}(\eta)$ is multiplicity free. 
}

 
\def\rr{r_1}
\subsection{Conservation relation: general cases}\label{thm:conservationproof}
Recall that we proved \Cref{thm:conservation} for theta cuspidal representations in \Cref{Conservation}.
In this subsection, we prove \Cref{thm:conservation} in the general cases.

The following lemma is easy to verify and will be used in the proof of \Cref{thm:conservation}.  
\begin{lemma}\label{pull back iota}
We have the following commutative diagram 
\begin{equation*}
\begin{tikzcd}
\Rep(\sfH_{l,\mu}) \ar[r, "\kappa^*"] \ar[d,"\otimes_R \bC_1"]
& \Rep(\sfH_{l,-\mu}) \ar[d,"\otimes_R \bC_1"]\\
\Rep(\sfW_l)\ar[r, "\otimes \varepsilon_l"] & \Rep(\sfW_l). 
\end{tikzcd}
\end{equation*}
 Here $\kappa$ is defined in \Cref{Different normalization} and 
$\kappa^*$ denotes the pull-back via $\kappa$. 
\end{lemma}
Retrieve the notations in \Cref{theta cuspidal} and \Cref{result interplation}.
For each $\gamma\in \Irr(\sfS_k)$, we define 
\begin{equation}\label{pieri rule}
 \rr (\gamma) := \max\Set{\text{non-negative integer } i |  \sfX_{i}^* (\gamma)\neq 0}.   
\end{equation}

\def\cc{c}
\begin{proof}[Proof of \Cref{thm:conservation}]  
We may assume  $\pi\in \cE(\rU(V_{l}),\sigma)$ for a theta cuspidal representation $\sigma$ of $\rU(V_0)$. 
We want to show  
\begin{equation*}\label{eq:conservation relation pro}
n_{\cV'_{V_l}}(\pi)+n_{\tcV'_{V_l}}(\pi)+c(\pi)  = 
2\dim V_l + \deltaan.
\end{equation*}
Let $(V'_0,\sigma')$ and $(\wtV'_0,\wtsigma')$ be  the first occurrences of $\sigma$ in the generalized Witt towers $\cV'_{V_0}$ and $\tcV'_{V_0}$ respectively.
Then both $\sigma'$ and $\wtsigma'$ are theta cuspidal.
Suppose $\pi$ corresponds to $\alpha\times \beta \in \sfW_l$ for some $\cc\in \bN$, $\alpha\in \Irr \sfS_\cc$ and  $\beta\in \Irr \sfS_{l-\cc}$ under the isomorphism $\cH_l\cong \sfH_{l,\half(\dim V'_0-\dim \wtV'_0),\nu=q}$. 
Then $n_{\cV'_{V_l}}(\pi)=\dim V'_{l-\rr(\beta)}$ by \Cref{cor:AMR}.
As for the generalized Witt tower $\tV'_V$, by \Cref{pull back iota} and \Cref{flip},  $\pi$ corresponds to $\beta \times \alpha \in \sfW_l$ under the isomorphism $\cH_l\cong \sfH_{l,\half(\dim \wtV'_0-\dim V'_0),\nu=q}$.
Similarly,  $ n_{\tcV'_{V_l}}(\pi)  =\dim \wtV'_{l-\rr(\alpha)} $ by \Cref{cor:AMR}.  
Since we have 
\[
\dim V'_0+\dim \widetilde{V}'_0= 2\dim V_0+\delta
\]
by the theta cuspidal case of \Cref{thm:conservation}, the proof is completed by the following lemma.
\end{proof}

\begin{lemma}\label{induction Wa}
We have 
$
\rr (\alpha)+\rr (\beta)=c(\pi).
$
\end{lemma}
\begin{proof}
Suppose that $\dd\in \bN$ and $\pi_0\in \Irr(\rU(V_{l-\dd}))$ such that 
\[
\Hom_{\rU(V_l)}\left(\pi, \Ind_{Q_\dd}^{\rU(V_l)}((\chi\circ \det)\boxtimes \pi_0\right)\neq 0,
\]
 then we have $\pi_0\in \cE(\rU(V_{l-\dd}),\sigma)$.
 Assume that $\pi_0$ corresponds to $\alpha_0\times \beta_0\in \Irr(\sfW_{l-\dd})$ under the isomorphism $\cH_{l-\dd}\cong \sfH_{l-\dd,\half(\dim V'-\dim \wtV'),\nu=q}$. 
 By  \Cref{lem:Hind} and Tits' deformation, we have
\[
\Hom_{\sfW_l}\left(\alpha\times \beta, \Ind_{\sfS_\dd\times \sfW_{l-\dd}}^{\sfW_l}\left(\bbone_{\dd}\boxtimes (\alpha_0\times \beta_0)\right)\right)\neq 0.
\]
By the Frobenius reciprocity and dimension counting, we have the following branching formula 
\[
\Ind_{\sfS_\dd}^{\sfW_\dd}(\bbone_{\dd})\cong \bigoplus_{a+b=\dd}\Ind_{\sfW_a\times \sfW_b}^{\sfW_{\dd}} (\bbone_{a}\boxtimes \varepsilon_{b})=\bigoplus_{a+b=\dd} \bbone_{a}\times \bbone_{b}.
\]
From this, we deduce that 
\[
\begin{split}
& \Ind_{\sfS_\dd\times \sfW_{l-\dd}}^{\sfW_l}\left(\bbone_{\dd}\boxtimes (\alpha_0\times \beta_0)\right)) \\
&\cong  \bigoplus_{a+b=\dd} \Ind_{\sfW_\dd\times \sfW_{l-\dd}}^{\sfW_l}\left( (\bbone_{a}\times \bbone_{b})\boxtimes (\alpha_0\times \beta_0)\right)) \\
& \cong  \bigoplus_{a+b=\dd} \sfX_{a}(\alpha_0) \times \sfX_{b}(\beta_0).
\end{split}
\]
The rest then follows. 
\end{proof}

\subsection{The comparison with Aubert-Michel-Rouquier \cite{AMR} and Pan \cite{PanUni}}\label{comparison normalization}

\Cref{cor:AMR}  covers the results of Aubert-Michel-Rouquier \cite[Theorem 3.10]{AMR} and 
Pan \cite[Corollary 3.5]{PanUni} for theta lifts 
between unipotent representations and quadratic-unipotent representations (called $\theta$-representations in \cite{LW}*{Theorem~3.3}). The formulation of \Cref{cor:AMR} is slightly different from loc. cit. We will explain these differences in Case $A$ and other cases are similar. 

Let us first recall the results in \cite{AM} on theta lifts of cuspidal unipotent representations in Case $A$. By Lusztig (see \cite{L}*{Section~8} and \cite{AM}*{Theorem~5.1}), if $G=\rU(V)$ is a unitary group, 
then $G$ has a cuspidal unipotent representation if and only if $\dim V=\frac{m(m+1)}{2}$ for some non-negative integers $m$.
In this case, there is a unique irreducible cuspidal unipotent representation and we denote it by $\sigma(m)$.
In \cite{AM}, it is proved that the first occurrence of $\sigma(m)$ in two different Witt towers are $\sigma(m-1)$ and $\sigma(m+1)$ respectively.
Combining these with \Cref{thm:norm}, we know that $|\mu(\sigma(m))|=m+\half$. This coincides with Lusztig's computation of $|\mu(\sigma(m))|$ in \cite[\S 4.6]{L}.

\trivial[h]{
\zjl{use one sentence}
\begin{thm}[\cite{AM}]\label{AMR1}
The\qcl{} first occurrence of $\sigma(m)$ appears on $\rU(V')$ and $\rU(\widetilde{V'})$ for $\dim V'=\frac{(m-1)m}{2}$ and $\dim \widetilde{V}'=\frac{(m+1)(m+2)}{2}$. We have
\[
\Theta_{V,V'}(\sigma(m))=\sigma(m-1)\quad \mbox{and}\quad \Theta_{V,\widetilde{V'}}(\sigma(m))=\sigma(m+1).
\]
\end{thm}
\begin{remark}
Combining \Cref{AMR1} with \Cref{thm:norm}, we know that $|\mu(\sigma(m))|=m+\half$. This coincides with Lustig's computation of $|\mu(\sigma(m))|$ in \cite[\S 4.6]{L}. 
\end{remark}
}

Let $\pi$ be a unipotent representation of a unitary group. By the classification of unipotent representations (\cite{L2}, also see \cite{AMR}*{Section~3}), there are some non-negative integers $l$ and $m$ such that $\pi\in \cE(\rU(V_l),\sigma(m))$ with $\dim V_0=\frac{m(m+1)}{2}$. Let $V'$ be an $\epsilon'$-Hermitian spaces of dimension $\frac{(m+1)(m+2)}{2}$ and set $(\sigma,\sigma'):=(\sigma(m),\sigma(m+1))$.  
Consider the theta lifts of $\pi$ to the generalized Witt tower $\cV'_{V}$ containing $V'$.
We know that $\Theta_{V_l,V'_{l'}}(\pi)$ lies in the Harish-Chandra series $\cE(\rU(V'_{l'}),\sigma(m+1))$. 
Appying \Cref{cor:AMR}, \Cref{pull back iota} and \Cref{flip}, we have the  following theorem due to Auber-Michel-Rouquier \cite[Theorem 3.10]{AMR} (cf. \cite[Corollary 3.5]{PanUni}).  
\begin{thm}\label{AMRu}
Adopt Lusztig's normalization (see \cite[\S 4.6]{L}): $\cH_l\cong \sfH_{l,m+\half,\nu=q}$ and $\cH_{l'}\cong \sfH'_{l,m+\tfrac{3}{2},\nu=q}$. 
Suppose $\pi\in \cE(\rU(V_l),\sigma(m))$
and $\pi$ corresponds to $\alpha\times\beta\in \Irr(\sfW_l)$. 
Then $\Theta_{V_l,V'_{l'}}(\pi)$ is a multiplicity free combination of  representations in $\cE(\rU(V'_{l'}),\sigma(m+1))$, and  corresponds to 
\[
\sum_{k=0} ^{\min\set{l,l'}}  \sfX_{l'-k} (\alpha) \times \sfX_{l-k}^* (\beta). 
\]
\end{thm}

\appendix
\section{Explicit computation  of $\cT_{\sftt_{l}}\cdot \indbasisgen$}\label{appendix A}
We prove the first equation in \Cref{lem:TlAk}. 
It suffices to compute the action of unnormalized Hecke operator $T_{\sftt_{l}}$ on $\indbasisgen$. By \Cref{actformula}, we have 
 \begin{equation}\label{firststep}
 \begin{split}
T_{\sftt_{l}}\cdot \indbasisgen =&
T_{\sftt_{l}}\cdot\left( I_{\cAkk}\otimes I_{0}\otimes \iota_0 \right)=\sum_{u\in U_{\sftt_{l}}} \omega_{V_l\otimes V_{l'}'}(u t_{l} )\left( I_{\cAkk}\otimes I_{0}\otimes \iota_0 \right).
 \end{split}
\end{equation}
By \Cref{Nkl} and \Cref{shape of U_t}, we have an exact sequence 
\[
\begin{tikzcd}[row sep=0em,ampersand replacement=\&]
1 \ar[r] \& \Herm(\langle v_{-l} \rangle , \langle v_l \rangle )  \ar[r, "{c\mapsto u(0,c)}"] \& U_{\sftt_{l}}  \ar[r,"{u(b,c)\mapsto b}"] \&\Hom(V_0,\langle v_l \rangle  )\ar[r] \& 1.
\end{tikzcd}
\]
For each $b\in  \Hom(V_0,\langle v_l \rangle )$, we choose an element $u(b,\lambda_b)\in U_{t_l}$. Then we may write 
\begin{equation}\label{shape of U_t 2}
U_{t_l}= \{u(b,\lambda_b)\com u(0,c)|b\in \Hom(V_0,\langle v_l \rangle ), c\in \Herm(\langle v_{-l} \rangle ,\langle v_l \rangle )\}.
\end{equation} 
 It follows by \Cref{firststep} and \Cref{shape of U_t 2} that 
\[
T_{\sftt_{l}}\cdot \indbasisgen 
=\hspace{-1em} \sum_{b\in \Hom(V_0,\langle v_l \rangle  )}\sum_{c\in  \Herm(\langle v_{-l} \rangle , \langle v_l \rangle )}\hspace{-2em} \omega_{V_l\otimes V_{l'}'}(u(b,\lambda_b)u(0,c)t_{l}) \left( I_{\cAkk}\otimes I_{0}\otimes \iota_0 \right). 
\]

 \trivial[h]{
Change to adapte the case when $p=2$. 
By \Cref{actformula} and \Cref{shape of U_t}, we have 
 \begin{equation*}\label{firststep}
 \begin{split}
T_{\sftt_{l}}\cdot \indbasisgen =&
T_{\sftt_{l}}\cdot\left( I_{\cAkk}\otimes I_{0}\otimes \iota_0 \right)=\sum_{u\in U_{\sftt_{l}}} \omega_{V_l\otimes V_{l'}'}(u t_{l} )\left( I_{\cAkk}\otimes I_{0}\otimes \iota_0 \right)\\
=&\sum_{b\in \Hom(V_0,\langle v_l \rangle  )}\sum_{c\in  \Herm(\langle v_{-l} \rangle , \langle v_l \rangle )} \omega_{V_l\otimes V_{l'}'}(u(b))\omega_{V_l\otimes V_{l'}'}(u(0,c)) \omega_{V_l\otimes V_{l'}'}(t_{l})\cdot \left( I_{\cAkk}\otimes I_{0}\otimes \iota_0 \right). 
 \end{split}
\end{equation*}
}
 \trivial[h]{
 We compute this in three steps:  first compute the action of $\omega_{V_l\otimes V_{l'}'}(t_{l})$ first and then compute the truncation \qcl{this word is not clear to reader} of $\Herm(\langle v_{-l} \rangle , \langle v_l \rangle )$ and $\Hom(V_0,\langle v_l \rangle)$. 
 }

\SSS
First, we compute $\omega_{V_l\otimes V_{l'}'}(t_{l})\, \left( I_{\cAkk}\otimes I_{0}\otimes \iota_0 \right)$.
Since $\omega_{V_l\otimes V_{l'}'}(t_{l})$  acts by the partial Fourier transform on $\Hom (\langle v_{l}\rangle, V_{l'}')$ as in  \Cref{Partial Fourier}, to simplify the computation, we  write $\cA_k$ as the difference of two ``hyperplanes'' as follows.  
Via the isomorphism \Cref{isomorphism V_l+}, define
\[
\Ak' :=\sum_{1\leq i\leq k-1}v_{-(l-k+i)}\otimes v_{i}' \in \Hom(V_l^+, V'_{l'}).
\]
  (Compare with $\Akk$ defined in \eqref{Akkdef}). 
  Let $\cAk':= P'_{l'}\com   \Ak'\com   B(V_l^+)$. 
Define 
\[
\begin{split}
\barcAkk & := \cAk'  + \langle v_{-l}\rangle \otimes V_{k}'^+,\quad
\quad\quad  \pcAkk   :=  \cAk'  + \langle v_{-l}\rangle \otimes V_{k-1}'^+,\\
\barcAkk^\perp & :=  \cAk'  + \langle v_{-l}\rangle \otimes (V_{k}'^+)^\perp,\quad
\pcAkk^\perp   := \cAk' +\langle  v_{-l}\rangle \otimes (V_{k-1}'^+)^\perp.
\end{split}                     
\]
Then $\cAkk = \barcAkk - \pcAkk$ and $I_{\cAkk} = I_{\barcAkk}-I_{\pcAkk}$.  Applying \Cref{Partial Fourier}, we have 
\begin{equation*}\label{action of tl}
\begin{split}
&\omega_{V_l\otimes V_{l'}'}(t_{l})\, \big( I_{\cAkk}\otimes I_{0}\otimes \iota_0 \big)\\
& = \gamma_{V_0'} q^{\dim V_{k}'^+-\half \dim  V_{l'}'} \big( I_{\barcAkk^\perp}\otimes I_{0}\otimes \iota_0 \big) -\gamma_{V_0'} q^{\dim V_{k-1}'^+-\half \dim  V_{l'}'} \big( I_{\pcAkk^\perp}\otimes I_{0}\otimes \iota_0 \big) \\
& = \gamma_{V_0'} q^{k-l'-\half\dim V_0'}\left( \big(I_{\barcAkk^\perp}\otimes I_{0}\otimes \iota_0 \big)- q^{-1} \big(I_{\pcAkk^\perp}\otimes I_{0}\otimes \iota_0 \big) \right).
\end{split}
\end{equation*}

\SSS We now consider the action of $u(0,c)$.
Let 
\begin{equation}\label{CN'}
\cN':=\set{v'\in V'_{l'}| \langle v', v'\rangle_{V'_{l'}}=0}.
\end{equation}
It follows from \Cref{mix model equation} that 
\begin{equation}\label{action of c 1}
\begin{split}
& \sum_{c\in  \Herm(\langle v_{-l} \rangle , \langle v_l \rangle )}u(0,c)\left(I_{\barcAkk^\perp}\otimes I_{0}\otimes \iota_0 \right)\\
& = q^{\half\deltaan}\left(I_{\cAk'+\langle v_{-l} \rangle \otimes \left((V_{k}'^+)^\perp\cap \cN' \right)}\otimes I_{0}\otimes \iota_0 \right),
\end{split}
\end{equation}
and 
\begin{equation}\label{action of c 2}
\begin{split}
& \sum_{c\in  \Herm(\langle v_{-l} \rangle , \langle v_l \rangle )}u(0,c)\left(I_{\pcAkk^\perp}\otimes I_{0}\otimes \iota_0 \right) \\
& =  q^{\half\deltaan} \left(I_{\cAk'+\langle v_{-l} \rangle \otimes \left((V_{k-1}'^+)^\perp\cap \cN' \right)}\otimes I_{0}\otimes \iota_0 \right).
\end{split}
\end{equation}

\trivial[h]{
It follows from \Cref{action of tl}, \Cref{action of c 1} and \Cref{action of c 2} that the support of 
\[
\sum_{c\in  \Herm(\langle v_{-l} \rangle , \langle v_l \rangle )} \omega_{V_l\otimes V_{l'}'}(u(0,c)) \omega_{V_l\otimes V_{l'}'}(t_{l})\cdot \left( I_{\cAkk}\otimes I_{0}\otimes \iota_0 \right). 
\]
is contained in $\cA_{k,0}+\langle v_{-l} \rangle \otimes \left((V_{k-1}'^+)^\perp\cap \cN' \right)$.
}

\SSS
Note that the action of   $u(b,\lambda_b)$ on a general $f\in \left(\Weil{V_l',V_l}\right)_{\sigmaP\times \sigmaPp}$
will not enlarge the support of $f$ as a function on $\Hom(V^+_l,V'_l)$ by \Cref{mix model equation}.  
In view of \Cref{support on a single orbit}, we now classify the $B_{V^+_l}\times P'_{l'}$-orbits in $\cAk'+\langle v_{-l} \rangle \otimes \left((V_{k-1}'^+)^\perp\cap \cN' \right)$. 
Recall that we have a partial flag  
\[
0=V_0'^+\subseteq V_1'^+\subseteq \cdots \subseteq V_{l'}'^+ \subseteq (V_{l'}'^+)^\perp \subseteq\cdots\subseteq (V_{1}'^+)^\perp \subseteq (V_{0}'^+)^\perp= V_{l'}'.
\]
For each $i=1,\cdots, l'$, let 
\begin{equation}\label{Nk and N-k}
\cN'_i:= (V_{i}'^+\backslash V_{i-1}'^+)\subset \cN' 
\quad \mbox{and} \quad 
\cN'_{-i}:=\left((V_{i-1}'^+)^\perp\backslash (V_{i}'^+)^\perp\right)\cap \cN'.
\end{equation}
For each $i=k,\cdots, l'$, define 
\[
\cA_{k,i}\coloneqq  \cAk'+\langle v_{-l} \rangle \otimes \cN'_i \quad \mbox{and}\quad
\cA_{k,-i}\coloneqq  \cAk'+\langle v_{-l} \rangle \otimes \cN'_{-i} .
\]
Note that $\cA_{k,k}=\cA_{k}$ by \eqref{eq:Akk}. We also define 
\[\cN'_{\emptyset} := \left( (V_{l'}'^+)^\perp  \backslash V_{l'}'^+\right )\cap \cN',\quad
\cA_{k,\emptyset}\coloneqq   \cAk'+\langle v_{-l} \rangle \otimes \cN'_{0} .
\]

{\trivial[h]{\begin{lemma}\label{geometry of orbits 1}
\qcl{state where each Aki is. E.g. Note that $\cA_{k,0}\in \cZ_{k-1}$. }

\begin{enumerate}
    \item For each $i=k,\cdots, l'$, $\cA_{k,i}$ and $\cA_{k,-i}$ are single $B_{V^+_l}\times P'_{l'}$-orbits.
    \item \qcl{Don't start with a symbol.} The set $\cA_{k,V_0'}$ is non-empty if and only if $\dim V'_0\neq 0$, in which case $\cA_{k,V_0'}$ is a single $B_{V^+_l}\times P'_{l'}$-orbit.
    \item We have 
\[
\cA_{k,0}+\langle v_{-l} \rangle \otimes \left((V_{k-1}'^+)^\perp\cap \cN' \right)=\left( \coprod_{i=k}^{l'} \cA_{k,i}\right)\coprod 
\left(\coprod_{i=k}^{l'} \cA_{k,-i} \right) \coprod \cA_{k,0}\coprod \cA_{k,V_0'} 
\]
and 
\[
\cA_{k,0}+\langle v_{-l} \rangle \otimes \left((V_{k}'^+)^\perp\cap \cN' \right)=\left( \coprod_{i=k}^{l'} \cA_{k,i}\right)\coprod 
\left(\coprod_{i=k+1}^{l'} \cA_{k,-i} \right) \coprod \cA_{k,0}\coprod \cA_{k,V_0'}. 
\]
\end{enumerate}
\end{lemma}
We give a prove here: 
First we give another description of $\cB_k$: Let $\cN'_{-k}$ be the orbit of $P'$ on $v_{-k}'$. Then we have 
\[
\cB_k= \cA_{k,0}+\langle v_{-l} \rangle \otimes \cN'_{-k}. 
\]
The proof is based on \Cref{B_k as an orbit} and some linear algebra (omit here). Then we need to prove 
\[
\cN'_{-k}= \left( (V_{k-1}'^+)^\perp\backslash (V_{k}'^+)^\perp\right) \cap \cN'
\]
We give a proof by parameterize the $P'$-orbits of $\cN'$. It follows from the Witt theorem that $\cN'\backslash \{0\}$ is a single $\rU(V'_{l'})$-orbit. The stabilizer of $v_1\in (\cN'\backslash \{0\})$ is the parabolic subgroup $Q_1'$. Therefore, we have a bijection between the sets of $P'$-orbits of $\cN'\backslash \{0\}$ and the double coset $P'\backslash \rU(V_{l'}')/ Q_1'$, which is also a double coset $W_{L'}\backslash W_{\rU(V_{l'}')} /W_{M_1}$. 

\begin{lemma}\label{P' orbits of nipotent cone}
The double coset $W_{L'}\backslash W /W_{M_1}$ has cardinality $2l'$ if $\dim V'=0$ and $2l'+1$ if $\dim V'\geq 1$. 
\end{lemma}
\begin{proof}
We give a proof when $\rU(V_{l'}')$ is of type $C$, the other case is similar (or maybe the same). In this case, assume that $\dim V'=2r$. We have $W_{\rU(V_{l'}')}\cong W_{l'+r}, W_{M_1}\cong W_{l'+r-1}$ and $W_L\cong W_{r}$. The embedding of $W_{l'+r-1}$ and $ W_{r}$ to $W_{l'+r}$ is by blabla. The rest follows by direct computation. 
\end{proof}
On the other hand, we know that $P'$ stabilizes the flags 
\begin{equation}\label{flag}
0\subseteq V_1'^{+}\subseteq \cdots V_{l'}'^{+}\subseteq (V_{l'}'^+)^\perp\subseteq \cdots (V_{1}'^+)^\perp \subseteq V_{l'}
\end{equation}
where $(V_{l'}'^+)^\perp=V_{l'}'^+\oplus V$. Therefore, $(V_{l'}'^+)^\perp=V_{l'}'^+$ if $\dim V'=0$ and $=V_{l'}'^+\nsubseteq (V_{l'}'^+)^\perp$ if $\dim V'\geq 1$. The difference of any two successive subspace intersects with $\cN'$ is stable under the action of $P'$. Combining \cref{P' orbits of nipotent cone}, we deduce that any two successive subspace intersects with $\cN'$ is a single orbit.
}
}
The following lemma is routine to check. 
\begin{lemma}\label{geometry of orbits 1}
\begin{enumerate}[label=(\arabic*),wide=0em]
    \item We have 
 $\cAk'= \cO^{k-1}_{\sfss_{l-k+1,l},1, 1}.$ For $i=k,\cdots, l'$,  
 \[
    \cA_{k,i}=\cO^k_{1,\sfss'_{k,i},1}\quad  \text{and} \quad \cA_{k,-i}=\cO^k_{1,\sftt'_i\sfss'_{k,i},1}.
 \]
    \item   The set $\cA_{k,\emptyset}\neq \emptyset$   if and only if $\dim V'_0\neq 0$, in which case $\cA_{k,\emptyset}$ is a single $B_{V^+_l}\times P'_{l'}$-orbit. Moreover $\cA_{k,\emptyset}$ is not  of the form $\cO^i_{\sfdd_1,\sfdd_2,\sfxx}$ for any $0\leq i\leq \min\set{l,l'}$ and $(\sfdd_1,\sfdd_2,\sfxx)\in \sfD_i\times \sfD'_i\times \sfS_i$. 
    \item We have the following decompositions:
\begin{equation}\label{orbit decoms 1}
\cAk'+\langle v_{-l} \rangle \otimes \left((V_{k-1}'^+)^\perp\cap \cN' \right)=\left( \bigsqcup_{i=k}^{l'} \cA_{k,i}\right)\sqcup
\left(\bigsqcup_{i=k}^{l'} \cA_{k,-i} \right) \sqcup \cAk'\sqcup \cA_{k,\emptyset}; 
\end{equation}
\begin{equation}\label{orbit decoms 2}
\cAk'+\langle v_{-l} \rangle \otimes \left((V_{k}'^+)^\perp\cap \cN' \right)=\left( \bigsqcup_{i=k}^{l'} \cA_{k,i}\right)\sqcup
\left(\bigsqcup_{i=k+1}^{l'} \cA_{k,-i} \right) \sqcup \cAk'\sqcup\cA_{k,\emptyset}.
\end{equation}
\end{enumerate}
\end{lemma}
\trivial[h]{
\begin{itemize}
    \item Define $\sfss'_{k,i}= \sfss'_{i-1}\cdots s'_k$ for $i\geq k+1$. Then it is easy to check that $\sfss'_{k,i}\in \sfD'_k$ and $\cA_{k,i}=\cO^k_{1,\sfss'_{k,i},1}$ for $i=k+1,\cdots,l'$. In another words
    \[
    \cA_{k,i}=\cO^k_{\sfdd_1,\sfdd_2,\sfxx}\quad \mbox{for} \quad d_1=1,d_2=(1,\cdots,k-1,i,k,\cdots i-1,i+1,\cdots, l'), x=1.
    \]
    \item Define $\sftt'_{k,i}= \sftt'_{i}\sfss'_{i-1}\cdots s'_k$ for $i\geq k+1$. Then it is easy to check that $\sftt'_{k,i}\in \sfD'_k$ and $\cA_{k,i}=\cO^k_{1,\sftt'_{k,i},1}$ for $i=k+1,\cdots,l'$ and $\cA_{k,-k}=\cO^k_{1,\sftt'_{k},1}$. In another words
     \[
    \cA_{k,i}=\cO^k_{\sfdd_1,\sfdd_2,\sfxx}\quad \mbox{for} \quad d_1=1,d_2=(1,\cdots,k-1,-i,k,\cdots i-1,i+1,\cdots, l'), x=1.
    \]
    \item $\cA_{k,V_0'}$ is not of the form $\cO^k_{\sfdd_1,\sfdd_2,\sfxx}$ for some $(\sfdd_1,\sfdd_2,\sfxx)\in \sfD_k\times \sfD'_k\times \sfS_k$. 
\end{itemize}



Now it follows from \Cref{geometry of orbits 1} that 
\begin{equation}\label{break up 1}
I_{\cA_{k,0}+\langle v_{-l} \rangle \otimes \left((V_{k-1}'^+)^\perp\cap \cN'\right)}\otimes I_0\otimes \iota_0=  \sum_{i=k}^{l'} (I_{\cA_{k,i}}\otimes I_0\otimes \iota_0)+ \sum_{i=k}^{l'} (I_{\cA_{k,-i}}\otimes I_0\otimes \iota_0)+I_{ \cA_{k,0}}\otimes I_0\otimes \iota_0 + I_{\cA_{k,V_0'}}\otimes I_0\otimes \iota_0 
\end{equation}
and 
\begin{equation}\label{break up 2}
I_{\cA_{k,0}+\langle v_{-l} \rangle \otimes \left((V_{k}'^+)^\perp\cap \cN'\right)}\otimes I_0\otimes \iota_0=  \sum_{i=k}^{l'} (I_{\cA_{k,i}}\otimes I_0\otimes \iota_0)+ \sum_{i=k+1}^{l'} (I_{\cA_{k,-i}}\otimes I_0\otimes \iota_0)+ I_{\cA_{k,0}}\otimes I_0\otimes \iota_0 + I_{\cA_{k,V_0'}}\otimes I_0\otimes \iota_0 
\end{equation}
}

\SSS
The right-hand sides of 
\Cref{action of c 1} and \Cref{action of c 2} decomposes as constant functions supported on     $B_{V_l^+}\times P'_{l'}$-orbits appearing in \Cref{orbit decoms 1} and \Cref{orbit decoms 2}.
The following lemma computes the average over $\sum u(b,\lambda_b)$  on each term in the decomposition.
\begin{lemma}\label{action of u(b)}
We have 
\begin{equation}\label{eq:ubaction}
\begin{split}    
&  \sum_{b\in \Hom(V_0, \langle v_l \rangle )} \Weil{V_l,V'_{l'}}(u(b,\lambda_b))  \left(I_{\star}\otimes  I_{0}\otimes \iota_0 \right)\\
& = 
\begin{cases}
q^{\dim V_0} \II^k_{1,\sfss'_{k,i},1}\quad & \mbox{if $\star=\cA_{k,i}$ and $i\in \set{k,\cdots l'}$},\\
 -q^{\dim V_0'+1-\half \deltaan} \II^k_{1,\sftt'_i\sfss'_{k,i},1}\quad & \mbox{if $\star=\cA_{k,-i}$ and $i\in \set{k,\cdots, l'}$},\\
 q^{\dim V_0} \II^{k-1}_{\sfss_{l-k+1,l},1, 1} \quad & \mbox{if $\star=\cAk'$},\\
0\quad & \mbox{if $\star=\cA_{k,\emptyset}$ }.
\end{cases}
\end{split}
\end{equation}
\trivial[h]{
We have 
 \begin{enumerate}
    \item 
   For $i=k,\cdots, l'$    
    \[
    \sum_{b\in \Hom(V_0, \langle v_l \rangle )} \Weil{V_l,V'_{l'}}(u(b))\cdot  (I_{\cA_{k,i}}\otimes  I_{0}\otimes \iota_0 )=  q^{\dim V_0} \II^k_{1,\sfss'_{k,i},1}.
    \]
    \item
   For $i=k,\cdots, l'$,   
    \[
\sum_{b\in \Hom(V_0, \langle v_l \rangle )} \Weil{V_l,V'_{l'}}(u(b))\cdot ( I_{\cA_{k,-i}}\otimes I_{0}\otimes \iota_0)=-q^{\dim V_0'+1-\half \deltaan} \II^k_{1,\sftt'_i\sfss'_{k,i},1}.
    \]
    \item 
    \[
    \sum_{b\in \Hom(V_0, \langle v_l \rangle )} \Weil{V_l,V'_{l'}}(u(b))\cdot  (I_{\cAk'}\otimes  I_{0}\otimes \iota_0) =  q^{\dim V_0} \II^{k-1}_{\sfss_{l-k+1,l},1, 1}.
    \]
    \item 
       \[
\sum_{b\in \Hom(V_0, \langle v_l \rangle )} \Weil{V_l,V'_{l'}}(u(b))\cdot  (I_{\cA_{k,0}}\otimes I_{0}\otimes \iota_0)= 0.
    \]
 \end{enumerate}
 }
\end{lemma}
\begin{proof}
We only prove the first case and the proofs for the rest are similar. 
 Since the $u(b,\lambda_b)$ action does not enlarge the support,  
 it follows from \Cref{support on a single orbit} and \Cref{indbasisrealzation } that the left-hand side of \eqref{eq:ubaction} must be a scalar multiple of $\II^k_{1,\sfss'_{k,i},1}$.
 To find the scalar, we evaluate both sides of the equation at the representative element $A_{k,i}\coloneqq\sfss'_{k,i}\com \Akk$.
 By \Cref{mix model equation}, 
we have 
    \begin{equation*} 
    \begin{split}
     & \Weil{V_l,V'_{l'}}(u(b,\lambda_b)) \left(I_{\cA_{k,i}}\otimes I_{0}\otimes \iota_0 \right)(A_{k,i})\\
     =& \rho_{V_0,V'_{l'}}\left(\left(A_{k,i}\com b,-\iinn{A_{k,i}\com \lambda_b}{A_{k,i}}_{\Hom(V_l,V'_{l'})} \right)\right) \left(I_{0}\otimes \iota_0 \right)\\
     =&\rho_{V_0,V'_{l'}}((A_{k,i}\com b,0)) \left(I_{0}\otimes \iota_0 \right).
     \end{split}
    \end{equation*}
  Note that $A_{k,i}\com b \in \Hom(V_0, \langle v'_{i} \rangle ) \cong V_0\otimes\langle v'_{i} \rangle $. Then by \Cref{Heisenberg representation}, we have 
    \[
    \rho_{V_0,V'_{l'}}((A_{k,i}\com b,0))\, \left(I_{0}\otimes \iota_0 \right)= \psi\left( \tr_{F/F_0}\left(
    \inn{ A_{k,i}\com b}{0}_{V_0\otimes V'_{l'}}\right) \right)\left(I_{0}\otimes \iota_0 \right)=I_{0}\otimes \iota_0.
    \]
    
So the value of the left-hand side of \eqref{eq:ubaction} at $A_{k,i}$ is $q^{\dim V_0} \left(I_{0}\otimes \iota_0\right)$.
For the right-hand side of \eqref{eq:ubaction},  it follows from \Cref{indbasisrealzation } that  
\[ 
\II_{1,\sfss'_{k,i},1}(A_{k,i})=\II_0\otimes \iota_0 
\] 
since $\iota(\sfss'_{k,i})=0$ and $\bY_{\sfss'_{k,i}}^-=0$. 
This finished the proof of the first case. 
\end{proof}

\SSS
Finally, the previous discussion and our normalization of $\cT_{\sftt_l}$ in \Cref{Normalized t_k2} give the formula of   $
\cT_{\sftt_{l}}\cdot \indbasisgen$ in \Cref{lem:TlAk}.



\section{The explicit computation of $\cT_{\sftt_{l-k}}\cdot \indbasisgen$}\label{appendix B}
We prove the second equation in Proposition \ref{lem:TlAk}.     
It suffices to compute the action of the unnormalized Hecke operator $T_{\sftt_{l-k}}$.   
\trivial[h]{
By \Cref{actformula},  we have 
 \begin{equation*}\label{firststep 22}
 \begin{split}
T_{\sftt_{l-k}}\cdot \indbasisgen=&
T_{\sftt_{l}}\cdot\left( I_{\cAkk}\otimes I_{0}\otimes \iota_0 \right)
T_{\sftt_{l-k}}\cdot\left( I_{\cAkk}\otimes I_{0}\otimes \iota_0 \right)=\sum_{u\in U_{\sftt_{l-k}}} \omega_{V_l\otimes V_{l'}'}(u t_{l-k} )\left( I_{\cAkk}\otimes I_{0}\otimes \iota_0 \right).
 \end{split}
\end{equation*}
}
By \eqref{actformula}, we have 
 \begin{equation*}
 \begin{split}
T_{\sftt_{l-k}}\cdot \indbasisgen 
& = T_{\sftt_{l-k}}\cdot\left( I_{\cAkk}\otimes I_{0}\otimes \iota_0 \right) \\
&=\sum_{u\in U_{\sftt_{l-k}}} \omega_{V_l\otimes V_{l'}'}(u \, t_{l-k}) \left( I_{\cAkk}\otimes I_{0}\otimes \iota_0 \right). 
 \end{split}
\end{equation*}
Note that $U_{\sftt_{l-k}}\subseteq N(V_{l-k}^+,V_l)$ and we can write 
\[
U_{\sftt_{l-k}}=\Set{u(b,c)|\begin{array}{l} b\in \Hom(\widehat{V}_k,\langle v_{l-k}\rangle),\\
c\in \Hom(\langle v_{-(l-k)}\rangle, \langle v_{l-k}\rangle),\\ c+c^*+bb^*=0 
\end{array}}
\]
according to  \Cref{Nkl}. Using the decomposition 
\[
\widehat{V}_k=\widehat{V}^+_k\oplus V_0\oplus \widehat{V}^-_k
\]
as in \Cref{Vl-k}, we may also write 
\[
\begin{split}
U_{\sftt_{l-k}}=\Set{u(b_2,\lambda_{b_2})u(b_1,0)u(b_3,0)u(0,c)|
\begin{array}{l}
b_1\in \Hom(\wtV_{k}^+,\langle v_{l-k}\rangle ),\\
b_2\in \Hom(V_0,\langle v_{l-k}\rangle  ),\\
b_3\in \Hom(\wtV_{k}^-,\langle v_{l-k} \rangle  ),\\
c\in \Herm(\langle v_{l-k}\rangle , \langle v_{l-k}\rangle )
\end{array}
}, 
\end{split}
\]
where $u(b_2,\lambda_{b_2})$ is a fixed element in $U_{\sftt_{l-k}}$ for each $b_2\in \Hom(V_0,\langle v_{l-k}\rangle  )$.
Then we have 
 \begin{equation*}
 \begin{split}
& T_{\sftt_{l-k}}\cdot \indbasisgen \\
&=  \sum_{b_2\in \Hom(V_0,\langle v_{l-k}\rangle  )}\sum_{b_1\in \Hom(\widehat{V}_{k}^+,\langle v_{l-k}\rangle ) } \sum_{b_3\in \Hom(\widehat{V}_{k}^-,\langle v_{l-k}\rangle )} \\ 
& \ \  \sum_{c\in  \Herm(\langle v_{-(l-k)} \rangle , \langle v_{l-k} \rangle )}   
\hspace{-3.5em} \omega_{V_l\otimes V_{l'}'}(u(b_2,\lambda_{b_2})\, u(b_1,0)\, u(b_3,0)\, u(0,c)\, t_{l-k})
\left( I_{\cAkk}\otimes I_{0}\otimes \iota_0 \right). 
 \end{split}
\end{equation*}
\SSS
First, we compute $\omega_{V_l\otimes V_{l'}'}(t_{l-k}) \left( I_{\cAkk}\otimes I_{0}\otimes \iota_0 \right)$. Since $I_{\cAkk}$ is supported on $0$ when restricted to $v_{-(l-k)}\otimes V_{l'}'$, it follows from \Cref{Partial Fourier} that 
\begin{equation*}\label{action of tl-k}
 \omega_{V_l\otimes V_{l'}'}(t_{l-k}) \left( I_{\cAkk}\otimes I_{0}\otimes \iota_0 \right) = \gamma_{V_{0}'} q^{-\half \dim  V_{l'}'}  \left( I_{\cAkk+ \langle v_{-(l-k)}\rangle \otimes V_{l'}'}\otimes I_{0}\otimes \iota_0 \right).   
\end{equation*}
\SSS
It follows from \Cref{mix model equation} that 
\begin{equation*}\label{truncation of c}
\begin{split}
& \sum_{c\in \Herm(\langle v_{-(l-k)}\rangle , \langle v_{(l-k)}\rangle )}\omega_{V_l\otimes V_{l'}'}(u(0,c)) \left( I_{\cAkk+ \langle v_{-(l-k)}
\rangle \otimes V_{l'}'}\otimes I_{0}\otimes \iota_0 \right) \\
& =q^{\half \deltaan} \left( I_{ \cAkk + \langle v_{-(l-k)}\rangle \otimes \cN'} \otimes I_{0}\otimes \iota_0 \right).
\end{split}
\end{equation*}
Here $\cN'$ is defined in \eqref{CN'}.

\SSS
 By a routine computation, 
we deduce from the fourth equation of \Cref{mix model equation} that 
\begin{equation*}\label{truncation of b3}
\begin{split}
& \sum_{b_3\in \Hom(\widehat{V}_{k}^-,\langle v_{l-k}\rangle )}\omega_{V_l\otimes V_{l'}'}(u(b_3,0)) \left( I_{ \cAkk + \langle v_{-(l-k)}\rangle \otimes \cN'} \otimes I_{0}\otimes \iota_0 \right) \\
& =q^k \left( I_{ \cAkk + \langle v_{-(l-k)}\rangle \otimes \left( (V_k'^{+})^\perp\cap \cN'\right)} \otimes I_{0}\otimes \iota_0 \right).
\end{split}
\end{equation*}
\trivial[h]{$\bullet$ \underline{The truncation of $\Hom(\widehat{V}_{k}^-,\langle v_{l-k}\rangle )$}: It follows from \Cref{mix model equation} that  
\[
\begin{split}
&\sum_{b_3\in \Hom(\widehat{V}_{k}^-,\langle v_{l-k}\rangle )}\omega_{V_l\otimes V_{l'}'}(u(b_3))\cdot \left( I_{ \cAkk + \langle v_{-(l-k)}\rangle \otimes \cN'} \otimes I_{0}\otimes \iota_0 \right)\\
=&\sum_{B\in \cAkk + \langle v_{-(l-k)}\rangle \otimes \cN'}\sum_{b_3\in \Hom(\widehat{V}_{k}^-,\langle v_{l-k}\rangle )}\psi(\langle u(b_3)\cdot B,B \rangle)\cdot \left( I_{B} \otimes I_{0}\otimes \iota_0 \right).
\end{split}
\]
Write $B=A+v_{-(l-k)}\otimes v'$ for $A\in \cAkk$ and $v'\in \cN'$. Then it is easy to check that the character \[
b_3 \mapsto \psi(\langle u(b_3)\cdot B,B \rangle)
\]
of $\Hom(\widehat{V}_{k}^-,\langle v_{l-k}\rangle ) )$ is trivial if and only if $v'\in \cN'\cap (V_k'^{+})^\perp$.
}
 

\SSS
The action of $u(b_1,0)$ is more complicated. It follows from the second equation of  
\Cref{mix model equation} that  
\begin{equation}\label{trucation of Vk+}
\begin{split}
&\sum_{b_1\in \Hom(\widehat{V}_{k}^+,\langle v_{l-k}\rangle ) )}\omega_{V_l\otimes V_{l'}'}(u(b_1,0)) \left( I_{ \cAkk + \langle v_{-(l-k)}\rangle \otimes \left( (V_k'^{+})^\perp\cap \cN'\right)} \otimes I_{0}\otimes \iota_0 \right)\\
=&\sum_{b_1\in \Hom(\widehat{V}_{k}^+,\langle v_{l-k}\rangle )}\sum_{B\in \cAkk + \langle{v_{-(l-k)}}\rangle \otimes \left( (V_k'^{+})^\perp\cap \cN'\right)}
\left( I_{B \cdot u(b_1,0)^{-1}} \otimes I_{0}\otimes \iota_0 \right).
\end{split}
\end{equation}
To compute \eqref{trucation of Vk+},
we consider $B(V_l^+)\times P'_{l'}$-orbits in 
  the image of 
\begin{equation}\label{geometric action of b_1}
\begin{split}
\Hom(\whV_{k}^+,\langle v_{l-k}\rangle ) \times \Big(\cAkk + \langle v_{-(l-k)}\rangle \otimes \big( (V_k'^{+})^\perp\cap \cN'\big)\Big)&\longrightarrow \cZ\\
(b_1, B) &\mapsto B\cdot u(b_1)^{-1}.
\end{split}
\end{equation}

For each $i=1,\cdots, l'$, define 
\[
\cB^\circ_{k,i} :=   \cAkk+\langle v_{-(l-k)} \rangle \otimes  \cN'_i \quad \mbox{and} \quad  \cB^\circ_{k,-i} :=   \cAkk+\langle v_{-(l-k)} \rangle \otimes  \cN'_{-i}, 
\]
where $\cN'_i$ and $\cN'_{-i}$ are defined in \Cref{Nk and N-k}. We also define $\cB^\circ_{k,0}:=\cAkk$ and 
\[
\cB^\circ_{k,\emptyset}:=\cAkk+\langle v_{-(l-k)} \rangle \otimes  \cN'_{\emptyset}. 
\]
Note that $\cB^\circ_{k,\emptyset}$ is non-empty if and only if $\dim V'_0\neq 0$. Then we have 
\[
\cAkk+\langle v_{-(l-k)} \rangle \otimes \left((V_{k}'^+)^\perp\cap \cN' \right)=\left( \bigsqcup_{i=0}^{l'} \cB^\circ_{k,i}\right)\sqcup
\left(\bigsqcup_{i=k+1}^{l'} \cB^\circ_{k,-i} \right) \sqcup \cB^\circ_{k,\emptyset}. 
\]
The following proposition is easy to check.
\begin{prop}\label{orbits Bki}
The restriction of \Cref{geometric action of b_1} to $\Hom(\whV_{k}^+,\langle v_{l-k}\rangle)  \times \cB^\circ_{k,i}$ is 
\[
\begin{cases}
\mbox{a $q^{k-i}(q-1)$ to $1$ map}\quad & \mbox{if $i=1,\cdots,k$,}\\
\mbox{a $q^{k}$ to $1$ map}\quad & \mbox{if $i=0$,}\\
\mbox{an injection}\quad  & \mbox{if $i=\pm(k+1),\cdots,\pm l'$ or $\emptyset$.}\\
\end{cases}
\]
\end{prop}
For each $i$ above, we define $\cB_{k,i}$ to be the image of $\Hom(\whV_{k}^+,\langle v_{l-k}\rangle)\times \cB^\circ_{k,i}$. Then a standard calculation shows that 
\[
\cB_{k,i}= \begin{cases}
\cO^k_{\sfss_{l-k,l-k+i},1,\sfss'_{1,i}}\quad &\mbox{if $i=1,\cdots,k$},\\
\cO^{k+1}_{1,\sfss'_{k,i},\sfss'_{1,k}}\quad &\mbox{if $i=k+1,\cdots, l'$,}\\
\cO^{k+1}_{1,\sftt'_{i}\sfss'_{k,i},\sfss'_{1,k}}\quad &\mbox{if $i=-(k+1),\cdots, -l'$}\\
\cO^k_{1,1,1}\quad &\mbox{if $i=0$}.
\end{cases}
\]
Note that $\cB_{k,\emptyset}$ is not of the form $\cO^i_{\sfdd_1,\sfdd_2,\sfxx}$ for any $1\leq i\leq \min\set{l,l'}$ and $(\sfdd_1,\sfdd_2,\sfxx)\in \sfD_i\times \sfD'_i\times \sfS_i$.

Now  \Cref{trucation of Vk+}  equals
\begin{equation}\label{truncation of ub1 equation}
\begin{split}
&q^k \left( I_{\cAkk} \otimes I_{0}\otimes \iota_0\right)+ \sum_{i=1}^k q^{k-i}(q-1) \left( I_{\cB_{k,i}} \otimes I_{0}\otimes \iota_0\right) \\
& + \sum_{i=k+1}^{l'}\left( I_{\cB_{k,i}} \otimes I_{0}\otimes \iota_0\right)
+ \sum_{i=k+1}^{l'}\left( I_{\cB_{k,-i}} \otimes I_{0}\otimes \iota_0\right)+
 \left( I_{\cB_{k,\emptyset}} \otimes I_{0}\otimes \iota_0\right) 
\end{split}  
\end{equation}
 
\SSS
The final computation is to average the action of $u(b_2,\lambda_{b_2})$ to each term in \Cref{truncation of ub1 equation}. 
Similar to \Cref{action of u(b)}, we have  
\[
\begin{split}
 &  \sum_{b_2\in \Hom(V_{0},\langle v_{l-k}\rangle )}\omega_{V_l\otimes V_{l'}'}(u(b_2,\lambda_{b_2})) \left(I_{\star}\otimes I_{0}\otimes \iota_0 \right)\\
 & =  
\begin{cases}
q^{\dim V_0}\II^k_{\sfss_{l-k,l-k+i},1,\sfss'_{1,i}}\quad & \mbox{if $\star=\cB_{k,i}$ and $i\in \set{1,\cdots, k}$},\\
q^{\dim V_0}\II^{k+1}_{1,\sfss'_{k,i},\sfss'_{1,k}}\quad & \mbox{if $\star=\cB_{k,i}$ and $i\in \set{k+1,\cdots, l'}$},\\
-q^{\dim V_0'+1-\half \deltaan}\II^{k+1}_{1,\sftt'_{i}\sfss'_{k,i},\sfss'_{1,k}} \quad & \mbox{if $\star=\cB_{k,-i}$ and $i\in \set{k+1,\cdots, l'}$},\\
q^{\dim V_0}\II^k_{1,1,1}\quad & \mbox{if $\star=\cB_{k,0}=\cAkk$ },\\
0\quad & \mbox{if $\star=\cB_{k,\emptyset}$ }.
\end{cases}
\end{split}
\]
\trivial[h]{
\begin{enumerate}
    \item For $i=1,\cdots, l'$, 
    \[
     \begin{split}
    \sum_{b_2\in \Hom(V_{0},\langle v_{l-k}\rangle )}\omega_{V_l\otimes V_{l'}'}(u(b_2,\lambda_{b_2}))\cdot \left(I_{\cB_{k,i}}\otimes I_{0}\otimes \iota_0 \right)
    = q^{\dim V_0}\II_{\cB_{k,i}}.
     \end{split}
    \]
    \item For $i=k+1,\cdots, l'$, 
      \[
    \sum_{b_2\in \Hom(V_{0},\langle v_{l-k}\rangle )}\omega_{V_l\otimes V_{l'}'}(u(b_2))\cdot \left(I_{\cB_{k,-i}}\otimes I_{0}\otimes \iota_0 \right)
    = \II_{\cB_{k,-i}}.
    \]
    \item 
    \[
    \sum_{b_2\in \Hom(V_{0},\langle v_{l-k}\rangle )}\omega_{V_l\otimes V_{l'}'}(u(b_2))\cdot \left(I_{\cAkk}\otimes I_{0}\otimes \iota_0 \right)
    = q^{\dim V_0}\indbasisgen.
    \]
    \item 
      \[
     \sum_{b_2\in \Hom(V_{0},\langle v_{l-k}\rangle )}\omega_{V_l\otimes V_{l'}'}(u(b_2))\cdot \left(I_{\widehat{\cB}_{k,i}}\otimes I_{0}\otimes \iota_0 \right)
    = 0.
    \]
\end{enumerate}
}
 \trivial[h]{
\begin{proof}
We prove (1) here and the rest are similar. Since $\cB_{k,i}$ is a single $B_{V^+_{l}}\times P'_{l'}$-orbits, by \Cref{support on a single orbit}, it is enough to compute blabla on it's representative elements. \zjl{up to here}. We have 
    \[
    \begin{split}
    &\omega_{V_l\otimes V_{l'}'}(u(b_2))\cdot \left(I_{B_{k,i}}\otimes I_{0}\otimes \iota_0 \right)=I_{B_{k,i}}\otimes \rho_{V_0\otimes V_{l'}'}( B_{k,i}\cdot u(b_2)^{-1})\cdot\left(
    I_{0}\otimes \iota_0 \right)\\
    =&I_{B_{k,i}}\otimes \psi(\langle B_{k,i}\cdot u(b_2)^{-1}, 0\rangle )\left(
    I_{0}\otimes \iota_0 \right)= \left(I_{B_{k,i}}\otimes 
    I_{0}\otimes \iota_0 \right).
    \end{split}
    \]
\end{proof}
}

\SSS
Combining all these computations and our normalization of $\cT_{\sftt_{l-k}}$ in \Cref{Normalized t_k2}, we deduce  the formula of   $\cT_{\sftt_{l}}\cdot \indbasisgen$ in \Cref{lem:TlAk} by keeping track of the scalars.






\begin{bibdiv}
\begin{biblist}
\bib{AM}{article}{
author={Adams, Jeffrey},
author={Moy, Allen},
title={Unipotent representations and reductive dual pairs over finite
fields},
journal={Trans. Amer. Math. Soc.},
volume={340},
date={1993},
number={1},
pages={309--321},
issn={0002-9947},
review={\MR{1173855}},
doi={10.2307/2154558},
}

\bib{AMR}{article}{
author = {Aubert, Anne-Marie},
author = {Michel, Jean},
author = {Rouquier, Rapha\"el},
doi = {10.1215/S0012-7094-96-08312-X},
fjournal = {Duke Mathematical Journal},
journal = {Duke Math. J.},
month = {05},
number = {2},
pages = {353--397},
publisher = {Duke University Press},
title = {Correspondance de Howe pour les groupes réductifs sur les corps finis},
url = {https://doi.org/10.1215/S0012-7094-96-08312-X},
volume = {83},
year = {1996}
}

 \bib{Be}{article}{
   title={Representations of p-adic groups, Fall 1992},
   author={Bernstein, J},
   journal={Lectures by Joseph Bernstein, Written by Karl E. Rumelhart}
 }




 


\bib{BK}{book}{
title={The admissible dual of $\GL(N)$ via compact open subgroups},
author={Bushnell, Colin John},
author = {Kutzko, Philip C},
number={129},
year={1993},
publisher={Princeton University Press},
}
 
\bib{BKt}{article}{
author = {Bushnell,CJ},
author = {Kutzko,PC},
title = {Smooth representations of reductive $p$-adic groups: structure theory via types},
fjournal = {Proceedings of the London Mathematical Society},
journal= {Proc. LMS}, 
volume = {77},
issue = {03},
month = {11},
year = {1998},
pages = {582--634},
numpages = {53},
}


\bib{C2}{book}{
author={Carter, Roger W.},
title={Finite groups of Lie type},
series={Wiley Classics Library},
note={Conjugacy classes and complex characters;
Reprint of the 1985 original;
A Wiley-Interscience Publication},
publisher={John Wiley \& Sons, Ltd., Chichester},
date={1993},
pages={xii+544},
isbn={0-471-94109-3},
review={\MR{1266626}},
}
  
\bib{Fu}{book}{
  title={Representation theory: a first course},
  author={Fulton, William},
  author={Harris, Joe},
  volume={129},
  year={2013},
  publisher={Springer Science \& Business Media}
}
 


\bib{Geck}{article}{
AUTHOR = {Geck, Meinolf},
TITLE = {A note on {H}arish-{C}handra induction},
JOURNAL = {Manuscripta Math.},
FJOURNAL = {Manuscripta Mathematica},
VOLUME = {80},
YEAR = {1993},
NUMBER = {4},
PAGES = {393--401},
ISSN = {0025-2611},
MRCLASS = {20C33 (20G05)},
MRNUMBER = {1243154},
MRREVIEWER = {Robert B. Howlett},
DOI = {10.1007/BF03026560},
URL = {https://doi-org.libproxy1.nus.edu.sg/10.1007/BF03026560},
}


\bib{GP}{book} {
   
    AUTHOR = {Geck, Meinolf},
     AUTHOR = {Pfeiffer, G\"{o}tz},
     TITLE = {Characters of finite {C}oxeter groups and {I}wahori-{H}ecke
              algebras},
    SERIES = {London Mathematical Society Monographs. New Series},
    VOLUME = {21},
 PUBLISHER = {The Clarendon Press, Oxford University Press, New York},
      YEAR = {2000},
     PAGES = {xvi+446},
      ISBN = {0-19-850250-8},
   MRCLASS = {20C15 (20C08 20F55)},
  MRNUMBER = {1778802},
MRREVIEWER = {Jian-yi Shi},
}


\bib{Ger}{article}{
title = {Weil representations associated to finite fields},
author = {G\'erardin, Paul},
journal = {Journal of Algebra},
volume = {46},
number = {1},
pages = {54-101},
year = {1977},
}

 


 

 

 
\bib{HL}{article}{
author={Howlett, R. B.},
author={Lehrer, G. I.},
title={Induced cuspidal representations and generalised Hecke rings},
journal={Inuent. Math.},
volume={58},
date={1980},
number={1},
pages={37--64},
issn={0020-9910},
review={\MR{570873}},
doi={10.1007/BF01402273},
}
 




 


\bib{Ku}{article}{  
author={Kudla, Stephen S.},
year={1986},
journal={Inuentiones Math.},
volume={83},
number={2},
title={On the local theta-correspondence},
pages={229-255},
}


\bib{L}{article}{
author={Lusztig, G.},
title={Irreducible representations of finite classical groups},
journal={Inuent. Math.},
volume={43},
date={1977},
number={2},
pages={125--175},
issn={0020-9910},
review={\MR{463275}},
doi={10.1007/BF01390002},
}


\bib{L2}{book}
{
AUTHOR = {Lusztig, George},
TITLE = {Characters of reductive groups over a finite field},
SERIES = {Annals of Mathematics Studies},
VOLUME = {107},
PUBLISHER = {Princeton University Press, Princeton, NJ},
YEAR = {1984},
PAGES = {xxi+384},
ISBN = {0-691-08350-9; 0-691-08351-7},
MRCLASS = {20G05 (14L20 20C15)},
MRNUMBER = {742472},
MRREVIEWER = {Bhama Srinivasan},
DOI = {10.1515/9781400881772},
URL = {https://doi-org.libproxy1.nus.edu.sg/10.1515/9781400881772},
}
 
 

\bib{LW}{article}{
author={Liu, Dongwen},
author={Wang, Zhicheng},
title={Remarks on the theta correspondence over finite fields},
journal={Pacific J. Math.},
volume={306},
date={2020},
number={2},
pages={587--609},
issn={0030-8730},
review={\MR{4129413}},
doi={10.2140/pjm.2020.306.587},
} 


 



\bib{MVW}{book}{
author={Moeglin, C.},
author={Vigneras, M. F.},
author={Waldspurger, J.-L.},
title={Correspondances de Howe sur un corps $p$-adique},  
series={Lecture Notes in Mathematics},
volume={1291},
year={1987},
publisher={Springer Verlag, Berlin},
}
 

 
\bib{Pan02J}{article}{
journal = {J. Math. Soc. Japan},
author = {Pan, Shu-Yen},
doi = {10.2969/jmsj/1191591993},
month = {10},
number = {4},
pages = {793--845},
publisher = {Mathematical Society of Japan},
title = {Local theta correspondence of depth zero representations and theta dichotomy},
url = {http://dx.doi.org/10.2969/jmsj/1191591993},
volume = {54},
year = {2002},
}
 

\bib{PanUni}{article}{
title={Howe Correspondence of Unipotent Characters for a Finite Symplectic/Even-orthogonal Dual Pair},
author={Pan, Shu-Yen},
year={2019},
eprint={1901.00623},
archivePrefix={arXiv},
primaryClass={math.RT}
}

 
\bib{Pan22}{article}
{title={Generalized Preservation Principle in Finite Theta Correspondence},
  author={Pan, Shu-Yen},
  journal={arXiv preprint arXiv:2207.06688},
  year={2022}
}

\bib{SZ}{article}{
author = {Sun, Binyong},
author = {Zhu, Chen-Bo},
title = {Conservation relations for local theta correspondence},
journal = {J. Amer. Math. Soc.},
fjournal = {Journal of the American Mathematical Society},
volume = {28},
year = {2015},
number = {4},
pages = {939--983},
issn = {0894-0347},
mrclass = {22E50 (11F27 22E46)},
mrnumber = {3369906},
mrreviewer = {Ivan Mati\'c},
doi = {10.1090/S0894-0347-2014-00817-1},
url = {https://doi.org/10.1090/S0894-0347-2014-00817-1},
}

\bib{Spr}{article}{
author={Springer, T. A.},
title={On the characters of certain finite groups},
conference={
title={Lie groups and their representations},
address={Proc. Summer School, Bolyai J\'{a}nos Math. Soc., Budapest},
date={1971},
},
book={
publisher={Halsted, New York},
},
date={1975},
pages={621--644},
review={\MR{0396779}},
}

\bib{SW}{book}
 {AUTHOR = {Scharlau, Winfried},
     TITLE = {Quadratic and {H}ermitian forms},
    SERIES = {Grundlehren der mathematischen Wissenschaften [Fundamental
              Principles of Mathematical Sciences]},
    VOLUME = {270},
 PUBLISHER = {Springer-Verlag, Berlin},
      YEAR = {1985},
     PAGES = {x+421},
      ISBN = {3-540-13724-6},
   MRCLASS = {11Exx (11-02 12D15 15A63 16A16 16A28)},
  MRNUMBER = {770063},
MRREVIEWER = {R. Ware},
       DOI = {10.1007/978-3-642-69971-9},
       URL = {https://doi-org.libproxy1.nus.edu.sg/10.1007/978-3-642-69971-9},
}
 

\bib{Ti}{article} {
    AUTHOR = {Tits, J.},
     TITLE = {Normalisateurs de tores. {I}. {G}roupes de {C}oxeter \'{e}tendus},
   JOURNAL = {J. Algebra},
  FJOURNAL = {Journal of Algebra},
    VOLUME = {4},
      YEAR = {1966},
     PAGES = {96--116},
      ISSN = {0021-8693},
   MRCLASS = {20.75},
  MRNUMBER = {206117},
MRREVIEWER = {Rimhak Ree},
       DOI = {10.1016/0021-8693(66)90053-6},
       URL = {https://doi-org.libproxy1.nus.edu.sg/10.1016/0021-8693(66)90053-6},
}


 
\bib{Z}{book}{
author={Zelevinsky, Andrey V.},
title={Representations of finite classical groups},
series={Lecture Notes in Mathematics},
volume={869},
note={A Hopf algebra approach},
publisher={Springer-Verlag, Berlin-New York},
date={1981},
pages={iv+184},
isbn={3-540-10824-6},
review={\MR{643482}},
}



\end{biblist}
\end{bibdiv}
 \end{document}